\documentclass{amsart}

\usepackage{amssymb}
\usepackage{amsfonts}
\usepackage{eepic}
\usepackage{color}

\allowdisplaybreaks

\font\bbold=bbold12

\newcommand{\ind}{\textrm{\bbold{1}}}

\numberwithin{equation}{section}

\newtheorem{theorem}{Theorem}[section]
\newtheorem{lemma}[theorem]{Lemma}
\newtheorem{proposition}[theorem]{Proposition}
\newtheorem{corollary}[theorem]{Corollary}
\newtheorem{remark}[theorem]{Remark}

\newtheorem{Atheorem}{Theorem}

\newtheorem{Aremark}[Atheorem]{Remark}

\newtheorem{Btheorem}{Theorem}

\newtheorem{TheoA}{Theorem A}
\newtheorem{TheoB}{Theorem B}
\newtheorem{TheoC}{Theorem C}
\newtheorem{TheoD}{Theorem D}

\newcommand{\N}{\mathbb{N}}
\newcommand{\Z}{\mathbb{Z}}
\newcommand{\R}{\mathbb{R}}
\newcommand{\C}{\mathbb{C}}

\newcommand{\summ}{\sum\nolimits}

\theoremstyle{definition}

\newcommand{\cL}{\mathcal{L}}
\newcommand{\cB}{\mathcal{B}}

\newcommand{\supp}{{\rm supp}}
\newcommand{\cM}{\mathcal{M}}
\newcommand{\cH}{\mathcal{H}}
\newcommand{\cR}{\mathcal{R}}
\newcommand{\cE}{\mathcal{E}}
\newcommand{\cZ}{\mathcal{Z}}
\newcommand{\cT}{\mathcal{T}}
\newcommand{\nphi}{\mathfrak{n}_\varphi}
\newcommand{\mphi}{\mathfrak{m}_\varphi}

\def\G{\mathrm{G}}
\def\Gd{\mathrm{G}_{\mathrm{disc}}}
\def\1{\mathbf{1}}
\def\H{\mathrm{H}}

\def\M{\mathcal{M}}

\def\RR{\mathcal{R}}

\def\V{\mathrm{\mathcal{L}G}}

\newcommand{\dem}{\noindent {\bf Proof. }}
\newcommand{\demA}{\noindent {\bf Proof of Theorem A: Unimodular case. }}
\newcommand{\demAA}{\noindent {\bf Proof of Theorem A: Nonunimodular case. }}
\newcommand{\demB}{\noindent {\bf Proof of Theorem B. }}
\newcommand{\demC}{\noindent {\bf Proof of Theorem C. }}

\newcommand{\demD}{\noindent {\bf Proof of Theorem D i) and ii). }}
\newcommand{\demDD}{\noindent {\bf Proof of Theorem D iv). }}

\newcommand{\fin}{\hspace*{\fill} $\square$ \vskip0.2cm}
\newcommand{\tens}{\otimes}

\addtocounter{tocdepth}{-1}

\begin{document}

\null

\vskip-30pt

\null

\title[Noncommutative de Leeuw theorems]{Noncommutative de Leeuw theorems}

\author[Caspers, Parcet, Perrin, Ricard]
{Martijn Caspers, Javier Parcet \\ Mathilde Perrin and \'Eric Ricard}

\maketitle

\null

\vskip-30pt

\null

\begin{abstract}
Let $\H$ be a subgroup of some locally compact group $\G$. Assume $\H$ is approximable by discrete subgroups and $\G$ admits neighborhood bases which are \lq almost-invariant\rq${}$ under conjugation by finite subsets of $\H$. Let $m: \G \to \C$ be a bounded continuous symbol giving rise to an $L_p$-bounded Fourier multiplier (not necessarily cb-bounded) on the group von Neumann algebra of $\G$ for some $1 \le p \le \infty$. Then, $m_{\mid_\H}$ yields an $L_p$-bounded Fourier multiplier on the group von Neumann algebra of $\H$ provided the modular function $\Delta_\H$ coincides with $\Delta_\G$ over $\H$. This is a noncommutative form of de Leeuw's restriction theorem for a large class of pairs $(\G,\H)$, our assumptions on $\H$ are quite natural and recover the classical result. The main difference with de Leeuw's original proof is that we replace dilations of gaussians by other approximations of the identity for which certain new estimates on almost multiplicative maps are crucial. Compactification via lattice approximation and periodization theorems are also investigated. 
\end{abstract}


\addtolength{\parskip}{+1ex}

\section*{{\bf Introduction}}

In 1965, Karel de Leeuw proved three fundamental theorems for Euclidean Fourier multipliers. Given a bounded continuous symbol $m: \R^n \to \C$, let us consider the corresponding multiplier $$\widehat{T_mf}(\xi) \, = \, m(\xi) \widehat{f}(\xi),$$ $$T_mf(x) \, = \, \int_{\R^n} m(\xi) \widehat{f}(\xi) e^{2 \pi i \langle x, \xi \rangle} \, d\xi.$$ The main results in \cite{dL} may be stated as follows:
\begin{itemize}
\item[{\bf i)}] {\bf Restriction.} If $m$ is continuous and $T_m$ is $L_p(\R^n)$-bounded $$\hskip25pt  T_{m_{\mid_{\H}}}: \int_{\H} \widehat{f}(h) \chi_h \, d\mu(h) \ \mapsto \ \int_{\H} m(h) \widehat{f}(h) \chi_h \, d\mu(h)$$ extends to a $L_p(\widehat{\H})$-bounded multiplier for any subgroup $\H \subset \R^n$ where the $\chi_h$'s stand for the characters on the dual group and $\mu$ is the Haar measure.

\vskip5pt

\item[{\bf ii)}] {\bf Periodization.} Given $\H \subset \R^n$ any closed subgroup and $m_q: \R^n / \, \H \to \C$ bounded, let $m_\pi: \R^n \to \C$  denote its $\H$-periodization which is defined by $m_\pi(\xi) = m_q(\xi + \H)$. Then we find $$\hskip25pt \big\| T_{m_\pi}: L_p(\R^n) \to L_p(\R^n) \big\| \, = \, \big\| T_{m_q}: L_p(\widehat{\R^n / \, \H}) \to L_p(\widehat{\R^n / \, \H}) \big\|.$$

\vskip5pt

\item[{\bf iii)}] {\bf Compactification.} Let $\R^n_{\mathrm{bohr}}$ be the Pontryagin dual of $\R^n_{\mathrm{disc}}$ equipped with the discrete topology. Given $m: \R^n \to \C$ bounded and continuous, the $L_p(\R^n)$-boundedness of $T_m$ is equivalent to the boundedness in $L_p(\R^n_{\mathrm{bohr}})$ of the multiplier with the same symbol $$\hskip25pt T_{m}: \summ_{\R^n_{\mathrm{disc}}} \widehat{f}(\xi) \chi_\xi \mapsto \summ_{\R^n_{\mathrm{disc}}} m(\xi) \widehat{f}(\xi) \chi_\xi.$$
\end{itemize}

Together with Cotlar's work \cite{Co}, de Leeuw theorems may be regarded as the first form of transference in harmonic analysis, prior to Calder\'on and Coifman/Weiss contributions \cite{C,CW}. The combination of the above-mentioned results produces a large family of previously unknown $L_p$-bounded Fourier multipliers ---a sample of them will appear in Appendix A--- and restriction/periodization are nowadays very well-known properties of Euclidean Fourier multipliers. Although not so much known, the compactification theorem was the core result of \cite{dL}. 

Our goal is to study these results within the context of general locally compact groups. Shortly after de Leeuw, Saeki \cite{S} extended these theorems to locally compact abelian (LCA) groups with an approach which relies more on periodization and the structure theorem of LCA groups. On the contrary, no analog transference results in the frequency group seem to exist for nonabelian groups, see \cite{Der1,DG,Herz,Weiss} for a dual approach. This gap is partly justified by the noncommutative nature of the spaces involved. Namely, the action in de Leeuw theorems occurs in the frequency groups and the Fourier multipliers must be defined in the corresponding duals. The dual of a nonabelian locally compact group can only be understood as a quantum group whose underlying space is a noncommutative von Neumann algebra. If $\mu_\G$ denotes the left Haar measure on a locally compact group $\G$ and $\lambda_\G: \G \to \mathcal{U}(L_2(\G))$ stands for the left regular representation on $\G$, the group von Neumann algebra $\V$ is the weak-$*$ closure in $\mathcal{B}(L_2(\G))$ of operators of the form $$f \, = \, \int_\G \widehat{f}(g) \lambda_\G(g) \, d\mu_\G(g) \quad \mbox{with} \quad \widehat{f} \in \mathcal{C}_c(\G).$$ The Plancherel weight is determined by $\tau_\G(f) = \widehat{f}(e)$ for $\widehat{f}$ in $\mathcal{C}_c(\G) \ast \mathcal{C}_c(\G)$ and $L_p(\widehat{\G})$ denotes the noncommutative $L_p$ space on $(\V,\tau_\G)$. Although very natural in operator algebra and noncommutative geometry, group von Neumann algebras are not yet standard spaces in harmonic analysis. The early remarkable work of Cowling/Haagerup \cite{CH,H} on approximation properties of these algebras was perhaps the first contribution in the line of harmonic analysis. The $L_p$-theory was not seriously considered until \cite{Ha}. However, only during very recent years a prolific series of results have appeared in the literature \cite{CdlS,JM,JMP1,JMP2,LdlS,NR,PRo}. 

In contrast with \cite{dL,S} where compactification and periodization took the lead respectively, we will first put the emphasis on restriction. Assume in what follows that our groups are second countable. We say that a locally compact group $\H$ is approximable by discrete subgroups ($\H \in \mathrm{ADS}$) when there exists a family of lattices $(\Gamma_j)_{j \ge 1}$ in $\H$ and associated fundamental domains $(\mathrm{X}_j)_{j \ge 1}$ which form a neighborhood basis  of the identity. On the other hand, we say that $\G$ has  small almost-invariant neighborhoods with respect to $\H$ ($\G \in [\mathrm{SAIN}]_\H$) if for every $\mathrm{F} \subset \H$ finite, there is a basis $(V_j)_{j \ge 1}$ of symmetric neighborhoods of the identity with  $$\lim_{j \to \infty} \frac{\mu_\G \big( (h^{-1} V_j h) \, \triangle \, V_j \big)}{\mu_\G(V_j)} = 0 \quad \mbox{for all} \quad h \in \mathrm{F}.$$

\begin{TheoA} \label{RestrictionThm}
Let $\H$ be a subgroup of some locally compact group $\G$. Assume $\H \in \mathrm{ADS}$ and $\G \in [\mathrm{SAIN}]_\H$. Let $m: \G \to \C$ be a bounded continuous symbol giving rise to an $L_p$-bounded multiplier for some $1 \le p \le \infty$. Then $$\big\| T_{m_{\mid_\H}}: L_p(\widehat{\H}) \to L_p(\widehat{\H}) \big\| \, \le \, \big\| T_m: L_p(\widehat{\G}) \to L_p(\widehat{\G}) \big\|$$ provided the modular function $\Delta_\H$ coincides with the restriction of $\Delta_\G$ to $\H$.
\end{TheoA}

A natural difficulty for the proof of Theorem A comes from the fact that we only assume boundedness of our multipliers. Indeed, when $\G$ is amenable, cb-bounded analogs easily follow from the recent transference results in \cite{CdlS,NR} between Fourier and Schur multipliers. It should however be noted that, for $L_p$-bounded multipliers or even for cb-bounded multipliers over nonamenable groups, our approach requires a different strategy which does not rely on previously known techniques. 

Pairs $(\G,\H)$ satisfying Theorem A include restriction onto Heisenberg groups and other classical nilpotent groups. In fact, amenable ADS subgroups of locally compact groups with $\Delta_\H = \Delta_{\G_{\mid_\H}}$ also fulfill the hypotheses. Other nonamenable pairs will be considered in the paper. Our assumptions are indeed natural for this degree of generality. The condition $\G \in [\mathrm{SAIN}]_\H$ has its roots in de Leeuw's original argument. Although not explicitly mentioned, a key point in his proof is the use of an approximate identity intertwining with the Fourier multiplier. In the Euclidean setting of \cite{dL}, this was naturally achieved by using dilations of the gaussian, which is fixed by the Fourier transform. In our general setting, the heat kernel must be replaced by other approximations and the SIN condition ---small invariant neighborhoods, which have been studied in the literature--- yields certain approximations intertwining with the Fourier multiplier. Our jump from SIN to the more flexible almost-invariant class SAIN requires a more functional analytic approach which circumvents the technicalities required for a heat kernel approach in such a general setting. The crucial novelty are certain estimates for almost multiplicative maps of independent interest. Surprisingly, our argument is equally satisfactory and much cleaner. We will prove a limiting intertwining behavior of our approximation of the identity as a consequence of the following result. 

\begin{TheoB} 
Let $(\M,\tau)$ be a semifinite von Neumann algebra equipped with a normal semifinite faithful trace. Let $T:\M\to \M$ be a subunital positive map with $\tau \circ T \leq \tau$. Then, given any $1 \le p \le \infty$ and $x \in L_{2p}^+(\M)$
$$\big\| T(x) - T(\sqrt{x})^2 \big\|_{2p} \le \frac12 \big\| T(x^2) - T(x)^2 \big\|_p^{\frac12}.$$
\end{TheoB}

We will use Haagerup's reduction method \cite{HJX} to extend the implications of Theorem B for type III von Neumann algebras. This will be the key subtle point in proving Theorem A for nonunimodular $\G$. Theorem B seems to provide new insight even in the commutative setting. Namely, arguing as in the proof of Theorem A we can use Theorem B to control the frequency support of a fractional power of a function in terms of the frequency support of the original function, up to certain small $L_p$ correction terms, we refer to Remark \ref{CommutativeApp} for further details.

Let us now go back to the other assumption in Theorem A. The ADS property of $\H$ was implicitly used in de Leeuw's original argument and could be a natural limitation for restriction of Fourier multipliers in this general setting, perhaps more powerful tools could be used for nice Lie groups \cite{VSC}. In our case, we will just prove the validity of Theorem A for discrete subgroups $\Gamma$ of a locally compact group $\G$ in the class $[\mathrm{SAIN}]_\Gamma$. Then, assuming $\H \in \mathrm{ADS}$ is approximated by $(\Gamma_j)_{j \ge 1}$, the complete statement follows from the inclusion $$[\mathrm{SAIN}]_\H \subset \bigcap_{j \ge 1} [\mathrm{SAIN}]_{\Gamma_j},$$ and the following noncommutative form of Igari's lattice approximation \cite{I,Jo}.

\begin{TheoC} \label{IgariThm}
If $\G \in \mathrm{ADS}$ is approximated by $(\Gamma_j)_{j \ge 1}$ $$\big\| T_m: L_p(\widehat{\G}) \to L_p(\widehat{\G}) \big\| \, \le \, \sup_{j \ge 1} \big\| T_{m_{\mid_{\Gamma_j}}}: L_p(\widehat{\Gamma}_j) \to L_p(\widehat{\Gamma}_j) \big\|$$ for any $1 \le p \le \infty$ and any bounded symbol $m: \G \to \C$ continuous $\mu_\G$--a.e.
\end{TheoC}

Apart from arbitrary discrete groups and many LCA groups, other nontrivial examples in the ADS class include again Heisenberg groups and other nilpotent groups. Although Theorem C is not very surprising, its proof is certainly technical and it becomes a key point in our compactification theorem. Let $\G$ be a locally compact group and write $\Gd$ for the same group equipped with the discrete topology.  Our next goal is to determine under which conditions $$\big\| T_m: L_p(\widehat{\G}) \to L_p(\widehat{\G}) \big\| \, \sim \, \big\| T_m: L_p(\widehat{\Gd}) \to L_p(\widehat{\Gd}) \big\|$$ for bounded continuous symbols. Of course, we may not expect that such an equivalence holds for arbitrary locally compact groups, since this would mean that Fourier multiplier $L_p$ theory reduces to the class of discrete group von Neumann algebras. Note also that restriction in the pair $(\G,\H)$ always holds when both group algebras admit $L_p$-compactification since restriction within the family of discrete groups follows by taking conditional expectations. This gives another evidence that compactification only holds under additional assumptions. We finally consider the periodization problem. Let $\H$ be a normal closed subgroup of some locally compact group $\G$. Consider any bounded symbol $m_q: \G / \H \to \C$ (not necessarily continuous) and construct the $\H$-periodization given by $m_\pi(g) = m_q(g\H)$. The periodization problem consists in giving conditions under which $$\big\| T_{m_\pi}: L_p(\widehat{\G}) \to L_p(\widehat{\G}) \big\|_p \, \sim \, \big\| T_{m_q}: L_p(\widehat{\G/\H}) \to L_p(\widehat{\G/\H}) \big\|_p.$$ 

\begin{TheoD} 
Let $\G$ be a locally compact unimodular group and $\H$ a normal closed subgroup of $\G$. Let us consider a bounded continuous symbol $m: \G \to \C$ and let $m_q: \G/\H \to \C$ be bounded with $\H$-periodization $m_\pi(g) = m_q(g\H)$. Then, the following inequalities hold for $1 \le p \le \infty$$\, :$
\begin{itemize}
\item[i)] If $\G$ is $\mathrm{ADS}$ $$\big\| T_m: L_p(\widehat{\G}) \to L_p(\widehat{\G}) \big\| \, \le \, \big\| T_m: L_p(\widehat{\Gd}) \to L_p(\widehat{\Gd}) \big\|.$$

\item[ii)] If $\Gd$ is amenable $$\big\| T_m: L_p(\widehat{\Gd}) \to L_p(\widehat{\Gd}) \big\| \, \le \, \big\| T_m: L_p(\widehat{\G}) \to L_p(\widehat{\G}) \big\|.$$

\item[iii)] If $\G$ is $\mathrm{LCA}$ $$\hskip3pt \big\| T_{m_\pi}: L_p(\widehat{\G}) \to L_p(\widehat{\G}) \big\| \le \big\| T_{m_q}: L_p(\widehat{\G/\H}) \to L_p(\widehat{\G/\H}) \big\|.$$

\item[iv)] If $\H$ is compact $$\hskip3pt \big\| T_{m_q}: L_p(\widehat{\G/\H}) \to L_p(\widehat{\G/\H}) \big\| \le \big\| T_{m_\pi}: L_p(\widehat{\G}) \to L_p(\widehat{\G}) \big\|.$$
\end{itemize}
\end{TheoD}

The unimodularity of $\G$ seems crucial for compactification, given the fact that $\Gd$ is always unimodular. The ADS condition is certainly natural to control Fourier multipliers on $\G$ by the same ones defined on $\Gd$. It is an interesting problem to decide whether this assumption is in fact necessary. As we will see the amenability in ii) and the commutativity in iii) (which goes back to Saeki) are very close to optimal. The inequality in iv) also holds for nonunimodular $\G$.

Our conditions above can be substantially relaxed for amenable groups in the assumption that our multipliers are completely bounded. This follows from the transference results in \cite{CdlS,NR} between Fourier and Schur multipliers ---which work in the cb-setting for amenable groups--- together with an approximation result for Schur multipliers from \cite{LdlS}. In particular, de Leeuw theorems hold in full generality in this context as we shall prove in the last section. The validity of our results for nonamenable groups is what forces us to find new arguments in this paper. We will close this article with two appendices. In Appendix A we analyze a certain family of idempotent Fourier multipliers in $\R$. By using restriction and lattice approximation we will relate these multipliers with Fefferman's theorem for the ball \cite{Fef} and solve a question from \cite{JMP1}. Appendix B contains an overview of what is known in the context of Jodeit's multiplier theorem \cite{Jo} for locally compact groups.  

\section{{\bf Almost multiplicative maps}} \label{Sect=AlmostMultiplicative}

In this section we shall prove Theorem B and some consequences  of it which will be crucial in our approach to noncommutative restrictions. Our results are of independent interest in the context of almost multiplicative maps on $L_p$. Along this section $(\M,\tau)$ will be a semifinite von Neumann algebra with a given normal semifinite faithful trace. We will need the following classical inequalities. They are well known for Schatten classes and can be found in Bhatia's book \cite[Theorems IX.4.5 and X.1.1]{Bhatia}. The proofs given there can be generalized to any semifinite von Neumann algebra, but we will provide more direct arguments. The second result is a one-sided generalization of the Powers-St{\o}rmer inequality.

\begin{lemma} \label{exchange}
Given $1 \le p \le \infty$, the identity $$\alpha^{\frac 12} \gamma \beta^{\frac 12} = \frac12 \int_\R \alpha^{-is}(\alpha \gamma+\gamma \beta)\beta^{is} \frac{ds}{\cosh (\pi s)}$$ holds in $L_p(\M)$ for any $\alpha,\, \beta,\, \gamma$ in $L_{2p}(\M)$ with $\alpha,\,\beta \ge 0$. In particular 
\begin{itemize}
\item[i)] $\big\|\alpha^{\frac 12}\gamma \beta^{\frac 12}\big\|_{p}\leq \frac 12 \big\|\alpha \gamma+\gamma \beta\big\|_p.$

\item[ii)] If $\gamma=\gamma^*$, $\big\|\alpha^{\frac 12} \gamma \alpha^{\frac 12}\big\|_p\leq\big\|\alpha \gamma\big\|_p$.
\end{itemize} 
\end{lemma}

\dem Inequalities i) and ii) follow from the first identity. By an approximation argument we may assume that $\alpha^\frac12$ and $\beta^\frac12$ have discrete spectrum, so that they are linear combinations of pairwise disjoint projections. By direct substitution this reduces the problem to $\alpha, \beta \in \R_+$ and $\gamma = 1$. Since the map $z \mapsto \alpha^{1-z}\beta^{z}$ is holomorphic on the strip $\Delta=\{0<{\rm Re}\, z<1\}$ and continuous on its closure, its value at $z = \frac12$ is given by the integral formula $$\alpha^\frac12 \beta^\frac12 = \int_{\partial\Delta} \alpha^{1-z} \beta^z \, d\mu(z)$$ where $\mu$ is the harmonic measure on $\partial \Delta$ relative to the point $z = 1/2$. Now, since this measure gives the probability for a random walk from the point $\frac12$ of hitting the boundary $\partial \Delta$, it coincides at both components $\partial_j = \{\mathrm{Re}z=j\}$ of the boundary ($j=0,1$). This means that there is a probability measure $\nu$ on $\R$ satisfying the identity $$\alpha^\frac12 \beta^\frac12 = \frac12 \Big( \int_\R \alpha^{1-is} \beta^{is} d\nu(s) + \int_\R \alpha^{-is} \beta^{1+is} \, d\nu(s)\Big) = \frac12 \int_\R \alpha^{-is}(\alpha + \beta) \beta^{is} \, d\nu(s).$$ An inspection of the harmonic measure in $\partial \Delta$ yields $d\nu(s) = ds /  \cosh (\pi s)$. Indeed, this can be obtained from the harmonic measure on the unit circle by means of a conformal map, see for instance \cite[p. 93]{BL}. The proof is complete. \fin   
 
\begin{lemma}\label{hardps} 
If $p\geq 1$, $0 < \theta \le 1$ and $x,\,y\in L_{\theta p}^+(\M)$ $$\big\|x^\theta-y^\theta\big\|_{p}\leq \big\|x-y\big\|_{\theta p}^\theta.$$
\end{lemma}

\dem Cutting $x$ and $y$ by some of their spectral projections we may assume that $(\M,\tau)$ is finite and $x,\, y \in \M$. We may also reduce the above estimate to the case $x\geq y\geq 0$. To that end, note that $$\|a-b\|_p^p\leq \|a\|_p^p+\|b\|_p^p$$ for $a,\,b \ge 0$. Indeed, let $q_+=\ind_{a-b\geq 0}$ and $q_-=\ind_{a-b< 0}$ then 
\begin{eqnarray*}
\|a-b\|_p^p & = & \| q_+(a-b)q_+ \|_p^p + \| q_-(b-a)q_- \|_p^p \\ & \le & \| q_+aq_+ \|_p^p + \| q_-bq_- \|_p^p \ \le \ \|a\|_p^p+\|b\|_p^p
\end{eqnarray*}
as $0\leq q_+(a-b)q_+\leq q_+aq_+$ and similarly for the other term. Now let $\delta_+, \delta_-$ be the positive and negative parts of $\delta = x-y = \delta_+ - \delta_-$. Let us consider the operators $$a=(x+\delta_-)^\theta-x^\theta \quad \mbox{and} \quad b=(y+\delta_+)^\theta-y^\theta.$$ Since $y + \delta_+ = x + \delta_-$ we deduce that $x^\theta - y^\theta = b-a$. Moreover, by operator monotonicity of $s\mapsto s^\theta$, $a$ and $b$ are positive. Then the result for $x + \delta_- \ge x \geq 0$ and $y + \delta_+ \ge y \geq 0$ yields $$\big\| x^\theta - y^\theta \big\|_p^p = \| a - b \|_p^p \le \|a\|_p^p + \|b\|_p^p \le \| \delta_- \|_{\theta p}^{\theta p} + \| \delta_+ \|_{\theta p}^{\theta p} =\| \delta \|_{\theta p}^{\theta p} = \big\| x-y \big\|_{\theta p}^{\theta p}.$$ Let us then prove the assertion when $x \ge y \ge 0$. We will also assume $y \geq \varepsilon \1$ to avoid unnecessary technical complications. Using the integral representation  for $s \in \M$ invertible $$s^\theta = c_\theta \int_{\R_+} \frac{t^\theta s}{s+t} \, \frac{dt}{t} \qquad \mbox{with} \qquad c_\theta = \Big( \int_{\R_+} \frac{u^\theta}{u(1+u)} \, du \Big)^{-1}.$$ Differentiating the above integral formula and putting $\delta = x-y$, we get $$x^\theta-y^\theta= c_\theta \int_0^1 \int_{\R_+} t^\theta (y+u\delta+t)^{-1} \delta (y+u\delta+t)^{-1} \, dt \, du.$$ Now, for a fixed $u \in [0,1]$, we consider the continuous function$$t \mapsto u_t = \frac{1}{\sqrt{\theta}} \, \frac{(y+u\delta)^{\frac{1-\theta}{2}}}{(y+u\delta+t)}$$ with positive values in the commutative algebra generated by $y+u \delta$. Moreover $$c_\theta \int_{\R_+} t^\theta u_t^2 \, dt = \frac{c_\theta}{\theta} \int_{\R_+} t^\theta \frac{(y +u\delta)^{1-\theta}}{(y+u\delta+t)^2} \, dt = \frac{c_\theta}{\theta} \int_{\R_+} \frac{t^\theta}{(1+t)^2} \, dt = 1.$$ Therefore, the map on $\M$ defined by $$z \mapsto c_\theta \int_{\R_+} t^\theta u_t z u_t \, dt$$ is unital, completely positive and trace preserving. In particular, it extends to a contraction on $L_p(\M)$ for all $1 \le p \le \infty$, see \cite[Remark 5.6]{HJX} for further details. We deduce 
\begin{eqnarray*}
\big\| x^\theta - y^\theta \big\|_p & \le & \theta \int_0^1 \big\| (y+u\delta\big)^{\frac {\theta-1}2}\delta (y+u\delta\big)^{\frac {\theta-1}2}\big\|_p \, du \\ & = & \theta \int_0^1 \big\|\delta^{\frac12}(y+u\delta\big)^{{\theta-1}}\delta^{\frac 12}\big\|_p \, du \ \le \ \theta \int_0^1 u ^{\theta-1} \|\delta^\theta\|_p \, du \ = \ \|\delta\|_{\theta p}^\theta,
\end{eqnarray*}
where the last inequality follows from the operator monotonicity of $s\mapsto s^{1-\theta}$. \fin

\demB Given $1 \le p \le \infty$, we claim that 
\begin{equation}\label{ThmB-claim}
\big\| R(x) - R(\sqrt{x})^2 \big\|_{2p} \le \frac12 \big\| R(x^2) - R(x)^2 \big\|_p^{\frac12}
\end{equation}
for any subunital positive map $R: \ell_\infty^n \to \M$ with values in $\M \cap L_1(\M)$ and any positive $x \in \ell_\infty^n$. The assumption above about the range of $R$ is to ensure that $R: \ell^n_\infty \to L_p(\M)$ is well-defined. Before proving this claim, let us show how this implies the assertion. Indeed, if $T: \M \to \M$ is a subunital positive map with $\tau \circ T \le \tau$ it follows from \cite[Remark 5.6]{HJX} that $T$ extends to a positive contraction on $L_p(\M)$. Now, when $x \in L_{2p}^+(\M)$ has a finite spectrum $x=\sum_{j=0}^n \lambda_j p_j$ with $\lambda_0=0$ and $p_j$ spectral projections, then $p_j \in \M \cap L_1(\M)$ for $j \ge 1$ and we may define $R: \ell_\infty^n \to \M$ by $R(e_j)=T(p_j)$, where $(e_j)_{j=1}^n$ denotes the canonical basis of $\ell_\infty^n$. The map $R$ clearly satisfies the assumptions of our claim and $R(z^\alpha) = T(x^\alpha)$ for $z = \sum_{j=1}^n  \lambda_j e_j$ and any $\alpha>0$. Hence \eqref{ThmB-claim} gives $$\big\| T(x) - T(\sqrt{x})^2 \big\|_{2p} \le \frac12 \big\| T(x^2) - T(x)^2 \big\|_p^{\frac12}$$ as desired. The general case $x \in L_{2p}^+(\M)$ follows by standard approximations. Let $$x_n \, = \, \sum_{k=1}^{n^2} \frac{k}{n} \, \ind_{[\frac{k}{n}, \frac{k+1}{n})}(x).$$ It is an exercise to show that for $\alpha\in\{1, 2, \frac12\}$, $x_n^\alpha\to x^\alpha$ in the appropriate $L_q$-space.

Let us now prove the claim \eqref{ThmB-claim}. As usual, $(e_{ij})_{i,j=1}^n$ will denote the canonical basis of the matrix algebra $\mathbb{M}_n$. We first use an explicit Stinespring's decomposition for $R$. Let $\pi: \ell_\infty^n \to \mathbb{M}_n$ be the usual diagonal inclusion and put $$u^* = \sum_{j=1}^n e_{j1} \tens R(e_j)^{\frac 12} \in \mathbb{M}_{n,1}(\M),$$ so that we have $R(x)=u\pi(x)u^*$. As $R$ is subunital $uu^* \le \1_\M$ and $u^*u \le \1_{\mathbb{M}_{n}(\M)}$. For any positive $y \in \ell_\infty^n$, we get $$R(y^2) - R(y)^2 \, = \, u \, \pi(y) \, \big( \1 - u^*u \big) \, \pi(y) \, u^* \, = \, \big| \sqrt{\1 - u^*u} \, \pi(y) \, u^* \big|^2.$$ Let us consider the following operators
\begin{eqnarray*}
a & = & R(\sqrt x) \ = \ u \, \pi(\sqrt x) \, u^* \ \in \ \M, \\ b & = & \sqrt{\1-u^*u} \, \pi(\sqrt x) \, \sqrt{\1-u^*u} \ \in \ \mathbb{M}_n(\M).
\end{eqnarray*}
Then we find $z \in \mathbb{M}_{n,1}(\M)$ with $\|z\|_\infty \le 1$ so that $\sqrt{\1-u^*u} \, \pi(\sqrt x) \, u^*= b^{\frac 12} z a^{\frac 12}$ and $$\sqrt{\1-u^*u} \, \pi(x) \, u^* = \sqrt{\1-u^*u} \, \pi(\sqrt x) \, \big( (\1-u^*u)+u^*u \big) \, \pi(\sqrt x) \, u^* = b^{\frac 32} z a^{\frac 12} + b^{\frac 12} z a^{\frac 32}.$$ We apply Lemma \ref{exchange} twice to conclude. First $$\big\| R(x)-R(\sqrt x)^2\big\|_{2p} = \big\|b^{\frac 12} z a^{\frac12} \big\|_{4p}^2 = \big\|a^{\frac 12} z^* b z a^{\frac12}\big\|_{2p}\leq\big\|a z^* b z \big\|_{2p} \le \big\|a z^* b \big\|_{2p}.$$ Then, taking $(\alpha, \gamma, \beta) = (a, a^{\frac 12} z^*b^{\frac12}, b)$ we obtain $$\big\|a z^* b\big\|_{2p}\leq \frac 12 \big\|a^{ \frac 32} z^* b^{\frac 12}+a^{ \frac 12} z^* b^{\frac 32} \big\|_{2p}=\frac 12\big\| R(x^2)-R(x)^2\big\|_p^{\frac 12}.$$ This completes the proof of our claim and also the proof of Theorem B. \fin

\begin{corollary} \label{Tx-x} 
Let $T: \M \to \M$ be a subunital positive map with $\tau\circ T\leq \tau$. Then there exists a universal constant $C > 0$ such that the following inequality holds for any $x \in L_2^+(\M)$ and any $0 < \theta \le 1$ $$\big\| T(x^\theta) - x^\theta \big\|_{\frac{2}{\theta}} \, \le \, C \, \big\| T(x) - x \big\|_2^{\frac{\theta}{2}} \, \|x\|_2^{\frac{\theta}{2}}.$$
\end{corollary}

\dem Given $x \in L_2^+(\M)$, note that $$\|T(x) - x\|_2^2 \le \|Tx\|_2^2 + \|x\|_2^2 \le \tau (T(x^2)) + \|x\|_2^2 \le 2 \, \|x\|_2^2$$ by Kadison's inequality for $T$ and $\tau \circ T \le \tau$. In particular, the result is trivially true for $\theta = 1$ with constant $2^{\frac14}$. We claim the assertion holds for $\theta = 2^{-n}$ with constant $\frac32$. We will proceed by induction since we know it holds for $n=0$. By Lemma \ref{hardps} and $n+1$ applications of Theorem B
\begin{eqnarray*}
\lefteqn{\hskip-20pt \big\|T(x^{2^{-(n+1)}})-x^{2^{-(n+1)}}\big\|_{2^{n+2}}^2} \\ [6pt] & \le & \big\|T(x^{2^{-(n+1)}})^2-x^{2^{-n}}\big\|_{2^{n+1}} \\ [6pt] & \le &   \big\| T(x^{2^{-(n+1)}})^2-T(x^{2^{-n}}) \big\|_{2^{n+1}} +\big\|T(x^{2^{-n}}) - x^{2^{-n}}\big\|_{2^{n+1}}\\ & \le & \underbrace{\prod_{j=0}^n 2^{-2^{-j}}}_{C_n = 2^{2^{-n}-2}} \big\|T(x^2) - T(x)^2 \big\|_1^{2^{-(n+1)}} +\big\|T(x^{2^{-n}}) - x^{2^{-n}}\big\|_{2^{n+1}}.
\end{eqnarray*}
On the other hand, Kadison's inequality and $\tau \circ T \le \tau$ yield 
\begin{eqnarray*}
\big\|T(x^2)-T(x)^2\big\|_1 & = & \tau \big( T(x^2)-T(x)^2 \big) \ \le \ \tau \big( x^2-T(x)^2 \big) \\ & \le & \big\|x^2-T(x)^2\big\|_1 \ \le \ 2 \big\| T(x) - x \big\|_2 \|x\|_2
\end{eqnarray*}
since $T$ extends to a contraction on $L_2(\M)$ by \cite[Remark 5.6]{HJX}. Combining the above estimates with the induction hypothesis for $\theta=2^{-n}$ we finally deduce that $$\big\|T(x^{2^{-(n+1)}})-x^{2^{-(n+1)}}\big\|_{2^{n+2}}^2 \, \le \, \Big[ \frac{3}{2} + 2^{2^{-(n+1)}} C_n \Big] \, \|x\|_2^{2^{-(n+1)}} \big\|T(x) - x \big\|_2^{2^{-(n+1)}}.$$ However, the constant in the right hand side is less than $9/4$ and the result follows for $\theta = 2^{-(n+1)}$ which completes the induction argument. For other values of $0 < \theta < 1$ we write $\theta = \alpha 2^{-(n+1)}$ for some $\alpha \in (1,2)$. Recall that $s \mapsto s^\alpha$ is operator convex and $s \mapsto s^{\frac \alpha2}$ is operator concave on $\R_+$, so that $$T \big(x^{2^{-(n+1)}}\big)^\alpha \, \le \, T \big(x^\theta \big) \, \le \, T\big(x^{2^{-n}}\big)^{\frac{\alpha}{2}}.$$ In conjunction with Lemma \ref{hardps}, Theorem B and our result for $\theta = 2^{-n}$, this yields
\begin{eqnarray*}
\big\| T(x^\theta) - x^\theta \big\|_{\frac{2}{\theta}} & \le & \big\| T(x^\theta) - T(x^{2^{-n}})^{\frac{\alpha}{2}} \big\|_{\frac{2}{\theta}} + \big\| T(x^{2^{-n}})^{\frac{\alpha}{2}} - x^\theta \big\|_{\frac{2}{\theta}} \\ [3pt] & \le & \big\| T(x^{2^{-n}})^{\frac{\alpha}{2}} \hskip1pt - \hskip1pt T(x^{2^{-(n+1)}})^{\alpha} \big\|_{\frac{2}{\theta}} \hskip1pt + \big\| T(x^{2^{-n}})^{\frac{\alpha}{2}} \hskip1pt - \hskip1pt x^\theta \big\|_{\frac{2}{\theta}} \\ [2pt] & \le & \big\| T(x^{2^{-n}}) - T(x^{2^{-(n+1)}})^2 \big\|_{2^{n+1}}^{\frac{\alpha}{2}} + \big\| T(x^{2^{-n}}) - x^{2^{-n}} \big\|_{2^{n+1}}^{\frac{\alpha}{2}} \\ & \le & \Big[ \Big( 2^{2^{-(n+1)}} C_n \Big)^{\frac{\alpha}{2}} + \Big( \frac32 \Big)^{\frac{\alpha}{2}} \Big] \big\| Tx-x \big\|_2^{\frac{\theta}{2}} \|x\|_2^{\frac{\theta}{2}}.
\end{eqnarray*}
Hence, a simple calculation shows that the result follows for some $C \le \frac{3+\sqrt{2}}{2}$. \fin

\begin{corollary} \label{Ty-y} 
Let $T: \M \to \M$ be a subunital positive map with $\tau\circ T\leq \tau$. Then there exists a universal constant $C > 0$ such that the following inequality holds for any self-adjoint $y \in L_2(\M)$ with polar decomposition $y = u|y|$ and any $0 < \theta \le 1$ $$\big\| T(u|y|^\theta) - u|y|^\theta \big\|_{\frac{2}{\theta}} \, \le \, C \, \big\| T(y) - y \big\|_2^{\frac{\theta}{4}} \, \|y\|_2^{\frac{3\theta}{4}}.$$
\end{corollary}

\dem Let us write $y = y_+ - y_-$ for the decomposition of $y$ into its positive and negative parts, so that $u|y|^\theta = y_+^\theta - y_-^\theta$. By positivity of the trace and $T$, we have $$\tau \big(T(y_+)y_+\big) + \tau \big(T(y_-)y_-\big) \, \ge \, \tau \big(T(y)y\big).$$ In particular $$\big\| T(y_+)-y_+\big\|_2^2 + \big\| T(y_-)-y_-\big\|_2^2 \, \le \, 2 \|y\|_2^2 - 2 \tau \big(T(y)y\big) \, \le \, 2 \big\| T(y)-y \big\|_2 \|y\|_2.$$ Using this and Corollary \ref{Tx-x} we deduce 
\begin{eqnarray*}
\lefteqn{\hskip-20pt \big\| T(u|y|^\theta) - u|y|^\theta \big\|_{\frac{2}{\theta}}^{\frac{4}{\theta}}} \\ & \le & \Big[ \big\| T(y_+^\theta) - y_+^\theta \big\|_{\frac{2}{\theta}} + \big\| T(y_-^\theta) - y_-^\theta \big\|_{\frac{2}{\theta}} \Big]^{\frac{4}{\theta}} \\ [2pt] & \le & 2^{\frac{4}{\theta}-1}  \Big[ \big\| T(y_+^\theta) - y_+^\theta \big\|_{\frac{2}{\theta}}^{\frac{4}{\theta}} + \big\| T(y_-^\theta) - y_-^\theta \big\|_{\frac{2}{\theta}}^{\frac{4}{\theta}} \Big] \\ & \le & \frac{(2C)^{\frac{4}{\theta}}}{2} \Big[ \big\| T(y_+) - y_+ \big\|_2^2 \|y_+\|_2^2 + \big\| T(y_-) - y_- \big\|_2^2 \|y_-\|_2^2 \Big] \\ & \le & \frac{(2C)^{\frac{4}{\theta}}}{2} \Big[ \big\| T(y_+) - y_+ \big\|_2^2 + \big\| T(y_-) - y_- \big\|_2^2 \Big] \|y\|_2^2 \\ [2pt] & \le & \ (2C)^{\frac{4}{\theta}}  \big\| T(y) - y \big\|_2 \|y\|_2^3.
\end{eqnarray*}
The assertion follows taking powers $\theta/4$ at both sides of the above estimate. \fin 

\begin{remark} \label{CommutativeApp}
\emph{The above corollary will be useful to localize the frequency support of square roots for elements in $L_p(\widehat{\G})$. This is even interesting in the commutative case where we may control the frequency support of a fractional power $f^\theta$ in terms of the frequency support of $f$, up to small $L_p$ corrections. As an illustration, if the Fourier transform of $f \in L_2(\R)$ is supported by $(-\alpha, \alpha)$, we may consider the positive definite functions $$\zeta_\beta(x) \, = \, \Big( 1 - \frac{|x|}{2\beta} \Big)_+ \quad \mbox{for} \quad \beta > 0.$$ The associated Fourier multipliers are positive, unital and trace preserving, so that we are in position to apply our results above. When $\beta = \alpha / 2\varepsilon$, we obtain that $\mathrm{supp} \, \zeta_\beta \subset \frac{1}{\varepsilon} (- \alpha, \alpha)$ and $1- \zeta_\beta(x) \le \varepsilon$ for $|x| \le \alpha$. This yields for $p \ge 2$ $$\Big\| \Big( \zeta_\beta (f^{\frac{2}{p}})^\wedge - (f^{\frac{2}{p}})^\wedge \Big)^{\vee} \Big\|_p^p \, \le \, C \big\| \zeta_\beta \widehat{f} - \widehat{f} \hskip2pt \big\|_2 \|f\|_2 \, \lesssim \, \varepsilon \|f\|_2^2.$$}
\end{remark}

\begin{remark}
\emph{It is well known that if $T:\M\to \M$ satisfies the above hypothesis and $\tau$ is finite, then its fixed points form a $*$-subalgebra. This is not true anymore when $\tau$ is semifinite, take for instance the map $x\mapsto s^*xs$ on $\mathcal{B}(\ell_2)$ where $s$ is a one-sided shift. Nevertheless, in general, it is not difficult to show using the generalized singular value decomposition that if $x \in L_1^+(\M) \cup L_2^+(\M)$ satisfies $T(x) = x$, then $T(x^\theta) = x^\theta$ for $\theta \in [0,1]$. Hence one could think of an ultraproduct argument to get perturbation results as given explicitly in Corollary \ref{Tx-x}, with an upper bound of the form  $F(\|T(x)-x\|_2)$ for certain continuous function $F$ vanishing at $0$. Unfortunately, semifiniteness is not preserved by ultraproduct and one would have to deal with type III von Neumann algebras. The situation is then much more intricate (even to define $T$ on $L_p(\M)$), that is why we choose to deduce the type III result from the semifinite one in Section \ref{Sect=MultiplicativeTypeIII}. The fact that there exists a unital completely positive map $T:(\M,\varphi)\to (\M,\varphi)$ with $\varphi\circ T= \varphi$ but $T$ does not commute with the modular group of $\varphi$ (think of a right multiplier on a quantum group with its left Haar measure) is an evidence that in the type III situation one need extra arguments as those of Corollary \ref{Tx-x arbitrary vNA}.}
\end{remark}

\section{{\bf Group algebras}} \label{GroupvNaSect}

Let $\G$ be a locally compact group equipped with its left Haar measure $\mu_\G$. Let $\lambda_\G: \G \to \cB(L_2(\G))$ be the left regular representation $\lambda_\G(g)(\xi)(h)=\xi(g^{-1}h)$ for any $\xi \in L_2(\G)$. When no confusion can arise, we shall write $\mu,\lambda$ for the left Haar measure and the left regular representation of $\G$. Recall the definition of the convolution in $\G$ $$\xi * \eta (g) \, = \, \int_\G \xi(h) \eta(h^{-1}g) \, d\mu(h).$$ We say that $\xi \in L_2(\G)$ is left bounded if the map $\eta \in \mathcal{C}_c(\G) \mapsto \xi \ast \eta \in L_2(\G)$ extends to a bounded operator on $L_2(\G)$, denoted by $\lambda(\xi)$. This operator defines the Fourier transform of $\xi$. The weak operator closure of the linear span of $\lambda(\G)$ defines the group von Neumann algebra $\cL \G$. It can also be described as the weak closure in $\mathcal{B}(L_2(\G))$ of operators of the form $$f \, = \, \int_\G \widehat{f}(g) \lambda(g) \, d\mu(g) \, = \, \lambda(\widehat{f} \, ) \quad \mbox{with} \quad \widehat{f} \in \mathcal{C}_c(\G).$$ The Plancherel weight $\tau_\G:\cL \G _+  \to [0,\infty]$ is determined by the identity $$\tau_\G(f^*f) = \int_\G | \widehat{f}(g) |^2 \, d\mu(g)$$ when $f = \lambda(\widehat{f} \, )$ for some left bounded $\widehat{f} \in L_2(\G)$ and $\tau_\G(f^*f) = \infty$ for any other $f \in \cL \G$. Again, we shall just write $\tau$ for $\tau_\G$ when the underlying group is clear from the context. After breaking into positive parts, this extends to a weight on a weak-$\ast$ dense domain within the algebra $\V$. It will be instrumental to observe that the standard identity $$\tau(f) \, = \, \widehat{f}(e)$$ applies for $f = \lambda(\widehat{f} \, ) \in \lambda(\mathcal{C}_c(\G)\ast \mathcal{C}_c(\G))$, see \cite[Section 7.2]{Ped} and \cite[Section VII.3]{TakII} for a detailed construction of the Plancherel weight. Note that for $\G$ discrete, $\tau$ coincides with the natural finite trace given by $\tau(f) = \langle f \delta_e, \delta_e \rangle$. It is clear that the Plancherel weight is tracial if and only if $\G$ is unimodular, which will be the case until Section \ref{non unimodular}. In the unimodular case, $(\cL \G, \tau)$ is a semifinite von Neumann algebra and we may construct the noncommutative $L_p$-spaces $$L_p(\cL \G, \tau) \, = \, L_p(\widehat{\G}) \, = \,\left\{\begin{array}{ll} \overline{\lambda(\mathcal{C}_c(\G)\ast \mathcal{C}_c(\G))}^{\| \ \|_p}&\quad \mbox{for } 1\leq p <2 \\ \overline{\lambda(\mathcal{C}_c(\G))}^{\| \ \|_p}&\quad \mbox{for } 2\leq p <\infty\end{array}\right.,$$ where the norm is given by $$\|f\|_p=\tau(|f|^p)^{1/p}$$ and the $p$-th power is calculated by functional calculus applied to the (possibly unbounded) operator $f$, we refer to Pisier/Xu's survey \cite{PX} for more details on noncommutative $L_p$-spaces. On the other hand, since left bounded functions are dense in $L_2(\G)$, the map $\lambda: \xi \mapsto \lambda(\xi)$ extends to an isometry from $L_2(\G)$ to $L_2(\widehat{\G})$. We will refer to it as the Plancherel isometry and use it repeatedly in the sequel with no further reference. Given a symbol $m:\G \to \C$, we may consider the associated multiplier $T_m$ defined by $$T_m(f) \, = \, \int_\G m(g) \widehat{f}(g) \lambda(g) \, d\mu(g) \quad \mbox{for} \quad \widehat{f} \in \mathcal{C}_c(\G) \ast \mathcal{C}_c(\G).$$ $T_m$ is called an $L_p$-Fourier multiplier if it extends to a bounded map on $L_p(\widehat{\G})$. 

The rest of this section will be devoted to collect some elementary results around amenability and Fell absorption principles that will be used in the sequel. We will also need the following result, which we prove for completeness.

\begin{lemma} \label{LHsubalgLG}
Let $\G$ be a second countable locally compact unimodular group. Then the group von Neumann algebra $\cL \H$ is a von Neumann subalgebra of $\cL \G$ for any closed unimodular subgroup $\mathrm{H}$ of $\G$.
\end{lemma}

\dem By the Effros-Mackey cross section theorem \cite[Theorem 5.4.2]{Sri}, there exists a Borel measurable map $\sigma: \mathrm{H} \! \setminus \!\! \G \to \G$ defined on the space of right cosets of $\G$. Hence, we have a Borel measurable correspondence between $\G$ and $\mathrm{H} \! \setminus \!\! \G \times \H$ given by $$\Upsilon : \G \ni g \mapsto (\H g,h(g)) \in \mathrm{H} \! \setminus \!\! \G \times \mathrm{H},$$ where $g = h(g) \sigma(\H g)$. According to \cite[Theorem 2.49]{F} for right cosets and since both $\G$ and $\H$ are unimodular, we know that there exists a $\G$-invariant Radon measure on right cosets. Therefore, the map $$\xi \mapsto \xi \circ \Upsilon^{-1}$$ defines an isometry between $L_2(\G)$ and $L_2(\mathrm{H} \! \setminus \!\! \G \times \H)$. This allows us to identify the algebras $\cB(L_2(\G))$ and $\cB(L_2(\mathrm{H} \! \setminus \!\! \G) \otimes_2 L_2(\H))$. Then, for any $h\in \H$ we get the identity
\begin{eqnarray*}
\big( id \otimes \lambda_\H(h) \big) (\xi \circ \Upsilon^{-1}) (\H g,h(g)) & = & \xi(h^{-1}h(g)\sigma(\H g)) \\ & = & \xi(h^{-1}g) \ = \ \lambda_\G(h)(\xi)(g)
\end{eqnarray*}
for $\xi \in L_2(\G)$ and $g \in \G$, which proves that $\cL \H \simeq \{\lambda_\G(h) \; : \; h \in \H\}'' \subset \cL \G$. \fin

In the sequel, if no confusion is possible and when Lemma \ref{LHsubalgLG} applies, we might use the notation $\lambda(h)$ to denote both $\lambda_\G(h)$ and $\lambda_\H(h)$. Let us now recall some well-known characterizations of amenability. Recall that amenability is stable under closed subgroups, quotients, direct products and group extensions. After that, we also give a formulation of Fell absorption principle in $L_p$ from \cite{PaPi}.

\begin{lemma} \label{amenable}
\emph{TFAE} for any locally compact group $\G$\emph{:}
\begin{itemize}
\item[i)] $\G$ is amenable

\vskip3pt

\item[ii)] \emph{F\o{}lner condition}. Given $\varepsilon > 0$ and $\mathrm{F} \subset \G$ finite, there exists $U_{\mathrm{F},\varepsilon} \subset \G$ of finite positive measure such that $\mu( U_{\mathrm{F}, \varepsilon }g \: \triangle \: U_{\mathrm{F}, \varepsilon } ) < \varepsilon \mu(U_{\mathrm{F},\varepsilon})$ for all $g \in \mathrm{F}$.

\vskip3pt

\item[iii)] \emph{Almost invariant vectors}. Given $\varepsilon > 0$ and $\mathrm{F} \subset \G$ finite, there exists a norm one function $\xi \in L_2(\G)$ such that $\|\lambda(g)\xi-\xi\|_{L_2(\G)}<\varepsilon$ for all $g\in \mathrm{F}$. 

\vskip3pt

\item[iv)] The inequality $$\Big\|\sum_{g\in \mathrm{F}} a_g \Big\|_{\mathcal{M}}  \, \le \, \Big\| \sum_{g\in \mathrm{F}} a_g \otimes \lambda(g)\Big\|_{\mathcal{M}\overline{\otimes}\cL \G}$$ holds for any finite $\mathrm{F} \subset \G$, any von Neumann algebra $\cM$ and $(a_g)_{g\in \mathrm{F}}\subset \mathcal{M}$.
\end{itemize}
\end{lemma}

\begin{lemma} \label{Fell}
Given a discrete group $\G$, we have\emph{:}
\begin{itemize}
\item[i)] If $\pi:\G\to \mathcal{U}(\mathcal{H})$ is strongly continuous, then $$\lambda \otimes \pi \simeq \lambda \otimes 1_{\mathcal{H}}$$ are unitarily equivalent with $1_{\mathcal{H}}$ the trivial representation on $\mathcal{H}$. 

\vskip2pt

\item[ii)] Let $\pi:\G\to \mathcal{U}(\mathcal{H})$ be strongly continuous and assume $\mathcal{N}=\pi(\G)''$ is finite. Then, given $1 \le p \le \infty$, any semifinite von Neumann algebra $\cM$ and any $a:\G\to L_p(\cM)$ continuous and compactly supported we have
\begin{eqnarray*}
\lefteqn{\hskip-70pt \Big\| \int_{\G} a(g) \otimes \lambda(g) \otimes \pi(g) \, d \mu (g)\Big\|_{L_p(\cM \overline{\otimes} \cL \G \overline{\otimes} \mathcal{N})}} \\ \null \hskip100pt \null & = & \Big\| \int_{\G} a(g) \otimes \lambda(g) \, d \mu (g) \Big\|_{L_p(\cM \overline{\otimes} \cL \G )}.
\end{eqnarray*}
\end{itemize}
\end{lemma}

\section{{\bf Lattice approximation}}\label{Sect:Igari}

In this section, we want to deduce the boundedness of an $L_p$-Fourier multiplier from the uniform boundedness of its restriction to certain families of lattices. As stated in Theorem C, this will be possible if $\G$ is approximated by these lattices in the sense $\G \in \mathrm{ADS}$ defined in the Introduction. Observe that if $\G \in \mathrm{ADS}$ is approximated by $(\Gamma_j)_{j\geq 1}$, then the union $\cup_j \Gamma_j$ of the approximating lattices is dense in $\G$. Indeed, let $g\in \G$ and $V$ be an open neighborhood of $g$. Then $Vg^{-1}$ is an open neighborhood of $e$ and for $j$ large enough we have $\mathrm{X}_j\subset Vg^{-1}$. Let $g_j$ be the representant of $g^{-1}$ in $\mathrm{X}_j$. In other words, there exists $\gamma_j \in \Gamma_j$ such that $g_j=\gamma_j g^{-1}$. This implies $\gamma_j=g_j g \in \mathrm{X}_j g \subset V$, so that $\Gamma_j  \cap V \neq \emptyset$ and we deduce the density result mentioned above. In the proof of Theorem C we shall need a couple of auxiliary results, which are stated below. 

\begin{lemma}\label{lattice-unimodular}
If $\G$ admits a lattice $\Gamma$ with $\Delta_{\G_{\mid_\Gamma}} = \Delta_\Gamma$, then $\G$ is unimodular. 
\end{lemma}

The result above implies that every ADS group is by definition unimodular. In particular, our preliminaries on group von Neumann algebras from Section \ref{GroupvNaSect} suffice for Theorem C. We need one more elementary result. 

\begin{lemma}\label{K-Omega}
Let $\G$ be a locally compact group and $\mathrm{K} \subset \Omega \subset \G$ with $\mathrm{K}$ compact and $\Omega$ open. Let $(V_j)_{j\geq 1}$ be a basis of neighborhoods of the identity. Then, there exists an index $j_0\geq 1$ such that for any $j \ge j_0$ $$\mathrm{K} \subset \bigcup_{g \in \mathrm{K}}gV_{j}\subset \Omega.$$
\end{lemma}

\dem Take a left invariant distance $d$ on $\G$, so that $d(K,\Omega^c)=\delta>0$. Since $\mathrm{diam}(V_j)\to 0$, any $j_0$ with $\mathrm{diam}(V_j) < \delta$ for $j\geq j_0$ satisfies the conclusion. \fin

\vskip2pt 

\demC The case $p=2$ is straightforward since $m$ is continuous almost everywhere and the union of lattices $\Gamma_j$ is dense in $\G$, so that the $L_\infty$-norm of $m$ is determined by lattice approximation. On the other hand, by a standard duality argument, we may assume that $p < 2$. Moreover, the case $p=1$ follows from the assertion for $1 < p < 2$ and the three lines lemma $$\| T_m \|_{1 \to 1} \, \le \, \lim_{p \to 1} \| T_m \|_{p \to p} \, \le \, \lim_{p \to 1} \sup_{j \ge 1} \| T_{m_{\mid_{\Gamma_j}}} \|_{p \to p} \, \le \, \sup_{j \ge 1} \| T_{m_{\mid_{\Gamma_j}}} \|_{1 \to 1}.$$ Therefore, we may and will assume in what follows that $1 < p < 2$. The strategy will be to approximate $T_m f$ weakly in $L_p$ by a sequence $S_j f$ constructed from a family $(S_j)_{j\geq 1}$ of uniformly bounded maps as follows. For each $j\geq 1$ we first define the map $$\Phi_j : \cL\Gamma_j \ni \lambda(\gamma) \mapsto h_j^*\lambda(\gamma)h_j \in \V,$$ where $h_j=\lambda (\ind_{\mathrm{X}_j}) \in \cL \G$. Since $\G$ is locally compact and $(\mathrm{X}_j)_{j\geq 1}$ is a basis of neighborhoods of $e$, we may assume that $\mathrm{X}_j$ lies in a compact set. In particular we have $0<\mu(\mathrm{X}_j)<\infty$. It is clear that $\Phi_j$ is completely positive and we may define the family of operators $$\Phi_j^p \, = \, \mu(\mathrm{X}_j)^{-2+\frac{1}{p}}\Phi_j.$$ Now we note the straightforward inequality $$\|\Phi_j^{\infty}(\mathbf{1})\|_{\cL \G} \, = \, \mu(\mathrm{X}_j)^{-2}\|h_j\|_{\cL \G}^2 \leq \mu(\mathrm{X}_j)^{-2}\|\ind_{\mathrm{X}_j}\|_{L_1(\G)}^2 \, = \, 1.$$ Moreover, since the sets $(\gamma \mathrm{X}_j)_{\gamma \in \Gamma_j}$ are disjoint, we also have for $\gamma \in \Gamma_j$ 
\begin{eqnarray*}
\tau\big(\Phi_j(\lambda(\gamma))\big) & = & \tau \big( h_j^*\lambda(\gamma)h_j \big) \ = \ \big\langle \lambda(\gamma)h_j, h_j \big\rangle_{L_2(\widehat{\G})} \\ & = & \big\langle \ind_{\gamma \mathrm{X}_j},\ind_{\mathrm{X}_j} \big\rangle_{L_2(\G)} \ = \ \mu(\mathrm{X}_j)\delta_{\gamma,e} \ = \ \mu(\mathrm{X}_j)\tau(\lambda(\gamma)).
\end{eqnarray*}
By complete positivity of $\Phi_j$, the first estimate implies that $\Phi_j^\infty: \mathcal{L} \Gamma_j \to \V$ is a contractive map. The second estimate implies that $\Phi_j^1$ is trace preserving and hence defines a contraction $L_1(\mathcal{L}{\Gamma}_j) \to L_1(\V)$ by means of \cite[Remark 5.6]{HJX}. Using interpolation of analytic families of operators, we get $$\big\| \Phi_j^p: L_p(\widehat{\Gamma}_j) \to L_p(\widehat{\G}) \big\| \, \le \, 1 \quad \mbox{for} \quad 1 \le p \le \infty.$$ On the other hand, the $L_2$-adjoints $\Psi_j = \Phi_j^*$ are given by $$\Psi_j(f) \, = \, \sum_{\gamma\in \Gamma_j}\tau\big(h_j^*\lambda(\gamma^{-1})h_j f\big) \lambda(\gamma)$$ for $f\in \cL\G$. Moreover, given $1\leq p \leq \infty$, consider the contractions $$\Psi_j^p = (\Phi^{p'}_j)^* = \mu(\mathrm{X}_j)^{-1-\frac{1}{p}}\Psi_j: L_p(\widehat{\G})\to L_p(\widehat{\Gamma}_j),$$ where $p'$ denotes the conjugate of $p$. We are finally ready to introduce the maps $$S_j = \Phi_j^p T_{m_{\mid_{\Gamma_j}}} \Psi_j^p = \mu(\mathrm{X}_j)^{-3} \Phi_j T_{m_{\mid_{\Gamma_j}}} \Psi_j: L_p(\widehat{\G}) \to L_p(\widehat{\G}),$$ which are uniformly bounded by $$C_p:=\sup_{j\geq 1}\big\| T_{m_{\mid_{\Gamma_j}}}: L_p(\widehat{\Gamma}_j) \to L_p(\widehat{\Gamma}_j) \big\|.$$ If we fix $f \in \lambda(\mathcal{C}_c(\G)\ast \mathcal{C}_c(\G))$, the sequence $(S_j f)_{j\geq 1}$ is uniformly bounded in $L_p(\widehat{\G})$ by $C_p\|f\|_p$ and it accumulates in the weak topology. We claim that $S_j f$ weakly converges to $T_mf = w\mbox{-}L_p\mbox{-}\lim_j S_j f$. The theorem will follow by the $L_p$-density of $\lambda(\mathcal{C}_c(\G)\ast \mathcal{C}_c(\G))$. To prove it, we can reduce ourselves to show 
\begin{equation}\label{claim-approxL2}
T_m f =L_2\mbox{-}\lim_{j \to \infty} S_j f \quad \mbox{for any} \quad f \in \lambda(\mathcal{C}_c(\G)\ast \mathcal{C}_c(\G)). 
\end{equation}
Indeed, if it holds true and $q$ is any $\tau$-finite projection $$\lim_{j \to \infty} \big\| qT_m f-qS_jf \big\|_p \, \le \, \|q\|_r \lim_{j \to \infty} \|T_m f-S_jf\|_2 \, = \, 0,$$ where $\frac{1}{p}=\frac{1}{r}+\frac{1}{2}$. Hence $$qT_m f \, = \, L_p\mbox{-}\lim_{j \to \infty} (qS_j f) \, = \, w \mbox{-}L_p\mbox{-}\lim_{j \to \infty} (qS_j f) \, = \, q \big( w\mbox{-}L_p\mbox{-}\lim_{j \to \infty} S_j f \big)$$ for any $\tau$-finite projection $q$, which implies $T_mf = w\mbox{-}L_p\mbox{-}\lim_j S_j f$. We now turn to the proof of the key result \eqref{claim-approxL2}. Let us introduce some notations. For $j\geq 1$ define $$L_j: L_2(\widehat{\G}) \ni f \mapsto \mu(\mathrm{X}_j)^{-1} h_j f \in L_2(\widehat{\G})$$ and $$P_j: L_2(\widehat{\G}) \ni f \mapsto \frac{1}{\mu(\mathrm{X}_j)} \sum_{\gamma \in \Gamma_j} \big\langle f,\lambda(\gamma) h_j \big\rangle_{L_2(\widehat{\G})} \lambda(\gamma) h_j \in L_2(\widehat{\G}).$$ Given $g \in \G$ and since $(\gamma \mathrm{X}_j)_{\gamma \in \Gamma_j}$ forms a partition of $\G$, there exists a unique $\gamma \in \Gamma_j$ such that $g \in \gamma \mathrm{X}_j$. Let us write $\gamma_j(g)$ for this element and consider the map $m_j: \G \to \mathbb{C}$ given by $m_j(g) = m(\gamma_j(g))$. We claim that 
\begin{itemize}
\item[i)] $S_j=L_j^*P_jT_{m_j}L_j$ on $L_2(\widehat{\G})$.

\vskip2pt

\item[ii)] $L_j, L_j^*, P_j: L_2(\widehat{\G})\to L_2(\widehat{\G})$ are contractive and $$\big\| T_{m_j}: L_2(\widehat{\G}) \to L_2(\widehat{\G}) \big\| \, \le \, \|m\|_\infty.$$

\vskip2pt

\item[iii)] Given $f \in \lambda(\mathcal{C}_c(\G))$, the following identity holds $$\hskip20pt \lim_{j \to \infty} \big\| L_jf - f \big\|_2 + \big\| L_j^*f - f \big\|_2 + \big\| P_jf - f \big\|_2 + \big\| T_{m_j}f - T_m f \big\|_2 \, = \, 0.$$ In fact, the first three summands also converge to $0$ for $f \in L_2(\widehat{\G})$.

\end{itemize}
The $L_2$-convergence \eqref{claim-approxL2} follows from this. Indeed, i) gives for $f \in \lambda(\mathcal{C}_c(\G))$
\begin{eqnarray*}
\big\| T_m f -S_j f \big\|_2 & \le & \big\| T_m f - L_j^*T_m f \big\|_2 \\ & + & \big\| L_j^*T_m f - L_j^*P_jT_m f \big\|_2 \\ & + & \big\| L_j^*P_jT_m f - L_j^*P_jT_{m_j} f \big\|_2 \\ & + & \big\| L_j^*P_jT_{m_j} f - L_j^*P_jT_{m_j}L_j f \big\|_2,
\end{eqnarray*}
which clearly tends to $0$ as $j \to \infty$ by ii) and iii). Therefore, it suffices to justify the properties i), ii) and iii). Let us start by noticing the following identity which follows from Plancherel's isometry $$\big\langle T_{m_j}f, \lambda(\gamma) h_j \big\rangle_{L_2(\widehat{\G})} \, = \, \big\langle m_j \widehat{f}, \ind_{\gamma \mathrm{X}_j} \big\rangle_{L_2(\G)} \, = \, m(\gamma) \big\langle f, \lambda(\gamma) h_j \big\rangle_{L_2(\widehat{\G})}.$$
Applying this to $h_jf$ we get
\begin{eqnarray*}
S_jf & = & \mu(\mathrm{X}_j)^{-3} \sum_{\gamma \in \Gamma_j}m(\gamma) \big\langle f,h_j^*\lambda(\gamma)h_j \big\rangle_{L_2(\widehat{\G})}h^*_j\lambda(\gamma)h_j \\ & = & \mu(\mathrm{X}_j)^{-3} \sum_{\gamma \in \Gamma_j} \big\langle T_{m_j}(h_jf),\lambda(\gamma) h_j \big\rangle_{L_2(\widehat{\G})} h_j^*\lambda(\gamma)h_j \ = \ L_j^*P_jT_{m_j}L_j f,
\end{eqnarray*}
which proves i). Claim ii) for $L_j$ follows from $\|\mu(\mathrm{X}_j)^{-1}h_j\|_\infty \le \|\mu(\mathrm{X}_j)^{-1} \ind_{\mathrm{X}_j}\|_1 \le 1$. The boundedness for $P_j$ is clear since it is the orthogonal projection onto the closed linear span of $(\lambda(\gamma) h_j)_{\gamma \in \Gamma_j}$. The last assertion in ii) is trivial since $\|m_j\|_\infty \le \|m\|_\infty$. Let us finally prove the convergence results in property iii). By \cite[Proposition 2.42]{F}, the family $\widehat{h}_j = \mu(\mathrm{X}_j)^{-1}\ind_{\mathrm{X}_j}$ forms an approximation of the identity, so that $$\lim_{j \to \infty} \big\| \widehat{h}_j \ast \xi - \xi \big\|_{L_2(\G)} \, = \, 0 \quad \mbox{for} \quad \xi \in L_2(\G).$$ By Plancherel's isometry, this yields vanishing limits for the first two terms in iii). Moreover, the third term will converge to $0$ if and only if the orthogonal projection $\widetilde{P}_j$ of $L_2(\G)$ onto $\mathrm{span}\{\ind_{\gamma \mathrm{X}_j} : \gamma \in \Gamma_j\}$ satisfies 
\begin{equation} \label{(iii)'}
\lim_{j \to \infty} \big\| \widetilde{P}_j \xi - \xi \big\|_2 \, = \, 0 \quad \mbox{for any} \quad \xi \in L_2(\G).
\end{equation} 
By the density of the simple functions in $L_2(\G)$, we may assume that $\xi=\ind_{\Omega}$ for a Borel subset $\Omega$ of $\G$ with finite measure. Moreover, since the Haar measure is outer regular, $\Omega$ can be assumed to be open. On the other hand, given any $\varepsilon > 0$ and since $\mu$ is inner regular on open sets, there exists a compact set $\mathrm{K} \subset \Omega$ such that $\mu(\Omega \setminus \mathrm{K}) \le \varepsilon$. By applying Lemma \ref{K-Omega} to the basis of neighborhoods of the identity given by $(\mathrm{X}_j^{-1}\mathrm{X}_j)_{j\geq 1}$, we obtain that there exists $j_0\geq 1$ such that for any $j\geq j_0$ $$\mathrm{K} \subset \bigcup_{g \in \mathrm{K}} \gamma_{j}(g)\mathrm{X}_{j} \subset \bigcup_{g \in \mathrm{K}} g \mathrm{X}_{j}^{-1}\mathrm{X}_{j} \subset \Omega.$$ The sets $(\gamma \mathrm{X}_{j})_{\gamma \in \Gamma_{j}}$ being disjoint, we can find a subset $\mathrm{F} \subset \mathrm{K}$ satisfying $$\mathrm{K} \subset \bigsqcup_{g \in \mathrm{F}} \gamma_{j}(g)\mathrm{X}_{j} \subset \Omega.$$ Moreover, since the sets $\Omega$ and $\gamma_{j}(g)\mathrm{X}_{j}$ are of finite measure, the set $\mathrm{F}$ has to be finite. Hence, the function $\eta=\sum_{g\in \mathrm{F}}\ind_{\gamma_{j}(g) \mathrm{X}_{j}}$ satisfies $\|\xi-\eta\|_2^2 \le \mu(\Omega\setminus \mathrm{K}) \le \varepsilon$ and the limit \eqref{(iii)'} is proved. It remains to consider the last term in iii). Let $\varepsilon >0$ and $f \in \lambda(\mathcal{C}_c(\G))$ be frequency supported by a compact set $\mathrm{K}$. Since $\gamma_j(g)\to g$ as $j\to \infty$ and $m$ is continuous $\mu$--a.e., we have $m_j \to m$ $\mu$-a.e. Moreover, by Egoroff's theorem \cite[Theorem 2.33]{F2}, there exists a set $\mathrm{E}\subset \mathrm{K}$ with $\mu(\mathrm{E}) < \varepsilon$ and such that $m_j \to m$ uniformly on $\mathrm{K} \setminus \mathrm{E}$. Pick $j_0 \ge 1$ satisfying $$\sup_{g \in \mathrm{K}\setminus \mathrm{E}}|m_j(g)-m(g)| \leq \varepsilon^{\frac12}$$ for all $j \geq j_0$. Then we get for $j \ge j_0$
\begin{eqnarray*}
\lefteqn{\hskip-15pt \big\| T_{m_j}f-T_mf \big\|^2_{L_2(\widehat{\G})} \ = \ \big\| (m_j-m) \widehat{f} \, \big\|^2_{L_2(\G)}} \\ & \le & \int_{\mathrm{K} \setminus \mathrm{E}} |\widehat{f}(g)|^2 \big| m_j(g)-m(g) \big|^2 d\mu(g) + \int_\mathrm{E} |\widehat{f}(g)|^2 \big| m_j(g)-m(g) \big|^2 d\mu(g) 
\end{eqnarray*}
which is dominated by $\varepsilon \big( \|\widehat{f}\|_2^2+4 \|m\|_{\infty}^2 \|\widehat{f}\|_\infty^2 \big)$ and proves iii) for the last term. \fin

\begin{remark}\label{Rk-Igari}
\emph{Modifying the proof above, we may extend Theorem C. Namely, if $\G\in \mathrm{ADS}$ is approximated by $(\Gamma_j)_{j \ge 1}$ and both $m: \G \to \C$ and $\widetilde{m}_j:\G\to \C$ are a.e. continuous symbols such that $\widetilde{m}_j\to m$ uniformly, then 
$$\big\| T_m: L_p(\widehat{\G}) \to L_p(\widehat{\G}) \big\| \, \le \, \sup_{j \ge 1} \big\| T_{(\widetilde{m}_j)_{\mid_{\Gamma_j}}}: L_p(\widehat{\Gamma}_j) \to L_p(\widehat{\Gamma}_j) \big\|$$ for any $1 \le p \le \infty$. Indeed, it suffices to define $$S_j=\mu(\mathrm{X}_j)^{-3} \Phi_j T_{(\widetilde{m}_j)_{\mid_{\Gamma_j}}} \Psi_j$$ and consider $m_j(g)=\widetilde{m}_j(\gamma_j(g))$ in the proof of Theorem C. Then the fourth summand in iii) will follow by means of Plancherel's isometry noticing that we have $|m_j(g)-m(g)| \le \|\widetilde{m}_j-m\|_\infty + |m(\gamma_j(g))-m(g)|$ and hence converges to $0$ a.e. Then we conclude as in the proof of Theorem C.}  
\end{remark}

\begin{remark} \label{examples-ADS}
\emph{We have not performed an extensive study of groups satisfying the ADS condition. Apart from discrete groups and many LCA groups, of particular interest to us is the Heisenberg group, defined as the set $\H_n = \R^n \times \R^n \times \R$ with inner law $(a,b,c) \cdot (a',b',c') = ( a+a', b+b', c+c' + \frac12 ( \langle a, b' \rangle - \langle a', b \rangle ) )$. It is a simple example of ADS group. Namely, take for instance the family of lattices $\Gamma_j = \frac{1}{j} \Z^n \times \frac{1}{j} \Z^n \times \frac{1}{2j^2} \Z$ which trivially satisfy the ADS condition. Other nilpotent groups satisfying the ADS condition are the groups $\H(\mathbb{K},n)$ of upper triangular matrices over the field $\mathbb{K}= \R, \C$ with $1$'s on the diagonal. In this case, a simple choice of lattices is $\Gamma_j = \mathrm{Id}_n + \langle j^{r-s} \Z \otimes e_{r,s} : 1 \le r < s \le n \rangle$.}
\end{remark}

\section{{\bf The restriction theorem}} \label{Sect=Restriction}

In this section we prove Theorem A for unimodular groups. In other words, we prove that under the SAIN condition, $L_p$-Fourier multipliers of a unimodular group $\G$ restrict to multipliers of any ADS subgroup $\mathrm{H}$, and this restriction map is norm decreasing. All our work so far will be needed in the proof. 

\demA Let us first reduce the proof to the particular case of discrete subgroups. Indeed, let $\mathrm{H} \in \mathrm{ADS}$ approximable by the family $(\Gamma_j)_{j \ge 1}$ and assume that $\G \in [\mathrm{SAIN}]_{\mathrm{H}}$. Since $\Gamma_j \subset \mathrm{H} \subset \G$ and $$[\mathrm{SAIN}]_\H \subset \bigcap_{j \ge 1} [\mathrm{SAIN}]_{\Gamma_j},$$ the pairs $(\G, \Gamma_j)_{j \ge 1}$ are under the hypotheses of Theorem A for discrete subgroups. Moreover, since $\H$ is ADS, it follows from Lemma \ref{lattice-unimodular} that $\H$ must be unimodular. Therefore, Theorem A for discrete subgroups in conjunction with Theorem C yields 
\begin{eqnarray*}
\lefteqn{\hskip-30pt \big\| T_{m_{\mid_\H}}: L_p(\widehat{\H}) \to L_p(\widehat{\H}) \big\|} \\ \hskip20pt & \le & \sup_{j \ge 1} \big\| T_{m_{\mid_{\Gamma_j}}}: L_p(\widehat{\Gamma}_j) \to L_p(\widehat{\Gamma}_j) \big\| \ \le \ \big\| T_m: L_p(\widehat{\G}) \to L_p(\widehat{\G}) \big\|.
\end{eqnarray*}
Hence, we shall consider in what follows a discrete subgroup $\Gamma$ of a locally compact unimodular group $\G$ satisfying $\G \in [\mathrm{SAIN}]_\Gamma$. Observe that by unimodularity of $\G$ and discreteness of $\Gamma$ our assumption ${\Delta_\G}_{\mid_\Gamma} = \Delta_\Gamma$ is superfluous. We have reduced the proof of Theorem A for $\G$ unimodular to this particular case. Arguing as we did at the beginning of the proof of Theorem C and using the continuity of the symbol for $p=2$, it suffices to consider $2<p<\infty$. Moreover, by density of the trigonometric polynomials in $L_p(\widehat{\Gamma})$, it is enough to prove that
\begin{equation}\label{aim-ThmA}
\| T_{m_{\mid_\Gamma}} f\|_{L_p(\widehat{\Gamma})} \, \le \, \big\| T_m: L_p(\widehat{\G}) \to L_p(\widehat{\G}) \big\| \, \|f\|_{L_p(\widehat{\Gamma})}
\end{equation}
for any trigonometric polynomial $f \in \cL \Gamma$. As we explained in the Introduction, the basic idea is to construct an approximation of the identity in $\G$ which intertwines with the pair $(T_m, T_{m_{\mid_\Gamma}})$ in the limit. Let us fix such a trigonometric polynomial $f_0 \in \mathcal{L} \Gamma$ and let $\mathrm{F} \subset \Gamma$ denote its frequency support $$f_0 \, = \, \sum_{\gamma \in \mathrm{F}} \widehat{f}_0(\gamma) \lambda(\gamma).$$ Let $(V_j)_{j\geq 1}$ be the symmetric neighborhood basis of the identity associated to $\mathrm{F}$ by the SAIN condition. Moreover, since $\Gamma$ is discrete, we may take $j$ large enough and assume that the sets $(\gamma V_j)_{\gamma \in \Gamma}$ are disjoint. Let us define the selfadjoint elements $h_j=\mu(V_j)^{-1/2}\lambda(\ind_{V_j})$ with polar decomposition $h_j=u_j|h_j|$, and set $$\Phi_j^q: \lambda(\gamma) \mapsto \lambda(\gamma)u_j|h_j|^{\frac{2}{q}} \quad \mbox{for} \quad \gamma \in \Gamma \ \mbox{ and } \ 2\leq q \le  \infty.$$ Then the proof will rely on the two following results. 

\noindent {\bf Claim A.} Given $2\leq q \leq \infty$, we have:
\begin{itemize}
\item[i)] $\Phi_j^q$ extends to a contraction $L_q(\widehat{\Gamma}) \rightarrow L_q(\widehat{\G})$. 

\item[ii)] Given any $f \in\cL \Gamma$ frequency supported by $\mathrm{F}$, we have $$\lim_{j \to \infty} \|\Phi_j^q f \|_{L_q(\widehat{\G})}=\|f\|_{L_q(\widehat{\Gamma})}.$$ 
\end{itemize}

\noindent {\bf Claim B.} Given $2\leq q <p$ and any trigonometric polynomial $f$ in $\cL \Gamma$ $$\lim_{j\to \infty} \| \Phi_{j}^q( T_{m_{\mid_\Gamma}} f) - T_m(\Phi_{j}^q f) \|_{L_q(\widehat{\G})} = 0.$$

Let us finish the proof of Theorem A before proving these two claims. Let $f_0$ be the trigonometric polynomial in $\cL \Gamma$ frequency supported by $\mathrm{F}$ that we have fixed above. The algebra $\cL\Gamma$ being finite, we have 
$$\| T_{m_{\mid_\Gamma}} f_0 \|_{L_p(\widehat{\Gamma})}=  \lim_{q \nearrow p} \| T_{m_{\mid_\Gamma}} f_0 \|_{L_q(\widehat{\Gamma})}.$$ Since $T_{m_{\mid_\Gamma}} f_0$ is also frequency supported by $\mathrm{F}$, Claims A and B yield
\begin{eqnarray*}
\Vert T_{m_{\mid_\Gamma}} f_0 \Vert_{L_p(\widehat{\Gamma})} & = & \lim_{q \nearrow p } \lim_{j \to \infty} \Vert \Phi_{j}^q (T_{m_{\mid_\Gamma}} f_0) \Vert_{L_q(\widehat{\G})} \\ & = &   \lim_{q \nearrow p } \lim_{j \to \infty} \Vert  T_{m }( \Phi_{j}^q  f_0) \Vert_{L_q(\widehat{\G})} \\ & \le & \lim_{q \nearrow p } \lim_{j \to \infty}  \big\| T_m: L_q(\widehat{\G}) \to L_q(\widehat{\G}) \big\|   \Vert  \Phi_{j}^q  f_0  \Vert_{L_q(\widehat{\G})} \\ & = &  \lim_{q \nearrow p } \| T_m \|_{q \to q} \Vert f_0 \Vert_{L_q(\widehat{\Gamma})} \le \big\| T_m:  L_p(\widehat{\G}) \to L_p(\widehat{\G}) \big\| \Vert f_0 \Vert_{L_p(\widehat{\Gamma})}.
\end{eqnarray*}
The latter inequality follows by interpolation since $$\|T_m\|_{q \to q} \, \le \, \|T_m\|_{2 \to 2}^{1-\theta} \|T_m\|_{p \to p}^\theta$$ for $\frac{1}{q} = \frac{1-\theta}{2} + \frac{\theta}{p}$, so that $\theta \to 1$ as $q \to p$. This completes the proof of \eqref{aim-ThmA}.

\vskip2pt

\noindent \textbf{Proof of Claim A.} Since the SAIN condition implies $\G$ second countable, we may consider $\cL \Gamma$ as a von Neumann subalgebra of $\cL \G$ by Lemma \ref{LHsubalgLG}. Thus Claim A i) is clear for $q=\infty$ by writing $\Vert \Phi_j^\infty f \Vert_{\cL \G} = \Vert f u_j \Vert_{\cL \G} \leq \Vert f \Vert_{\cL \G} = \Vert f \Vert_{\cL \Gamma}.$ Moreover, by Plancherel's isometry and disjointness of the sets $(\gamma V_j)_{\gamma \in \Gamma}$ we get
\begin{equation}\label{Phi2-isometry}
\| \Phi_j^2 f \|_{L_2(\widehat{\G})}^2 = \| f h_j \|_{L_2(\widehat{\G})}^2 =  \mu(V_j)^{-1}\Big\| \sum_{\gamma \in \Gamma} \widehat{f}(\gamma) \ind_{\gamma V_j} \Big\|_{L_2(\G)}^2  = \|f\|_{L_2(\widehat{\Gamma})}^2.
\end{equation}
Claim A i) then follows using interpolation for analytic families of operators, we leave the details to the reader. The upper estimate in Claim A ii) follows from i) and it suffices to show that $\lim_j \|\Phi_j^q(f)\|_q \ge \|f\|_q$ for trigonometric polynomials $f$ in $\cL \Gamma$ frequency supported by $\mathrm{F}$. Let $q^\ast$ be the $L_2$-conjugate index of $q$, so that $1/q + 1/q^\ast =1/2$. We have $$\Vert f \Vert_{L_q(\widehat{\Gamma})} \, = \, \Vert f^\ast \Vert_{L_q(\widehat{\Gamma})} \, = \, \sup_{\begin{subarray}{c} \|k\|_{L_{q^*}(\widehat{\Gamma})} \le 1 \\ k \ \mathrm{trigonometric \ polynomial} \end{subarray}} \Vert k f^\ast \Vert_{L_2(\widehat{\Gamma})}.$$ Fix such a polynomial $k = \sum_{\gamma \in \mathrm{M}} \widehat{k}(\gamma) \lambda(\gamma)$. Then, since $\Phi_j^2$ is an isometry by \eqref{Phi2-isometry}
\begin{eqnarray*}
\Vert k f^\ast \Vert_{L_2(\widehat{\Gamma})} & = & \Vert \Phi_j^2( k f^\ast ) \Vert_{L_2(\widehat{\G})} \ = \ \Vert k f^\ast h_j \Vert_{L_2(\widehat{\G} )} \\ & \leq & \Vert k f^\ast h_j - \Phi_j^{q^\ast}(k) u_j \Phi_j^q(f)^\ast \Vert_{L_2(\widehat{\G})} + \Vert  \Phi_j^{q^\ast}(k) u_j \Phi_j^q(f)^\ast \Vert_{L_2(\widehat{\G})}.
\end{eqnarray*}
By H\"older's inequality and Claim A i) $$\Vert  \Phi_j^{q^\ast}(k) u_j \Phi_j^q(f)^\ast \Vert_{L_2(\widehat{\G})} \leq  \Vert \Phi_j^{q^\ast}(k) \Vert_{L_{q^\ast}(\widehat{\G})} \Vert \Phi_j^q(f)^\ast \Vert_{L_q(\widehat{\G})} \leq  \Vert \Phi_j^q (f)\Vert_{L_q(\widehat{\G})}.$$ For the first summand, let us prove that $$\lim_{j\to \infty}\Vert k f^\ast h_j - \Phi_j^{q^\ast}(k) u_j \Phi_j^q(f)^\ast \Vert_{L_2(\widehat{\G})} = 0.$$ This will complete the proof of Claim A. Since $h_j$ is self-adjoint $$ \Phi_j^{q^\ast}(k) u_j \Phi_j^{q}(f)^\ast = k u_j \vert h_j \vert^{2/q^\ast} u_j \vert h_j \vert^{2/q} u_j^\ast f^\ast = k h_j f^\ast.$$ Then, by Plancherel's isometry and the Cauchy-Schwarz inequality we get 
\begin{eqnarray*}
\lefteqn{\hskip-20pt \Vert k f^\ast h_j - k h_j f^\ast \Vert_{L_2(\widehat{\G})}} \\ & = & \Big\Vert \sum_{\gamma' \in \mathrm{M}, \gamma \in \mathrm{F}} \widehat{k}(\gamma') \overline{\widehat{f}(\gamma)} \big(\lambda(\gamma') \lambda(\gamma^{-1}) h_j - \lambda(\gamma') h_j \lambda(\gamma^{-1}) \big) \Big\Vert_{L_2(\widehat{\G})} \\ & \le & \sum_{\gamma'\in \mathrm{M}, \gamma \in \mathrm{F}} \big| \widehat{k}(\gamma') \, \widehat{f}(\gamma) \big| \, \big\Vert \lambda(\gamma') \lambda(\gamma^{-1}) h_j - \lambda(\gamma') h_j \lambda(\gamma^{-1}) \big\Vert_{L_2(\widehat{\G})} \\ & = &  \sum_{\gamma'\in \mathrm{M}, \gamma \in \mathrm{F}} \big| \widehat{k}(\gamma') \, \widehat{f}(\gamma) \big| \, \mu(V_j)^{-1/2} \big\Vert \, \ind_{\gamma^{-1}V_j \gamma} - \ind_{V_j} \big\Vert_{L_2(\G)} \\ & = &  \sum_{\gamma'\in \mathrm{M}, \gamma \in \mathrm{F}} \big| \widehat{k}(\gamma') \, \widehat{f}(\gamma) \big| \Big( \frac{\mu(\gamma^{-1} V_j \gamma \triangle V_j) }{\mu(V_j)} \Big)^{\frac{1}{2}} \\ & \le & \Big( \sum_{\gamma'\in \mathrm{M}, \gamma \in \mathrm{F}} \vert \widehat{k}(\gamma') \widehat{f}(\gamma) \vert^2 \Big)^{\frac{1}{2}} \Big( \sum_{\gamma'\in \mathrm{M}, \gamma \in \mathrm{F}} \frac{\mu(\gamma^{-1} V_j \gamma \triangle V_j) }{\mu(V_j)} \Big)^{\frac{1}{2}} \\ & = & \Vert k \Vert_{L_2(\widehat{\G})} \Vert f \Vert_{L_2(\widehat{\G})} \vert \mathrm{M} \vert^{\frac{1}{2}} \Big( \sum_{ \gamma \in \mathrm{F}} \frac{\mu(\gamma^{-1} V_j \gamma \triangle V_j) }{\mu(V_j)} \Big)^{\frac{1}{2}},
\end{eqnarray*}
which converges to 0 as $j \to \infty$ since we assumed have that $\G \in [\mathrm{SAIN}]_\Gamma$.

\vskip2pt

\noindent \textbf{Proof of Claim B.} Without loss of generality we may assume that $f = \lambda(\gamma)$ for some $\gamma \in \Gamma$ by the triangle inequality in $L_q(\widehat{\G})$. Replacing $m$ by $m(\gamma\; \cdot)$ we may assume that $\gamma = e$ (we leave the details to the reader here). This means that we are reduced to prove 
\begin{equation}\label{Eqn=IntertwiningProof}
\lim_{j \to \infty} \big\Vert m(e) u_j \vert h_j \vert^{\frac{2}{q}} - T_m( u_j \vert h_j \vert^{\frac{2}{q}} ) \big\Vert_{L_q(\widehat{\G})} = 0.
\end{equation}
Given $\varepsilon > 0$ and since $m$ is continuous in $e \in \G$, there exists a neighborhood  $ U_\varepsilon$ of the identity such that $\vert m(g) - m(e) \vert < \varepsilon$ for every $g \in U_\varepsilon$. Since $\G$ is locally compact we may assume that $U_\varepsilon$ is relatively compact and so $\mu(U_\varepsilon) < \infty$. Let $W_\varepsilon$ be a symmetric neighborhood of $e$ with $W_\varepsilon^2 \subset U_\varepsilon$ and define $$\zeta (g) = \frac{\mu( W_\varepsilon \cap g W_\varepsilon) }{\mu(W_\varepsilon)} = \frac{\langle \lambda(g) \ind_{W_\varepsilon}, \ind_{W_\varepsilon}\rangle }{\mu(W_\varepsilon)}.$$ Hence, $\zeta $ is a coefficient function of the left regular representation and the coefficient is given by the positive vector state with respect to the vector $\mu(W_\varepsilon)^{-\frac{1}{2}}\ind_{W_\varepsilon} $. It is then standard that $\zeta $ is continuous, positive definite and $\zeta (e) = 1$. Furthermore by construction $\supp \: \zeta  \subset U_\varepsilon$. Let $T_\zeta $ be the associated Fourier multiplier, then $T_\zeta : \cL \G \rightarrow \cL \G$ is a normal, trace preserving, unital, completely positive map. This implies that it extends to a contraction $$T_\zeta : L_p(\widehat{\G}) \rightarrow L_p(\widehat{\G})$$ for every $1 \leq p \leq \infty$. By Plancherel isometry we have $$\big\| T_\zeta  h_j - h_j \big\|_{L_2(\widehat{\G})}^2 =  \big\| (\zeta -1) \mu(V_j)^{-\frac12}\ind_{V_j} \big\|_{L_2(\G)}^2 =  \frac{1}{\mu(V_j)} \int_{V_j} \vert \zeta (g) - 1 \vert^2 d\mu(g),$$ which converges to $0$ as $j \to \infty$ since $V_j \rightarrow \{ e\}$ and $\zeta $ is continuous at $e$. At this point we need our result on almost multiplicative maps. Indeed, since $h_j$ is a self-adjoint operator of $L_2$-norm one, we deduce from Corollary \ref{Ty-y} that
\begin{equation}\label{Eqn=LimitThatCostsALotOfWork}
\lim_{j \to \infty} \big\| T_\zeta ( u_j \vert h_j \vert^{\frac2q} ) - u_j \vert h_j \vert^{\frac2q} \big\|_{L_q(\widehat{\G})} \, = \, 0.
\end{equation}
Let us now prove \eqref{Eqn=IntertwiningProof}. Setting $z_j = u_j \vert h_j \vert^{2/q}$ we write 
\begin{eqnarray*}
\Vert m(e) z_j - T_mz_j \Vert_{L_q(\widehat{\G})} & \le & \Vert m(e) (z_j - T_\zeta z_j) \Vert_{L_q(\widehat{\G})} \\ & + & \Vert m(e) T_\zeta z_j - T_m(T_\zeta z_j) \Vert_{L_q(\widehat{\G})} \\ & + & \Vert T_m(T_\zeta z_j)  -  T_mz_j \Vert_{L_q(\widehat{\G})} \ = \ A_j+B_j+C_j.
\end{eqnarray*}
By \eqref{Eqn=LimitThatCostsALotOfWork}, $\lim_j A_j=\lim_j C_j=0$. By definition of $U_\varepsilon$ we have $$\big\| T_{(m(e) - m)\zeta } :L_2(\widehat{\G})\to L_2(\widehat{\G})\big\| =  \Vert (m(e) - m) \zeta  \Vert_{L_\infty(\G)} < \varepsilon.$$ On the other hand, since $\|T_\zeta: L_p(\widehat{\G})\to L_p(\widehat{\G})\|=1$ we get $$\big\| T_{(m(e) - m)\zeta } :L_p(\widehat{\G})\to L_p(\widehat{\G}) \big\| \le \vert m(e) \vert + \big\| T_{m } :L_p(\widehat{\G})\to L_p(\widehat{\G}) \big\|.$$ Applying the three lines lemma to the symbol $(m(e) - m)\zeta$ we obtain $$B_j  \, \le \, \varepsilon^{1-\theta} \Big(\vert m(e) \vert + \big\| T_{m } :L_p(\widehat{\G})\to L_p(\widehat{\G}) \big\| \Big)^\theta \quad \mbox{for} \quad \frac{1}{q} = \frac{1-\theta}{2} + \frac{\theta}{p}.$$ This implies \eqref{Eqn=IntertwiningProof}, which gives Claim B and completes the proof of Theorem A. \fin

We end this section by giving some examples of groups satisfying the conditions of Theorem A. We have already considered the ADS condition in the previous section, so let us analyze the SAIN condition. There are two general conditions which imply small almost invariant neighborhoods:
\begin{itemize}
\item $\G \in[\mathrm{SIN}]_\H$ (small invariant neighborhoods) if there exists a neighborhood basis of the identity of $\G$ consisting of open sets that are invariant under conjugation with respect to $\H$, see for instance \cite{GM, L} for this class of pairs $(\G, \mathrm{H})$ when $\G = \mathrm{H}$. Of course, we have $[\mathrm{SIN}]_\H \subset [\mathrm{SAIN}]_\H$.

\vskip3pt 

\item Another interesting class of pairs satisfying the SAIN condition is given by amenable discrete subgroups $\Gamma$ satisfying $\Delta_{\G_{\mid_\Gamma}} = \Delta_\Gamma$, see Theorem \ref{AmenableImpliesSAIN} below. As a consequence of it, we shall show that Theorem A holds for pairs $(\G, \mathrm{H})$ with $\mathrm{H}$ any ADS amenable group.
\end{itemize} 

\noindent Concrete examples (even in the nonunimodular setting) will be given in Section \ref{non unimodular}. 

\begin{remark}
\emph{Both properties above are strictly weaker than the SAIN condition since none of them is included in the other one. To see this, let us construct examples of pairs $(\G,\Gamma)$, where $\Gamma$ is a discrete subgroup of a unimodular, locally compact group $\G$, satisfying only one of these two properties:} 
\begin{itemize}
\item[\textbf{i)}] \emph{The free group with two generators $\mathbb{F}_2$ can be represented as a (non-closed) subgroup of $SO(3)$. This way $\mathbb{F}_2$ acts on $\mathbb{R}^3$ and the open balls $B_r(0) \subset \mathbb{R}^3$ with center 0 and radius $r$ are invariant under the action of $\mathbb{F}_2$. We may consider the semidirect product $\G = \mathbb{R}^3 \rtimes \mathbb{F}_2$, which is unimodular since the action of $\mathbb{F}_2$ is measure preserving. Then the sets $B_r(0)$ are naturally contained in $\G$ and in fact form a basis of neighborhood of the identity which are invariant under conjugation with respect to $\mathbb{F}_2$. Hence $\mathbb{R}^3 \rtimes \mathbb{F}_2 \in [{\rm SIN}]_{\mathbb{F}_2}$ but $\mathbb{F}_2$ is not amenable.}

\vskip4pt 

\item[\textbf{ii)}] \emph{Let $\G$ be the Heisenberg group in $\mathbb{R}^n$ and $\Gamma = \mathbb{Z}^n\times \{0\}\times \{0\} \subset \G$. Then $\Gamma$ satisfies our second property above but $\G \notin [{\rm SIN}]_\Gamma$. Indeed, let $U$ be a small neighborhood of $(0,0,0)$ invariant under conjugation by $\Gamma$. Assume that $U\subset \mathbb{R}^n \times \mathbb{R}^n \times [-L,L]$ for some $L>0$. Since conjugation in the Heisenberg group gives $$(-a,-b,-c) \cdot (x,y,t) \cdot (a,b,c) = \big( x,y, t - \langle a, y \rangle + \langle b, x \rangle \big),$$ we deduce that $(x,y,t) \in U \Rightarrow (x,y,t-ay)\in U$ for all $a\in \Z^n$. But we can find an element $(x,y,t)\in U$ with $y\neq 0$ and a sequence $(a_j)_{j\geq 1}$ in $\Z^n$ verifying $|a_j y|\to \infty$ which contradicts this property.}
\end{itemize}
\end{remark}

\begin{remark}
\emph{We already know from Lemma \ref{lattice-unimodular} that every ADS group must be unimodular. On the other hand, it also holds that $\G \in [\mathrm{SAIN}]_\H$ with $\Delta_\H = \Delta_{\G_{\mid_\H}}$ implies $\mathrm{H}$ unimodular since 
\begin{eqnarray*}
\Delta_\H (h) & = & \Delta_\G(h) \ = \ \lim_{j \to \infty} \frac{\mu(h^{-1} V_j h)}{\mu(V_j)} \\ & = & \lim_{j \to \infty} \frac{\mu(h^{-1} V_j h) - \mu(h^{-1} V_j h \setminus V_j)}{\mu(V_j)} \\ & = & \lim_{j \to \infty} \frac{\mu(h^{-1} V_j h \cap V_j)}{\mu(V_j)} \ = \ 1 - \lim_{j \to \infty} \frac{\mu(V_j \setminus h^{-1} V_j h)}{\mu(V_j)} \ = \ 1
\end{eqnarray*}
for every $h \in \H$. In particular, all our conditions in Theorem A point to the unimodularity of $\H$. As we shall see in Section \ref{cb}, this is not the case when we work with amenable groups in the category of operator spaces. We leave as an open problem to decide whether unimodularity is an essential assumption for restriction of Fourier multipliers.}  
\end{remark}

\section{{\bf The compactification theorem}}

We now extend de Leeuw's compactification theorem. In other words, given a locally compact group $\G$, let us write $\Gd$ to denote the same group equipped with the discrete topology. Under the conditions in Theorem D, we prove that the $L_p$-boundedness of a Fourier multiplier on $\G$ is equivalent to the $L_p$-boundedness of that multiplier defined on $\Gd$. In this section and for the sake of clarity, we will write $\lambda=\lambda_\G$ and $\lambda'=\lambda_{\Gd}$ for the left regular representation on $\G$ and $\Gd$ respectively. Moreover, we shall use a similar terminology for trigonometric polynomials in both $\V$ and $\cL \Gd$ $$f \, = \, \sum_{g \in \mathrm{F}} \widehat{f}(g) \lambda(g) \ \Leftrightarrow \ f' \, = \, \sum_{g \in \mathrm{F}} \widehat{f}(g)\lambda'(g).$$ Before proving the compactification theorem, let us first discuss the conditions on the group $\G$ that we impose. In de Leeuw's proof of the compactification theorem for $\R^n$, the following basic properties were crucial:
\begin{itemize}
\item[\textbf{P1)}] We have $$\R^n \, = \, \overline{\bigcup_{j \ge 1} 2^{-j} \Z^n}.$$ 

\vskip4pt

\item[\textbf{P2)}] There is an injective homomorphism $\Psi: \R^n \to \R^n_{\mathrm{bohr}}$ ---the dual to the canonical inclusion map $\R^n_{\mathrm{disc}} \to \R^n$--- with dense image and such that $f=f' \circ \Psi$ for any pair $(f, f') \in L_\infty(\R^n) \times L_\infty(\R_{\mathrm{bohr}}^n)$ of trigonometric polynomials with matching Fourier coefficients. In particular $$\hskip20pt \|f'\|_{L_\infty(\R^n_{\mathrm{bohr}})} = \sup_{\xi \in \R^n} |f' \circ \Psi(\xi)| = \sup_{\xi \in \R^n} |f(\xi)| = \|f\|_{L_\infty(\R^n)}.$$
\end{itemize}

Of course, we will replace P1) by our ADS condition. On the other hand, P2) is not a general property of locally compact groups. Indeed, according to Lemma \ref{amenable} iv) for $\M = \C$ (see the proof), if $\|f\|_{\cL \G}=\|f'\|_{\cL \Gd}$ for any trigonometric polynomial $f$ in $\cL \G$ then the amenability of $\G$ is equivalent to the amenability of $\Gd$. However, this is false in general. Consider for instance the group $\G=SO(3)$ which is compact, hence amenable. On the contrary, since the free group $\mathbb{F}_2$ is a subgroup of $\Gd=SO(3)_{\mathrm{disc}}$, the discretized group $\Gd$ is not amenable. In the following result we show that $\|f\|_{\cL \G}=\|f'\|_{\cL \Gd}$ when $\Gd$ is amenable. 

\begin{lemma}\label{lemma:normLG-LGdisc}
If $f$ is a trigonometric polynomial in $\cL \G$\emph{:}
\begin{enumerate}
\item[i)] We always have $\|f'\|_{\cL \Gd}\leq \|f\|_{\cL \G}$.

\item[ii)] The reverse inequality holds true whenever $\Gd$ is amenable. 
\end{enumerate}
\end{lemma}

\dem Let $(V_j)_{j\geq 1}$ be a symmetric basis of neighborhoods of the identity in $\G$ and let $\mathrm{F}\subset \G$ be finite. Then for $j \ge 1$ large enough and $h_j=\mu(V_j)^{-1/2}\lambda(\ind_{V_j})$ the following map is isometric  
\begin{equation} \label{L2-isometry}
L_{h_j}: \ell_2(\mathrm{F})\ni (a_g)_{g\in \mathrm{F}}\mapsto \Big(\sum_{g\in \mathrm{F}} a_g\lambda(g)\Big)h_j \in L_2(\widehat{\G}).
\end{equation}
Indeed, since $(gV_j)_{g\in \mathrm{F}}$ are disjoint for $j$ large enough $$\|L_{h_j}(a)\|_2^2 \, = \, \mu(V_j)^{-1}\Big\| \sum_{g \in \mathrm{F}} a_g\lambda(\ind_{gV_j})\Big\|_2^2 \, = \, \sum_{g\in \mathrm{F}} |a_g|^2 \, = \, \|a\|_{\ell_2(\mathrm{F})}^2.$$ To prove i), we first write $$\|f'\|_{\cL \Gd} \, = \, \sup \big\langle f' \xi_1, \xi_2 \big\rangle_{\ell_2(\Gd)}$$ where the supremum runs over all finite subset $\mathrm{X}\subset \G$ and all $\xi_1,\xi_2 \in \ell_2(\mathrm{X})$ with $\|\xi_1\|_2=\|\xi_2\|_2=1$. Pick any such $\mathrm{X}$ and $\xi_1, \xi_2$. Since $f' \xi_1$ is supported by $\mathrm{F}\mathrm{X}$ the inner product above can be taken in $\ell_2(\mathrm{S})$, where $\mathrm{S} = \mathrm{F}\mathrm{X} \cup \mathrm{X}$. Applying \eqref{L2-isometry} to this finite set $\mathrm{S}$, we may find an isometry $L_h: \ell_2(\mathrm{S}) \to L_2(\widehat{\G}).$ Since $L_h(f' \xi_1) = f L_h(\xi_1)$ $$\big\langle f' \xi_1, \xi_2 \big\rangle_{\ell_2(\mathrm{S})} \, = \, \big\langle L_h(f' \xi_1), L_h(\xi_2) \big\rangle_{L_2(\widehat{\G})} \, = \, \big\langle f L_h(\xi_1),L_h(\xi_2) \big\rangle_{L_2(\widehat{\G})} \, \le \, \|f\|_{\cL \G}.$$ Taking suprema we obtain i). If $\Gd$ is amenable, Lemma \ref{amenable} yields
\begin{eqnarray*}
\|f\|_{\cL \G} & = & \Big\| \sum_{g \in \mathrm{F}} \widehat{f}(g) \lambda(g) \Big\|_{\V} \\ & \le & \Big\| \sum_{g \in \mathrm{F}} \widehat{f}(g) \lambda(g) \otimes \lambda'(g) \Big\|_{\cL \G \overline{\otimes}\cL \Gd} = \ \|f'\|_{\cL \Gd},
\end{eqnarray*}
where the last equality comes from Fell's absorption principle in Lemma \ref{Fell} ii). \fin

\begin{remark}\label{Gdisc-G amenable}
\emph{It follows that $\Gd$ amenable $\Rightarrow \G$ amenable, but not reciprocally.}
\end{remark}

We can now prove Theorem D i) and ii), the noncommutative version of de Leeuw's compactification theorem. The first implication requires P1) and follows easily  from the lattice approximation in Theorem C. The second one requires an analogue of P2) ---$\Gd$ amenable, as suggested by Lemma \ref{lemma:normLG-LGdisc}--- and it follows by adapting our restriction argument in Theorem A.  

\demD If $\G \in \mathrm{ADS}$ is approximated by lattices $(\Gamma_j)_{j\geq 1}$, then $\Gamma_j \subset \Gd$ for $j\geq 1$. Since both groups are discrete, we may restrict by taking a conditional expectation. In conjunction with Theorem C, we obtain
\begin{eqnarray*}
\big\| T_m: L_p(\widehat{\G}) \to L_p(\widehat{\G}) \big\| & \le & \sup_{j \ge 1} \big\| T_{m_{\mid_{\Gamma_j}}} \hskip-5pt : L_p(\widehat{\Gamma}_j) \to L_p(\widehat{\Gamma}_j) \big\| \\ & \le & \big\| T_m: L_p(\widehat{\Gd}) \to L_p(\widehat{\Gd}) \big\|.
\end{eqnarray*}
This proves i). For the converse implication, we may and will assume as in the proof of Theorem A that $2<p<\infty$. Now, since $\Gd$ is amenable, we claim that $\G \in [\mathrm{SAIN}]_{\Gd}$. Namely, it follows by the exact same argument as in Theorem \ref{AmenableImpliesSAIN} since our proof there does not use the fact that the topology on the subgroup is induced by the topology of $\G$. Once we know that the SAIN condition holds, the goal is to show that $$\|T_m f'\|_{L_p(\widehat{\Gd})} \, \le \, \big\| T_m: L_p(\widehat{\G}) \to L_p(\widehat{\G}) \big\| \, \|f'\|_{L_p(\widehat{\Gd})}$$ for any trigonometric polynomial $f' \in \mathcal{L} \Gd$. Fix such a trigonometric polynomial $f_0' = \sum_{\gamma \in \mathrm{F}} \widehat{f_0'}(\gamma) \lambda'(\gamma) \in \V$ and let $\mathrm{F} \subset \G$ denote its frequency support. Let $(V_j)_{j\geq 1}$ be the neighborhood basis of the identity associated to $\mathrm{F}$ by the SAIN condition. Following the proof of Theorem A, define $h_j = \mu(V_j)^{-1/2}\lambda(\ind_{V_j})$ with polar decomposition $h_j=u_j|h_j|$. The main difference with the restriction theorem is that we may no longer assume that the sets $(gV_j)_{g \in \G}$ are disjoint. Then we cannot define properly the maps $\Phi_j^p$ for all $j$, since they are not contractive any longer. However, this still holds true at the limit. 

\noindent {\bf Claim A'.} Let $2\leq q \leq \infty$. Then 
\begin{enumerate}
\item[i)] If $f' \in \cL \Gd$ is any trigonometric polynomial $$\lim_{j\to \infty} \big\| fu_j|h_j|^{\frac2q} \big\|_{L_q(\widehat{\G})}\leq \|f'\|_{L_q(\widehat{\Gd})}.$$ 

\item[ii)] If $f' \in \cL \Gd$ is frequency supported by $\mathrm{F}$, we also have $$\lim_{j\to \infty} \big\| fu_j|h_j|^{\frac2q} \big\|_{L_q(\widehat{\G})}= \|f'\|_{L_q(\widehat{\Gd})}.$$ 
\end{enumerate}

\noindent The intertwining result we gave in Claim B of the proof of Theorem A ---restated conveniently without using the maps $\Phi_j^p$--- holds replacing $\Gamma$ by $\Gd$ with verbatim the same argument. Moreover, Theorem D ii) follows from it and Claim A' above exactly as in the proof of Theorem A. Thus, it suffices to justify this claim. 

\vskip2pt

\noindent \textbf{Proof of Claim A'.} Let $\varepsilon>0$ and let $f'$ be any trigonometric polynomial in $\cL \Gd$. Since interpolation cannot be used any longer in our case, Claim A' i) will simply follow from the three-lines lemma. Let $a=a(f',\varepsilon,q)$ be a trigonometric polynomial in $\cL \G$ such that
\begin{equation}\label{def-a}
\big\| |f|^{\frac{q}{2}} - a \big\|_{\cL \G} = \big\| |f'|^{\frac{q}{2}} - a' \big\|_{\cL \Gd} \, < \, \frac{1}{2} \, \varepsilon^{q/2},
\end{equation}
where the equality comes from Lemma \ref{lemma:normLG-LGdisc} (together with a standard approximation argument in the weak-$*$ topology) since $\Gd$ is amenable, and $a'$ denotes the trigonometric polynomial in $\cL \Gd$ associated to $a$. By \eqref{L2-isometry}, there exists an index $j_0=j_0(f',\varepsilon,q)$ such that 
\begin{equation}\label{def-j0}
\|ah_j\|_{L_2(\widehat{\G})} \, =Ê\, \|a'\|_{L_2(\widehat{\Gd})} \quad \mbox{ for any } \quad j \ge j_0.
\end{equation}
The map $F_j(z)=u|f|^{qz/2}u_j|h_j|^z$ ---where $f=u|f|$ is the polar decomposition of $f \in \cL \G$--- is holomorphic on the strip $\Delta=\{0<{\rm Re}\, z<1\}$ and continuous on its closure. Since $F_j(it)=u|f|^{iqt/2}u_j|h_j|^{it}$ is a partial unitary $$\sup_{t \in \R} \|F_j(it)\|_{\cL \G} \le 1.$$ On the other hand, by \eqref{def-a} and \eqref{def-j0} we get for all $t\in \R$
\begin{eqnarray*}
\lefteqn{\hskip-15pt \big\| F_j(1+it) \big\|_{L_2(\widehat{\G})}} \\ [4pt] & = & \big\| |f|^{\frac{q}{2}} h_j \big\|_{L_2(\widehat{\G})} \\ & \le & \big\| \big( |f|^{\frac{q}{2}} - a \big) h_j \big\|_{L_2(\widehat{\G})} + \|a h_j\|_{L_2(\widehat{\G})} \ \le \ \frac{\varepsilon^{q/2}}{2} + \|a'\|_{L_2(\widehat{\Gd})} \\ & \le & \frac{\varepsilon^{q/2}}{2} + \big\| a' - |f'|^{\frac{q}{2}} \big\|_{L_2(\widehat{\Gd})} + \big\| |f'|^{\frac{q}{2}} \big\|_{L_2(\widehat{\Gd})} \ \le \ \varepsilon^{q/2} + \|f'\|_{L_q(\widehat{\Gd})}^{q/2}.
\end{eqnarray*}
Therefore, the three-lines lemma implies that for any $j\geq j_0$
\begin{eqnarray}\label{ClaimA-i}
\|F_j(2/q)\|_{L_q(\widehat{\G})} & = & \|fu_j|h_j|^{\frac{2}{q}}\|_{L_q(\widehat{\G})} \\ \nonumber & \le & \Big(\varepsilon^{q/2} + \|f'\|_{L_q(\widehat{\Gd})}^{q/2} \Big)^{\frac{2}{q}} \ \le \ \varepsilon+ \|f'\|_{L_q(\widehat{\Gd})},
\end{eqnarray}
which proves Claim A' i). To prove Claim A' ii) we proceed exactly as in the proof of Theorem A, but using our version of Claim A' i). For a fixed trigonometric polynomial $f'$ in $\cL \Gd$ frequency supported by $\mathrm{F}$, let $k'=k'(f',\varepsilon,q)$ be another trigonometric polynomial in $\cL \Gd$ (frequency supported by $\mathrm{M}\subset \G$ finite) and satisfying $\|k'\|_{L_{q^\ast}(\widehat{\Gd})}=1$ with $$\|f'\|_{L_{q}(\widehat{\Gd})} \le \|k'f'^\ast\|_{L_{2}(\widehat{\Gd})} + \frac{\varepsilon}{2},$$ where $1/q + 1/q^\ast=1/2$. We may choose $j_0=j_0(f',\varepsilon,q)$ such that for any $j\geq j_0$
\begin{enumerate}
\item[i)]  $\|k'f'^\ast\|_{L_{2}(\widehat{\Gd})}=\|k'f'^\ast h_j \|_{L_{2}(\widehat{\G})}$,

\vskip7pt 

\item[ii)] $\big\| k u_j |h_j|^{2/q^\ast} \big\|_{L_{q^*}(\widehat{\G})} \hskip2pt \le \hskip1pt 1+\varepsilon$,

\vskip2pt 

\item[iii)] $\displaystyle\sum_{g \in \mathrm{F}} \displaystyle\frac{\mu \big( g^{-1} V_j g \, \triangle \, V_j \big)}{\mu(V_j)}\leq \displaystyle\frac{\varepsilon^2}{\|k\|_{2}^2\|f\|_2^2|\mathrm{K}|}$.
\end{enumerate}
Namely, the first property follows from \eqref{L2-isometry}, the second one from \eqref{ClaimA-i} and the third one from the SAIN condition. By the same argument as in the proof of Theorem A, we obtain that for any $j \ge j_0$ $$\|f'\|_{L_q(\widehat{\Gd})}\leq \varepsilon + (1+\varepsilon)\|fu_j|h_j|^{2/q}\|_{L_q(\widehat{\G})}.$$ Letting $\varepsilon \to 0^+$, this implies Claim A' ii) and completes Theorem D ii). \fin

\begin{remark}
\emph{According to Remark \ref{examples-ADS}, we know that the Heisenberg group $\H_n$ and the upper triangular matrix groups $\H(\mathbb{K},n)$ are ADS. Moreover, since they are nilpotent the same happens for their discretized forms, which implies in turn that the discretized forms are amenable. In summary, if $\G$ denotes any of these groups, it satisfies the two-sided compactification result in Theorem D i) and ii) for bounded continuous symbols $$\big\| T_m: L_p(\widehat{\G}) \to L_p(\widehat{\G}) \big\| \, = \, \big\| T_m: L_p(\widehat{\Gd}) \to L_p(\widehat{\Gd}) \big\|.$$}
\end{remark}

\section{{\bf The periodization theorem}} \label{sect:periodization}

We finish our collection of noncommutative de Leeuw's theorems in the Banach space setting for unimodular groups with the periodization theorem, nonunimodular groups and statements in the operator space setting will be considered below. In this section we consider a locally compact, unimodular, second countable group $\G$; a normal closed subgroup $\H$ of $\G$; a bounded symbol $m_q: \G/\H \to \C$ and its $\H$-periodization $m_\pi: \G \to \C$ given by $m_\pi(g) = m_q(g\H)$. As mentioned in the Introduction, the abelian case has been solved by Saeki \cite{S} but we cannot go further in the line of Theorem D iii). More precisely, in general $$T_{m_q}: L_p(\widehat{\G/\H}) \to L_p(\widehat{\G/\H}) \ \nRightarrow \ T_{m_\pi}: L_p(\widehat{\G}) \to L_p(\widehat{\G}).$$ Indeed, consider for instance the infinite permutation group $\mathrm{H} = \mathcal{S}_{\infty}$ and construct the cartesian product $\G = \mathbb{T} \times \mathcal{S}_{\infty}$, so that $\G/\H \simeq \mathbb{T}$. By \cite[Proposition 8.1.3]{Pis}, for $1<p\neq 2 <\infty$ we can find a bounded $m_q:\mathbb{T}\to \C$ giving rise to a Fourier multiplier which is bounded in $\ell_p(\Z)$ but not completely bounded. Then, its $\H$-periodization $m_\pi=m_q\otimes id$ cannot define a bounded Fourier multiplier on $$L_p(\widehat{\G})=\ell_p(\Z; L_p(\mathcal{R})),$$ where $\mathcal{R}=\cL \mathcal{S}_\infty$ denotes the hyperfinite II$_1$ factor. Hence, Theorem D iii) fails for this pair $(\G,\H)$. In fact, since Pisier's result on the existence of bounded/not cb multipliers has been extended to any infinite LCA groups \cite{Arh,Ha}, with that process we can construct a large class of counter-examples by taking any group of the form $\G=\mathrm{K} \times \H$ with $\mathrm{K}$ an infinite LCA group and $\H$ a group satisfying that $\cL \H$ contains arbitrarily large matrix algebras $\mathbb{M}_n$. This suggests that there is not so much to do in this direction outside the class of abelian groups. The result in Theorem D iii) was already proved by Saeki \cite{S}. Hence, we now focus on the reverse implication given in Theorem D iv) for $\G$ nonabelian and $\H$ compact. 

\demDD Assume $\H$ is compact and let $\mu_\H$ denote the normalized Haar measure on $\H$. By duality it is enough to consider the case $p \ge 2$. By Lemma \ref{LHsubalgLG}, we may see $\cL \H$ as a von Neumann subalgebra of $\cL \G$ and identify $\lambda_\G(h)$ and $\lambda_\H(h)$ for any $h \in \H$. Consider the operator $$\Pi = \int_\H \lambda(h) \, d\mu_\H(h) \in \cL \H \subset \cL \G.$$ Since $\H$ is a normal, compact (unimodular) subgroup of $\G$, we deduce that $\Pi$ is a central, $\H$-invariant projection of $\cL \G$ onto the functions of $L_2(\G)$ which are constant on $\H$-cosets, denoted by $$\mathcal{H} = \Pi L_2(\G) = \Big\{ \xi \in L_2(\G) : \xi(g) = \xi(g') \mbox{ when } g\H=g'\H \Big\}.$$ The map $\pi: \G \to \M:=(\cL\G)\Pi$ given by $\pi(g)=\lambda(g) \Pi$ defines a $\ast$-representation of $\G$ over the Hilbert space $\mathcal{H}$. Moreover, $\pi$ is invariant on cosets, hence this yields a $\ast$-representation of the quotient $\G /\H$ still denoted by $\pi:\G/\H \to \M$. Observe that $\pi(g\H)=v\lambda_{\G/\H}(g\H)v^\ast$, where the unitary $v:L_2(\G/\H) \to \mathcal{H}$ is the natural identification. Hence $\pi$ can be extended to a normal map $\pi:\cL (\G/\H)\to \M$ by setting $\pi(f)=vfv^\ast$. Since this map is isometric and surjective at the $L_\infty$ and $L_2$ levels, this yields by interpolation an isometric map $$\pi: L_p(\widehat{\G/\H}) \to L_p(\M) = L_p(\widehat{\G}) \Pi$$ for any $2 \le p \le \infty$. On the other hand, $\pi$ intertwines the Fourier multipliers so that $\pi \circ T_{m_q}=T_{m_\pi}\circ \pi$. Indeed, let $f \in  \lambda(\mathcal{C}_c(\G/\H))$. Since the $\G$-invariant measure on left cosets \cite[Theorem 2.49]{F} coincides with the Haar measure on the quotient group $\G/\H$ when $\H$ is normal we get 
\begin{eqnarray*}
\pi \circ T_{m_q}(f) & = & \int_{\G/\H} m_q(g\H) \widehat{f}(g\H) \lambda(g) \Pi \, d\mu_{\G/\H}(g\H) \\ [2pt] & = &\int_{\G/\H} m_q(g\H) \widehat{f}(g\H) \Big( \int_\H \lambda(gh) \, d\mu_\H(h) \Big) \, d\mu_{\G/\H}(g\H) \\ [2pt] & = & \int_\G m_\pi(g) \widehat{f}(g\H) \lambda(g) \, d\mu_{\G}(g) \, \Pi \ = \ T_{m_\pi}\circ \pi (f).
\end{eqnarray*}
Using this property, we conclude with the estimate 
\begin{eqnarray*}
\|T_{m_q}f\|_{L_p(\widehat{\G/\H})} & = & \|\pi \circ T_{m_q}(f)\|_{L_p(\M)} \ = \ \|T_{m_\pi}\circ \pi (f)\|_{L_p(\widehat{\G}) \Pi} \\ [2pt] & \le & \big\| T_{m_\pi}: L_p(\widehat{\G}) \to L_p(\widehat{\G}) \big\| \|\pi (f)\|_{L_p(\widehat{\G}) \Pi} \\ [4pt] & = & \big\| T_{m_\pi}: L_p(\widehat{\G}) \to L_p(\widehat{\G}) \big\| \|f\|_{L_p(\widehat{\G/\H})} 
\end{eqnarray*}
for $f\in L_p(\widehat{\G/\H})$. This completes the proof of Theorem D iv). \fin

\section{\bf Nonunimodular groups} \label{non unimodular}

This section is devoted to extend our results to nonunimodular groups. Again the main focus will be on restriction since compactification and periodization admit less generalizations, see Remark \ref{nonunimodularD}. When $\G$ is nonunimodular, the modular function $\Delta_\G$ is not trivial and the Plancherel weight ---defined in Section \ref{GroupvNaSect} and denoted by $\varphi$ in this section--- is not a trace. This forces to introduce noncommutative $L_p$ spaces associated with arbitrary von Neumann algebras. We will in fact consider two different such $L_p$ spaces, the Haagerup and the Connes-Hilsum ones \cite{Hilsum,TerpI} which turn to be isomorphic as we explain below. Recall that the proof of Theorem A in the unimodular case is based on crucial results derived from Theorem B. Thus we will need to extend these results to arbitrary von Neumann algebras by using Haagerup's reduction method. After that, we will derive Theorem A for nonunimodular groups and give some examples. 

\subsection{Haagerup's reduction for weights}\label{HaagerupReduction}

We start by recalling the reduction method from \cite{HJX} adapted to a von Neumann algebra $\cM \subset \cB(\cH)$ equipped with a fixed normal semifinite faithful (nsf) weight $\varphi$. Note that the constructions in \cite{HJX} are carried out with respect to a normal faithful state $\varphi$ instead of a weight. This is not sufficient for our purposes. The weight case is treated in an unpublished extended version of \cite{HJX} by Xu. For the sake of completeness, we will indicate below the technical modifications of the arguments in \cite{HJX} to obtain the analogous results for weights instead of states. In this paragraph, we consider the so-called Haagerup $L_p$-spaces defined in \cite{TerpI}, see also \cite{HJX} for a standard introduction of the concepts involved. Since they are only used in this auxiliary technical subsection and the next one, we will not detail the construction but refer to the above mentioned works. Let $\sigma^\varphi$ be the modular automorphism group of $\varphi$ and denote $$\mathfrak{n}_\varphi \, = \, \{x \in \mathcal{M}\; :\; \varphi(x^\ast x )<\infty\} \quad \mbox{and} \quad  \mathfrak{m}_\varphi \, = \, \mathfrak{n}_\varphi^\ast \mathfrak{n}_\varphi \, = \, \mathrm{span}\{y^\ast x \; :\; x,y \in  \mathfrak{n}_\varphi\}.$$ In this subsection we fix $\G = \cup_{n \geq 1}2^{-n} \mathbb{Z}$ with the discrete topology and consider the crossed product $\cR = \cM \rtimes_{\sigma^{\varphi}} \G$. Recall that $\cR$ is the von Neumann algebra acting on $L_2(\G, \cH)$ generated by the operators $$\big(\lambda(t)\xi\big)(s) = \xi(s-t)\quad \mbox{and} \quad \big(\pi(x)\xi\big)(s) = \sigma^\varphi_{-s}(x) \xi(s)$$ for $s,t  \in \G, x\in \cM$ and $ \xi \in L_2(\G,\cH)$. We define the unitary operator $$\big(w(\gamma)\xi\big)(s) = \overline{\gamma(s)} \xi(s)$$ for $(s, \gamma, \xi) \in \G \times \widehat{\G} \times L_2(\G,\cH)$ and $\widehat{\alpha}_\gamma(z) = w(\gamma) z w(\gamma)^\ast$ for $z \in \cR$. Then $\pi(\cM)$ is the fixed point algebra for $\widehat{\alpha}$ and the conditional expectation $\cE: \cR \rightarrow \cM$ is given by $\cE(x) = \int_{\widehat{\G}} \widehat{\alpha}_\gamma(x) d\gamma$. The dual weight $\widehat{\varphi}$ on $\cR$ is defined as $\widehat{\varphi} = \varphi \circ \pi^{-1} \circ \cE$. Let $\cR_{\widehat{\varphi}}$ be the centralizer of $\widehat{\varphi}$ in $\cR$ and denote by $\cZ(\cR_{\widehat{\varphi}})$ its center. Consider $$b_n = -i \hskip1pt {\rm Log}(\lambda(2^{-n})) \quad \mbox{and} \quad a_n = 2^n b_n,$$ with ${\rm Log}$ the principal branch of the logarithm, so that $0 \leq {\rm Im( Log}(z)) < 2 \pi$. Then $b_n \in \cZ(\cR_{\widehat{\varphi}})$ and $\varphi_n( \, \cdot \, ) = \widehat {\varphi}( e^{-a_n} \, \cdot \, )$ formally defines a nsf weight. More precisely, $\varphi_n$ has Connes cocycle derivative $(D\varphi_n/ D \widehat{\varphi})_s = e^{-i s a_n}$ for $s \in \mathbb{R}$.

\begin{theorem} \label{Thm=Reduction}
Let $\cR_n$ be the centralizer of $\varphi_n$ in $\cR$. The sequence $(\cR_n)_{n \geq 1}$ forms an increasing sequence of von Neumann subalgebras of $\cR$. Moreover, the following properties hold\emph{:}
\begin{enumerate}
\item[i)] $\cR_n$ is semifinite for each $n \geq 1$ with trace $\varphi_n$.

\vskip2pt

\item[ii)] There exist conditional expectations $\cE_n: \cR \rightarrow \cR_n$ such that $$\widehat{\varphi} \circ \cE_n = \widehat{\varphi} \quad \mbox{and} \quad \cE_n \circ \sigma_s^{\widehat{\varphi}} = \sigma_s^{\widehat{\varphi}} \circ \cE_n \quad \mbox{ for all } \quad s \in \mathbb{R}.$$

\vskip2pt

\item[iii)] $\cE_n(x) \rightarrow x$ $\sigma$-strongly for $x \in \mathfrak{n}_{\widehat{\varphi}}$ and $\bigcup_{n \geq 1} \cR_n$ is $\sigma$-strongly dense in $\cR$.
\end{enumerate}
\end{theorem}

\dem The proof is a mutatis mutandis copy of the arguments in \cite[Section 2]{HJX}. We indicate the main adaptations. Observe that \cite[Lemma 2.2]{HJX} does not remain valid. This lemma is applied only in two places, where the arguments need to be adapted. Firstly, it is needed to prove the uniqueness of $b_n$ in \cite[Lemma 2.3]{HJX}, but this does not play a role in the subsequent proofs. Secondly it is used in the proof of \cite[Lemma 2.6]{HJX}. However, we claim that the following fact still holds true: for every $x \in \mathfrak{n}_{\widehat{\varphi}}$ and every $\varepsilon > 0$ there exists a trigonometric polynomial $P$ on $\mathbb{T}$ with 
\begin{equation}\label{Eqn=PConvergence}
\big\Vert \big[b_n - P(\lambda(2^{-n})), x \big]  \big\Vert_{\widehat{\varphi}} \leq \varepsilon \quad \mbox{ for all } \quad n \in \N,
\end{equation}
where $[x,y]=xy-yx$ denotes the commutator of two operators $x$ and $y$ and $\|y\|_{\widehat{\varphi}}^2 =\widehat{\varphi}(y^\ast y)$ for any $y\in \cR$. This fact is what is actually needed. Let us now prove it. If $x \in \mathfrak{n}_{\widehat{\varphi}}$, then $$\big\Vert \big(b_n - P(\lambda(2^{-n})) \big) x \big\Vert^2_{\widehat{\varphi}} \, = \, \widehat{\varphi} \big( x^\ast \vert b_n -  P(\lambda(2^{-n})) \vert^2   x \big).$$ Now $\widehat{\varphi}( x^\ast \cdot x)$ is a normal functional on $\cR$ and hence it restricts to a normal functional $\omega$ on the von Neumann subalgebra generated by $\lambda(2^{-n})$, which equals $L_\infty(\mathbb{T})$. So $\omega$ corresponds to integration against a function in $L_1(\mathbb{T})$. Recalling that $b_n = -i \hskip1pt {\rm Log}(\lambda(2^{-n}) )$ we see that we may choose $P$ such that for every $n$ we have $\omega(   \vert b_n -  P(\lambda(2^{-n})) \vert^2 ) < \varepsilon$. On the other hand, we first consider $$x \in \cT_{\widehat{\varphi}} := \Big\{ x \in \cR \;:\; x \textrm{ is analytic for } \sigma^{\widehat{\varphi}} \textrm{ and } \sigma^{\widehat{\varphi}}_z(x) \in \nphi \cap \nphi^\ast \ \forall \, z \in \mathbb{C} \Big\}.$$ In that case, from Tomita-Takesaki theory we have $$\big\Vert x \big(b_n - P(\lambda(2^{-n}) )\big) \big\Vert^2_{ \widehat{\varphi}} \, = \, \widehat{\varphi}\big( x  \vert b_n - P(\lambda(2^{-n}) ) \vert^2 \sigma_{-i}^{\widehat{\varphi}} (x^\ast) \big),$$ and as above we may find $P$ such that for every $n$ this expression becomes smaller than $\varepsilon$. This proves our claim \eqref{Eqn=PConvergence} in case $x \in \cT_{\widehat{\varphi}}$. For a general operator $x \in \mathfrak{n}_{\widehat{\varphi}}$ the claim follows by taking a net $(x_j)_{j \in J}$ in $\cT_{\widehat{\varphi}}$ such that $\Vert x_j - x \Vert_{\widehat{\varphi}} \rightarrow 0$ (see for instance \cite{TakII}) and using that $\Vert b_n \Vert \leq 2 \pi$.

Let us now return to the constructions of \cite[Section 2]{HJX}. The statements and proofs of \cite[Lemmas 2.3, 2.4, 2.5]{HJX} remain unchanged except that $b_n$ might not be unique, which is not relevant for the proof. Note in particular that the restriction of $\varphi_n$ to its centralizer is semifinite. Then Lemma 2.6 remains true provided $x \in \mathfrak{n}_{\widehat{\varphi}}$ instead of general $x \in \cR$ and also Lemma 2.7 remains valid for $x \in \mathfrak{n}_{\widehat{\varphi}}$. Indeed, as in the proof of Lemma 2.7, this follows from Lemma 2.6 in case $x \in \mathfrak{n}_{\widehat{\varphi}}$ (and also in the weight case one invokes  Lemma 2.5 to derive strong convergence, which implies $\sigma$-strong convergence for a bounded net). This completes the proof. \fin

Let $L_p(\cM), L_p(\cR)$ and  $L_p(\cR_n)$ be the Haagerup $L_p$-spaces constructed from the weights $\varphi, \widehat{\varphi}$, and $\widehat{\varphi}$ restricted to $\cR_n$ respectively, see \cite {TerpI} or \cite[Section 1.2]{HJX}. The modular automorphism group $\sigma^{\widehat{\varphi}}$ restricted to $\cM \simeq \pi(\cM)$ equals $\sigma^\varphi$. By Theorem \ref{Thm=Reduction}, the restriction of $\widehat{\varphi}$ to $\cR_n$ is semifinite. This implies that the crossed products $\cM \rtimes_{\sigma^{\varphi}} \mathbb{R}$ and $\cR_n \rtimes_{\sigma^{\widehat{\varphi}}} \mathbb{R}$ are well-defined subalgebras of $\cR \rtimes_{\sigma^{\widehat{\varphi}}} \mathbb{R}$. Let $D$ be the generator of the left regular representation in each of these crossed products, then $D$ is the usual density operator in the Haagerup $L_p$-space $L_p(\cR)$. Recall that we have two $\widehat{\varphi}$-preserving conditional expectations $\mathcal{E}: \mathcal{R} \rightarrow \mathcal{M}$ and $\mathcal{E}_n: \mathcal{R} \rightarrow \mathcal{R}_n$. For $1\leq p <\infty$, by Remark 5.6 and Example 5.8 of \cite{HJX} we obtain  contractive projections $$\mathcal{E}^p: L_p(\mathcal{R}) \rightarrow L_p(\mathcal{M}) \quad \mbox{and} \quad \mathcal{E}_n^p: L_p(\mathcal{R}) \rightarrow L_p(\mathcal{R}_n)$$ given by $\mathcal{E}^p(D^{\frac{1}{2p}} x D^{\frac{1}{2p}} ) = D^{\frac{1}{2p}} \mathcal{E}(x) D^{\frac{1}{2p}}$ for any $x \in \mathfrak{m}_{\hat{\varphi}}$, and similarly for $\mathcal{E}_n$. More generally, for any $p\leq r,s\leq \infty$ such that $\frac{1}{r}+\frac{1}{s}=\frac{1}{p}$ we have $\mathcal{E}^p(D^{\frac{1}{r}} x D^{\frac{1}{s}} ) = D^{\frac{1}{r}} \mathcal{E}(x) D^{\frac{1}{s}}$ for $x \in \mathfrak{m}_{\hat{\varphi}}$, see \cite[Proposition 5.5]{HJX} for the proof in the state case.

\begin{remark}\label{Rmk=LpConvention}
{\rm The notation $D^{\frac{1}{r}} x D^{\frac{1}{s}}$ for $x \in \mathfrak{m}_{\hat{\varphi}}$ used in \cite{HJX} and which we keep using in the sequel is formal. If $x$ can be decomposed as a finite sum $x = \sum_{j} y_j^\ast z_j$ with $y_j, z_j \in \mathfrak{n}_{\hat{\varphi}}$, then the notation $D^{\frac{1}{r}} x D^{\frac{1}{s}}$ stands for $\sum_j D^{\frac{1}{r}} y_j^\ast \cdot [z_j D^{\frac{1}{s}}]$, which is a well-defined element of $L^p(\cR)$ by \cite[Theorem 26]{TerpII} and H\"older's inequality. Here $[\; \cdot\;]$ denotes the closure of a preclosed operator. Arguing as in \cite{GoldsteinLindsay} one can derive that this expression does not depend on the decomposition of $x$.}
\end{remark}

\begin{lemma} \label{Lem=LpApproxI} 
Given $1 \leq p < \infty$ and $x \in L_p(\cR)$ we have  $$\lim_{n \to \infty} \Vert \cE_n^p(x) - x \Vert_{p}  = 0.$$ 
\end{lemma}

\dem We first assume that $1 \leq p \leq 2$. Let $x \in \cR$ and $x', x''  \in \cR_m$ for some $m \geq 1$. Assume moreover that $x,  x''  \in \mathfrak{n}_{\widehat{\varphi}}, x'\in \mathfrak{n}_{\widehat{\varphi}}^\ast$ and $n \geq m$. Let $r$ be such that $\frac{1}{p}-\frac{1}{2} = \frac{1}{r}$. Since $\cE_n(x') = x'$ and $\cE_n(x'') = x''$, using the convention of Remark \ref{Rmk=LpConvention}, H\"older's inequality and \cite[Theorem 23]{TerpII} imply  
\begin{eqnarray*}
\big\| \cE_n^p(  D^{\frac{1}{r}} x'  x x'' D^{\frac{1}{2}}  )- D^{\frac{1}{r}} x' x x'' D^{\frac{1}{2}}  \big\|_p & = & \big\| D^{\frac{1}{r}} x' \cE_n(x) x'' D^{\frac{1}{2}} - D^{\frac{1}{r}} x'  x x'' D^{\frac{1}{2}}  \big\|_p \\ & \le & \Vert D^{\frac{1}{r}} x'\Vert_{r} \, \big\| (\cE_n(x) - x)   x'' D^{\frac{1}{2}} \big\|_{2} \\ [2pt] & = & \Vert D^{\frac{1}{r}} x'\Vert_{r} \, \big\| (\cE_n(x) - x) \Lambda( x'')  \big\|_{2}  \rightarrow 0,
\end{eqnarray*}
since $\cE_n(x) \rightarrow x$ strongly by Theorem \ref{Thm=Reduction}. Here $\Lambda$ denotes the canonical injection of $\mathfrak{n}_{\widehat{\varphi}}$ into its Hilbert space completion. We claim that the linear span of elements $D^{1/r} x'   x x'' D^{1/2}$ with $x, x', x''$ as above is dense in $L_p(\cM)$. Then the result will follow for any operator $x \in L_p(\cR)$ for $1\le p \leq 2$ by contractivity of $\cE_n^p$. By \cite[Theorem 26]{TerpII} the linear span of $ D^{1/r} y^\ast = [y D^{1/r}]^\ast$ with $y \in \mathfrak{n}_{\widehat{\varphi}}$ is dense in $L_r(\cM)$ for $2 \le r < \infty$. Then the H\"older inequality gives that $\mathrm{span}\{D^{1/r} x D^{1/2}\;: \; x \in \mathfrak{m}_{\widehat{\varphi}}\}$ is dense in $L_p(\cM)$. Let $(x'_j)_{j \in J}, (x''_j)_{j \in J}$ be nets in $\cR_n$  of elements that are analytic for $\sigma^{\widehat{\varphi}}$ and such that $$\sigma_z^{\widehat{\varphi}}(x_j'), \sigma_z^{\widehat{\varphi}}(x_j'')  \in \mathfrak{n}_{\widehat{\varphi}} \cap \mathfrak{n}_{\widehat{\varphi}}^\ast$$ for every $z \in \mathbb{C}$. Assume that $\sigma_{-i/r}^{\widehat{\varphi}}(x_j') \rightarrow 1$ and $\sigma_{i/2}^{\widehat{\varphi}}(x_j'') \rightarrow 1$ strongly. Then using \cite[Lemma 2.5]{GoldsteinLindsay} and \cite[Lemma 2.3]{JungeSherman} we get $$D^{\frac{1}{r}} x'_j x x''_j D^{\frac{1}{2}} = \sigma_{-i/r}^{\widehat{\varphi}}(x_j') \cdot  D^{\frac{1}{r}} x  D^{\frac{1}{2}} \cdot \sigma_{i/2}^{\widehat{\varphi}}(x_j'')  \rightarrow    D^{\frac{1}{r}} x  D^{\frac{1}{2}},$$ in the norm of $L_p(\cM)$. Here, the domains of the operators in the first equality are equal by an argument similar to the one we will use to prove Lemma \ref{Lem=LpCommutation}.This concludes our claim, and hence the Lemma for $1 \leq p \leq 2$. We now consider the case $p \ge 2$. Suppose that we have proved the Lemma for $p/2$. Take $x \in \cR$ such that $x \in \mathfrak{n}_{\widehat{\varphi}}$. By H\"older's inequality
\begin{eqnarray*}
\big\| (\cE_n(x) - x) D^{\frac{1}{p}} \big\|_p^2 & = & \big\| (\cE_n(x) - x) \cdot D^{\frac{2}{p}} (\cE_n(x) - x)^\ast  \big\|_{p/2} \\ & \le & \big\| \cE_n(x) - x  \Vert_\infty \Vert  (\cE_n(x) - x) D^{\frac{2}{p}}    \big\|_{p/2} \\ & = & \big\| \cE_n(x) - x  \Vert_\infty  \Vert  \cE_n^p(  x  D^{\frac{2}{p}} )  - x  D^{\frac{2}{p}} \big\|_{p/2}, 
\end{eqnarray*}
which goes to $0$ as $n$ tends to $\infty$. Therefore, the result for a general operator $x \in L_p(\cR)$ follows by density \cite[Theorem 26]{TerpII}, recalling that $\cE_n^p$ is contractive. \fin

\subsection{Almost multiplicative maps on arbitrary von Neumann algebras} \label{Sect=MultiplicativeTypeIII}

We now apply the reduction method detailed above to the results of Section \ref{Sect=AlmostMultiplicative} needed to prove Theorem A in the nonunimodular setting. Let $\cM$ be a von Neumann algebra with a nsf weight $\varphi$ and $T: \cM \rightarrow \cM$ be a positive map such that $\varphi \circ T \leq \varphi$. Given $1 \le p < \infty$ and according to \cite[Remark 5.6]{HJX}, the map $T$ induces a bounded map $T_p$ on the Haagerup $L_p$-space $L_p(\cM)$ determined by $$T_p(D_{\varphi}^{\frac{1}{2p}} x D_{\varphi}^{\frac{1}{2p}})=D_{\varphi}^{\frac{1}{2p}}T( x) D_{\varphi}^{\frac{1}{2p}}$$ for $x \in \mathfrak{m}_{\varphi}$, where $D_{\varphi}$ denotes the density operator of $\varphi$. With that notation, we can state and prove the following analogues of Corollary \ref{Tx-x}  and Corollary \ref{Ty-y} for arbitrary von Neumann algebras. 

\begin{corollary}  \label{Tx-x arbitrary vNA} 
Let $\cM$ be a von Neumann algebra equipped with a nsf weight $\varphi$ and let $T: \M \to \M$ be a subunital completely positive map with $\varphi \circ T \leq \varphi$ and $T \circ \sigma^\varphi_s = \sigma^\varphi_s \circ T$ for every $s \in \mathbb{R}$. Then there exists a universal constant $C > 0$ such that the following inequality holds for any $x \in L_2^+(\M)$ and any $0 < \theta \le 1$ $$\big\| T_{\frac{2}{\theta}}(x^\theta) - x^\theta \big\|_{\frac{2}{\theta}} \, \le \, C \, \big\| T_2(x) - x \big\|_2^{\frac{\theta}{2}} \|x\|_2^{\frac{\theta}{2}}.$$
\end{corollary}

\dem We use the notations of Section \ref{HaagerupReduction}. By \cite[Section 4]{HJX}, we know that the map $T$ admits a subunital completely positive normal extension, which is given by $$\widehat{T}: \cR \ni \pi(x)\lambda(s) \mapsto \pi(T(x))\lambda(s) \in \cR$$ for any $(s,x) \in \G \times \cM$. Note that $\cL \G$ is in the multiplicative domain of $\widehat{T}$. Moreover, we also have $$\widehat{\varphi} \circ \widehat{T} \leq \widehat{\varphi} \quad \mbox{and} \quad \sigma_t^{\widehat{\varphi}} \circ \widehat{T}=\widehat{T}\circ \sigma_t^{\widehat{\varphi}}.$$ Recall that $\sigma_t^{\varphi_n}$ and $\cE_n$ are defined in \cite{HJX} respectively by $$\sigma_s^{\varphi_n}(x) = e^{-isa_n}\sigma_s^{\widehat{\varphi}} (x) e^{isa_n} \quad \mbox{and} \quad \cE_n(x)=2^n\int_0^{2^{-n}}\sigma_s^{\varphi_n}(x)ds$$ for any $(x,s) \in \cR \times \R$. Note that these expressions were used in \cite{HJX} for states, although the same construction is valid for weights and the resulting conditional expectations commute with the action of the modular automorphism group. Since $e^{isa_n}\in \cL \G$, we deduce that $\widehat{T}$ commutes with $\cE_n$. Hence, we may consider its restriction to $\cR_n$ and deduce that we still have that $$\varphi_n\circ\widehat{T} \le \varphi_n.$$ By Theorem \ref{Thm=Reduction} i) $(\cR_n, \varphi_n)$ is semifinite and we may extend $\widehat{T}$ to a contractive map on the tracial $L_p$-space $L_p(\cR_n,\varphi_n)$. This extension does not depend on $p$. On the other hand, for $1 \le p < \infty$ the map given by $$\widehat{T}_p(D_{\widehat{\varphi}}^{1/2p}xD_{\widehat{\varphi}}^{1/2p})=D_{\widehat{\varphi}}^{1/2p}\widehat{T}(x)D_{\widehat{\varphi}}^{1/2p} \quad \mbox{for} \quad x \in \mathfrak{m}_{\widehat{\varphi}}$$ extends to a bounded map $\widehat{T}_p: L_p(\cR,\widehat{\varphi})\to L_p(\cR,\widehat{\varphi})$  by \cite[Remark 5.6]{HJX}. Since $\cE \circ \widehat{T}=T \circ \cE$, where $\cE:\cR \to \cM$ is the $\widehat{\varphi}$-preserving conditional expectation, the restriction of $\widehat{T}_p$ to $L_p(\cM)$ equals $T_p$. Moreover, it commutes with $\cE_n^p$ and we may consider the restriction $$\widehat{T}_p:L_p(\cR_n,\widehat{\varphi})\to L_p(\cR_n,\widehat{\varphi}).$$ As it is proved in \cite{TerpI}, we have $L_p(\cR_n,\widehat{\varphi}) \simeq L_p(\cR_n,\varphi_n)$ isometrically, and the isomorphism preserves positive elements. The two restriction maps $\widehat{T}$ and $\widehat{T}_p$ are compatible with respect to that isomorphism. Namely, let $\kappa_p:L_p(\cR_n,\widehat{\varphi}) \to L_p(\cR_n,\varphi_n)$ be the isometric isomorphism given by $\kappa_p(D_{\widehat{\varphi}}^{1/2p}xD_{\widehat{\varphi}}^{1/2p})=e^{\frac{a_n}{2p}}xe^{\frac{a_n}{2p}}$ for any $x\in  \mathfrak{m}_{\widehat{\varphi}}$, then $$\kappa_p \circ \widehat{T}_p=\widehat{T}\circ \kappa_p \quad \mbox{on} \quad L_p(\cR_n,\widehat{\varphi})$$ since $a_n$ lies in the multiplicative domain of $\widehat{T}$. Fix $x\in L_2^+(\cM)$, then by Lemma \ref{Lem=LpApproxI} iii)  and the fact that $\widehat{T}_p$ commutes with $\cE_n^p$ we can write
\begin{eqnarray*}
\big\|T_{\frac{2}{\theta}}(x^\theta)-x^\theta\big\|_{L_{\frac{2}{\theta}}(\cM)} & = & \big\|\widehat{T}_{\frac{2}{\theta}}(x^\theta)-x^\theta\big\|_{L_{\frac{2}{\theta}}(\cR)} \\ & = & \lim_{n\to \infty} \big\|\cE_n^{\frac{2}{\theta}}\circ \widehat{T}_{\frac{2}{\theta}}(x^\theta)-\cE_n^{\frac{2}{\theta}}(x^\theta)\big\|_{L_{\frac{2}{\theta}}(\cR_n,\widehat{\varphi})} \\ & = & \lim_{n\to \infty} \big\|\widehat{T}_{\frac{2}{\theta}}\big(\cE_n^{\frac{2}{\theta}}(x^\theta)\big)-\cE_n^{\frac{2}{\theta}}(x^\theta)\big\|_{L_{\frac{2}{\theta}}(\cR_n,\widehat{\varphi})} \\ & = & \lim_{n\to \infty} \Big\|\widehat{T}\Big(\kappa_{\frac{2}{\theta}}\big(\cE_n^{\frac{2}{\theta}}(x^\theta)\big)\Big)-\kappa_{\frac{2}{\theta}}\big(\cE_n^{\frac{2}{\theta}}(x^\theta)\big)\Big\|_{L_{\frac{2}{\theta}}(\cR_n,\varphi_n)}.
\end{eqnarray*}
By Corollary \ref{Tx-x} in $(\cR_n,\varphi_n)$ applied to the map $\widehat{T}$ and to $\kappa_{\frac{2}{\theta}}\big(\cE_{n}^{\frac{2}{\theta}}( x^\theta)\big)^{\frac{1}{\theta}}\in L_2^+(\cR_n,\varphi_n)$
\begin{eqnarray*}
\lefteqn{\big\| T_{\frac{2}{\theta}}(x^\theta) - x^\theta \big\|_{\frac{2}{\theta}}} \\ \!\!\!\! & \le & \!\!\!\! C \lim_{n \to \infty} \Big\| \widehat{T} \Big( \kappa_{\frac{2}{\theta}} \big(\cE_n^{\frac{2}{\theta}}(x^\theta)\big)^{\frac{1}{\theta}}\Big) - \kappa_{\frac{2}{\theta}} \big(\cE_{n}^{\frac{2}{\theta}}( x^\theta)\big)^{\frac{1}{\theta}} \Big\|_{L_2(\cR_n,\varphi_n)}^{\frac{\theta}{2}} \Big\| \kappa_{\frac{2}{\theta}} \big(\cE_{n}^{\frac{2}{\theta}}( x^\theta)\big)^{\frac{1}{\theta}} \Big\|_{L_2(\cR_n,\varphi_n)}^{\frac{\theta}{2}} \\ \!\!\!\! & = & \!\!\!\! C \lim_{n \to \infty} \big\|\widehat{T}_{2}\big(\cE_n^{\frac{2}{\theta}}(x^\theta)^{\frac{1}{\theta}}\big) - \cE_n^{\frac{2}{\theta}}(x^\theta)^{\frac{1}{\theta}} \big\|_{L_{2}(\cR_n,\widehat{\varphi})}^{\frac{\theta}{2}} \big\|\cE_n^{\frac{2}{\theta}}(x^\theta)^{\frac{1}{\theta}}\big\|_{L_{2}(\cR_n,\widehat{\varphi})}^{\frac{\theta}{2}}.
\end{eqnarray*}
We finally claim that $$\lim_{n \to \infty} \big\|\cE_n^{\frac{2}{\theta}}(x^\theta)^{\frac{1}{\theta}}-x \big\|_2,$$ which yields the result since $\Vert \widehat{T}_2(x) - x \Vert_2 = \Vert T_2(x) - x \Vert_2$. This claim follows from Lemma \ref{Lem=LpApproxI} iii) and the fact that for any operators $x,y \in L_2(\cR)$ such that $\|y\|_2\leq \|x\|_2$ and any parameter $0<\theta\leq 1$ we have $$ \|x-y\|_2 \, \le \, k \big\| x^\theta - y^\theta \big\|_{\frac{2}{\theta}}^{\frac{1}{\theta k}} \|x\|_2^{1-\frac{1}{k}}$$ for any integer $k\geq 1$ satisfying $\frac1k\leq \theta$. Indeed, we first observe that for $1 \le p \le \infty$ and $k \in \Z_+$ $$\big\| x^k - y^k \big\|_{p} \, \le \, k \|x-y\|_{pk} \|x\|_{pk}^{k-1} \quad \mbox{for} \quad x,y \in L^+_{pk}(\cR)$$ with $\|y\|_{pk} \le \|x\|_{pk}$. This easily follows from H\"older inequality and the identity $$x^k-y^k=\sum_{j=0}^{k-1} x^{k-j-1}(x-y)y^{j}.$$ Then we get $$\|x-y\|_2 \leq k\|x^{\frac{1}{k}}-y^{\frac{1}{k}}\|_{2k}\|x\|_2^{(k-1)/k} \leq k\|x^\theta-y^{\theta}\|_{2/\theta}^{1/\theta k}\|x\|_2^{(k-1)/k}.$$ The last inequality follows from the Powers-St{\o}rmer inequality Lemma \ref{hardps}. Note that we have not justified the validity of such inequality for type III algebras. It is however a simple exercise to deduce it from \cite[Proposition 7 and Lemma B]{Kosaki}. \fin

\begin{corollary} \label{Ty-y arbitrary vNA} 
Let $\cM$ be a von Neumann algebra equipped with a nsf weight $\varphi$ and let $T: \M \to \M$ be a subunital completely positive map with  $\varphi \circ T \leq \varphi$ and $T \circ \sigma^\varphi_s = \sigma^\varphi_s \circ T$ for every $s \in \mathbb{R}$. Then there exists a universal constant $C > 0$ such that the following inequality holds for any self-adjoint $y \in L_2(\M)$ with polar decomposition $y = u |y|$ and any $0 < \theta \le 1$ $$\big\| T_{\frac{2}{\theta}}(u|y|^\theta) - u|y|^\theta \big\|_{\frac{2}{\theta}} \, \le \, C \, \big\| T_2(y) - y \big\|_2^{\frac{\theta}{4}} \|y\|_2^{\frac{3\theta}{4}}.$$
\end{corollary}

\dem The proof is similar to the one of Corollary \ref{Ty-y}, details are omitted. \fin

\subsection{Connes-Hilsum $L_p$ spaces}

In this subsection we recall the construction for group von Neumann algebras of Connes-Hilsum $L_p$-spaces \cite{Hilsum}, since we shall use them in the proof of Theorem A. This construction will also be needed in the next section, in order to apply the transference results from \cite{CdlS} in the category of operator spaces. Since our proof of Theorem A will rely on the results derived from Theorem B established in Section \ref{Sect=MultiplicativeTypeIII} for the Haagerup $L_p$-spaces, we need to compare both constructions.

\noindent \textbf{The Connes-Hilsum construction for group algebras.} We shall follow the presentation of \cite{CdlS}. Let $\G$ be a locally compact group and let $\rho:\G \to \cB(L_2(\G))$ be the right regular representation $$\rho(g)(\xi)(h)=\Delta_\G(g)^{\frac12}\xi(hg)$$ for any $\xi \in L_2(\G)$ and $g, h\in \G$. Set $$\rho(\xi)= \int_\G \xi(g)\rho(g)d\mu(g) \quad \mbox{for any } \xi \in L_2(\G).$$ There exists a nsf weight $\varphi'$ on the commutant $\cL \G '=\rho(\G)''$ given by $$\varphi'(f^*f) \, = \, \int_\G |\xi(g)|^2 \, d\mu(g)$$ when $f = \rho(\xi)$ for some $\xi \in L_2(\G)$ and $\varphi'(f^*f) = \infty$ for any other $f \in \cL \G'$. For a nsf weight $\omega$ on $\cL \G$, the partial derivative $(d\omega/d\varphi')^{\frac{1}{2}}$ is the unique closed densely defined operator, whose domain consists of the left bounded functions in $L_2(\G)$ and such that $$\big\Vert (d\omega/d\varphi')^{\frac{1}{2}} \xi\big\Vert_{L_2(\G)}^2 = \omega(\lambda(\xi) \lambda(\xi)^\ast ) < \infty.$$ For $1\leq p <\infty$, the Connes-Hilsum noncommutative space $L_p(\widehat{\G})=L_p(\widehat{\G},\varphi')$ is then defined as the set of closed densely defined operators $f$ on $L_2(\G)$ with polar decomposition $f = u \vert f \vert$ such that $u \in \cL \G$ and $\vert f \vert^p$ equals $d\omega / d\varphi'$ for some $\omega \in \cL \G_\ast$. In that case $$\Vert f \Vert_{L_p(\widehat{\G})} \, = \, \omega(1)^{\frac1p} \, = \, \Vert \omega \Vert^{\frac1p}.$$ Equipped with this norm, $L_p(\widehat{\G})$ is a Banach space and the H\"older inequality holds by understanding the product of two operators as the closure of their product. For $\xi \in L_1(\G) \cap L_2(\G)$, we have the Plancherel formula $$[ \lambda(\xi) \Delta_\G^{\frac12}] \in L_2(\widehat{\G}) \qquad \textrm{ with } \qquad \big\Vert [ \lambda(\xi) \Delta_\G^{\frac12}] \big\Vert_{L_2(\widehat{\G})} = \Vert \xi \Vert_{L_2(\G)}.$$ In fact, such elements are dense in $L_2(\widehat{\G})$. Moreover, the set of operators $$\big\{[ \lambda(\xi) \Delta_\G^{\frac1p}]\;:\; \xi \in \mathcal{C}_c(\G) \big\}$$ is dense in $L_p(\widehat{\G})$ for $2 \le p < \infty$, see \cite[Theorem 26]{TerpII}. Connes-Hilsum $L_p$-spaces are compatible with interpolation, meaning that we may find a compatible structure so that the family $(L_p(\widehat{\G}))_{1\leq p \leq  \infty}$ forms an interpolation scale, further details can be found in \cite{CdlS}. Let $2 \leq p \leq \infty$ and consider any symbol $m \in L_\infty(\G)$. Given any $\xi \in \mathcal{C}_c(\G)$, we have $$[\lambda(\xi) \Delta_{\G}^{1/p}], [\lambda(m \xi) \Delta_{\G}^{1/p}] \in L_p(\widehat{\G}).$$ Then we may consider the associated multiplier $$T_m^p:  L_p(\widehat{\G}) \ni [\lambda(\xi) \Delta_{\G}^{1/p}] \mapsto [\lambda(m \xi) \Delta_{\G}^{1/p}] \in  L_p(\widehat{\G}),$$ which is called an $L_p$-Fourier multiplier if it extends  boundedly to $L_p(\widehat{\G})$ (to a normal map if $p=\infty$). For $1 \leq p \leq 2$ and a given bounded symbol $m$, we define the associated multiplier by $$T_m^1 := (T_{m_{\mathrm{op}}}^{\infty})_\ast \quad \mbox{and} \quad T_m^p = (T_{m_{\mathrm{op}}}^{p'})^\ast \quad \mbox{where} \quad m_{\mathrm{op}}(s) = m(s^{-1}).$$ 

\noindent \textbf{Relation between Haagerup and Connes-Hilsum spaces.} Let us fix some notation. We let $\cM$ be a von Neumann algebra equipped with nsf weight $\varphi$. Let $\varphi'$ be a nsf weight on the commutant $\cM'$. We let $L_p(\cM)$ be the Connes-Hilsum $L_p$ space constructed from $\varphi'$. Let $d = d\varphi/ d\varphi'$ be the spacial derivative. Let $x \in \mphi$ and write $x = \sum_i y_i^\ast z_i$ with $y_i, z_i \in \nphi$ (finite sum). Define $$j_p(x) = \summ_i d^{\frac{1}{2p}} y_i^\ast \cdot [z_id^{\frac{1}{2p}}] \in L_p(\cM).$$ The sum above does not depend on the choice of $y_i$ and $z_i$, see \cite{GoldsteinLindsay}. We let $L_p(\cM)_\rtimes$ be the Haagerup $L_p$-space constructed from $\varphi$ with density operator $D$. Define also $$j_{p, \rtimes}(x) = \summ_i D^{\frac{1}{2p}} y_i^\ast \cdot [z_i D^{\frac{1}{2p}}] \in L_p(\cM).$$ Let $\mathcal{S}_\varphi$ be the set of all $x \in \cM$ such that $x$ is analytic for $\sigma^\varphi$ and $\sigma^\varphi_z(x) \in \mphi$ for every $z \in \mathbb{C}$. Recall that if $a$ is a closed unbounded operator and $b$ is a bounded operator then $ab$ is automatically closed.

\begin{lemma}\label{Lem=LpCommutation}
For every $x \in \mathcal{S}_\varphi$ we have $$j_p(x) = d^{\frac{1}{p}} \sigma^\varphi_{\frac{i}{2p}}(x) = [\sigma^\varphi_{-\frac{i}{2p}}(x) d^{\frac{1}{p}}].$$
\end{lemma}

\dem Let $y_i, z_i \in \nphi$ be such that $$x = \summ_i y_i^\ast z_i.$$ Using \cite[Lemma 22]{TerpII} for the first inclusion and an elementary inclusion $$\sigma_{-\frac{i}{2p}}^\varphi(x) d^{\frac{1}{p}} \subseteq d^{\frac{1}{2p}} x d^{\frac{1}{2p}} \subseteq \summ_{i}  d^{\frac{1}{2p}} y_i^\ast \cdot [z_i d^{\frac{1}{2p}}] \in L_{p}(\cM).$$ Hence $$\big( \sigma_{-\frac{i}{2p}}^\varphi(x) d^{\frac{1}{p}} \big)^\ast \supseteq \summ_{i}  \big( d^{\frac{1}{2p}} y_i^\ast \cdot [z_i d^{\frac{1}{2p}}] \big)^\ast \in L_p(\cM).$$ By (the proof of) \cite[Theorem 4 (1)]{Hilsum} we in fact have an equality $$\big( \sigma_{-\frac{i}{2p}}^\varphi(x) d^{\frac{1}{p}} \big)^\ast = \summ_{i}  \big( d^{\frac{1}{2p}} y_i^\ast \cdot [z_i d^{\frac{1}{2p}}] \big)^\ast.$$ Therefore, taking adjoints yields the equality in the next line $$d^{\frac{1}{p}} \sigma_{\frac{i}{2p}}^\varphi(x) \supseteq [\sigma_{-\frac{i}{2p}}^\varphi(x) d^{\frac{1}{p}}] = \summ_{i}  d^{\frac{1}{2p}} y_i^\ast \cdot [z_i d^{\frac{1}{2p}}],$$ whereas the first inclusion follows from \cite[Lemma 22]{TerpII}. Because the right hand side is in $L_p(\cM)$ this inclusion is in fact an equality by \cite[Theorem 4 (1)]{Hilsum}. \fin

\begin{lemma}\label{Lem=LpCommutationII}
For every $x \in \mathcal{S}_\varphi$ we have $$j_{p,\rtimes}(x) = D^{\frac{1}{p}} \sigma^\varphi_{\frac{i}{2p}}(x) = [\sigma^\varphi_{-\frac{i}{2p}}(x) D^{\frac{1}{p}}].$$
\end{lemma}

\dem The proof is the same as of Lemma \ref{Lem=LpCommutation}. The only difference being that every time that we used \cite[Theorem 4 (1)]{Hilsum} one uses \cite[Proposition I.12]{TerpI}. \fin

\begin{proposition}
Let $T: \cM \rightarrow \cM$ be a completely bounded map with $\varphi \circ T \le \varphi$ that commutes with $\sigma^\varphi$. Let $T_{p, \rtimes}: L_p(\cM)_\rtimes \rightarrow L_p(\cM)_\rtimes$ be the extended map to the Haagerup $L_p$-space given in \cite[Remark 5.6]{HJX} and which is determined by the relation below for $x \in \mathcal{S}_\varphi$ $$T_{p, \rtimes}: j_{p, \rtimes}(x)\mapsto j_{p, \rtimes}(T(x)).$$ Then, the isometric isomorphism $$\kappa_p: L_p(\cM) \rightarrow L_p(\cM)_{\rtimes}$$ defined in \cite{TerpI} intertwines $T_p$ and $T_{p, \rtimes}$, where $$T_{p}: L_p(\cM) \ni j_p(x) \mapsto j_p(T(x)) \in L_p(\cM) \quad \mbox{for} \quad x \in \mathcal{S}_\varphi.$$
\end{proposition}

\dem Note that the statement above uses that  $T$ preserves the set $\mathcal{S}_\varphi$, which is clear from the definition. Let $u_0: L_2(\mathbb{R}, \cH) \rightarrow  L_2(\mathbb{R}, \cH)$ be the map defined by $(u_0\xi)(s) = d^{is} \xi(s)$ with $s \in \mathbb{R}$. Let $D_0$ be such that $D_0^{is} = \lambda(s)$, where $$(\lambda(s)\xi)(t) =  \xi(t-s)$$ is the left regular representation on $L_2(\R)$. It is proved in \cite[Proposition IV.3]{TerpI} that 
\begin{eqnarray*}
u_0 \pi(x) u_0^\ast & = & x \otimes \1, \\ u_0 \lambda(s) u_0^\ast & = & d^{is} \otimes \lambda(s),
\end{eqnarray*}
where $x \in \cM$ and $s \in \mathbb{R}$. Let $\kappa_p: L_p(\cM) \rightarrow L_p(\cM)_{\rtimes}$ be the isometric isomorphism defined in \cite{TerpI} by $a \mapsto u_0^\ast (a \otimes D_0^{1/p}) u_0$. Let $1 \leq p < \infty$ and $x \in \mathcal{S}_\varphi$. Consider the element $d^{1/p} x$ which is in $L_p(\cM)$ by Lemma \ref{Lem=LpCommutation}. Then
\begin{eqnarray*}
T_{p, \rtimes} \circ \kappa_p \, \big( d^{\frac{1}{p}} x \big) & = & T_{p, \rtimes} \Big( u_0^\ast \big( d^{\frac{1}{p}}x \otimes D_0^{\frac{1}{p}} \big) u_0 \Big) \\ [-2pt] & = & T_{p, \rtimes} \Big( u_0^\ast \big( d^{\frac{1}{p}} \otimes D_0^{\frac{1}{p}} \big) u_0 u_0^\ast \big( x \otimes \1 \big) u_0 \Big) \\ & = & T_{p, \rtimes} \big( D^{\frac{1}{p}} \pi(x) \big) \ = \ D^{\frac{1}{p}} \pi \big( T( x) \big) \\ & = & u_0^\ast \big( d^{\frac{1}{p}} \otimes D_0^{\frac{1}{p}} \big) u_0 u_0^\ast \big( T(x) \otimes \1 \big) u_0 \\ & = & u_0^\ast \big( d^{\frac{1}{p}} T(x) \otimes D_0^{\frac{1}{p}} \big) u_0 \ = \ \kappa_p \circ T_p \big( d^{\frac{1}{p}} x \big).
\end{eqnarray*}
Note in particular that at each instance we have an equality of domains and the fourth equality follows from Lemma \ref{Lem=LpCommutationII} and the definition of $T$. Similarly, the last equality follows from Lemma \ref{Lem=LpCommutation}. Since such elements $ d^{1/p} x$ are dense in $L_p(\cM)$ this proves that $T_{p, \rtimes} \circ \kappa_p = \kappa_p \circ T_p$. This completes the proof. \fin

\begin{remark} \label{RCHHLp}
\emph{In particular, Corollary \ref{Ty-y arbitrary vNA} is valid for Connes-Hilsum $L_p$-spaces.}
\end{remark}

\subsection{Nonunimodular restriction theorem}

We finish this section by sketching the proof of the restriction theorem in the nonunimodular setting, enlightening the main changes. Note that in the nonunimodular case, the Fourier multipliers depend on $p$. However, for the sake of clarity we just used the notation $T_m$ in the statement of Theorem A given in the Introduction. After the proof, we shall construct some natural examples illustrating Theorem A which complement what we did in Section \ref{Sect=Restriction}. We shall also give a brief discussion on Theorem D in Remark \ref{nonunimodularD}. 

\demAA The proof follows the same strategy as in the unimodular case, the main ingredient being that in this case the operator $h_j$ should be defined as $$h_j \, = \, \big\| \ind_{V_j} \Delta_{\G}^{-\frac{1}{4}} \big\|_{L_2(\G)}^{-1} \big[ \lambda( \ind_{V_j} \Delta_{\G}^{-\frac{1}{4}}) \Delta_{\G}^{\frac{1}{2}} \big] \in L_2(\widehat{\G}).$$ Note that $h_j$ is a self-adjoint operator. Indeed, according to \cite[Lemma 2.5]{GoldsteinLindsay} and the fact that $V_j$ is symmetric (recalling that $\xi^{\ast}(g) = \Delta_\G(g)^{-1}\overline{\xi(g^{-1})}$), we obtain the following identity
\begin{eqnarray*}
h_j^\ast & = & \big\| \ind_{V_j} \Delta_G^{-\frac{1}{4}} \big\|_{L_2(\G)}^{-1} \Delta_{\G}^{\frac{1}{2}} \lambda(\ind_{V_j} \Delta_\G^{-\frac{1}{4}})^\ast \\ & = & \big\| \ind_{V_j} \Delta_G^{-\frac{1}{4}} \big\|_{L_2(\G)}^{-1} \Delta_{\G}^{\frac{1}{2}} \lambda(\ind_{V_j} \Delta_\G^{-\frac{3}{4}}) \\ & = & \big\| \ind_{V_j} \Delta_G^{-\frac{1}{4}} \big\|_{L_2(\G)}^{-1} \big[ \lambda(\ind_{V_j}) \Delta_\G^{-\frac{1}{4}}) \Delta_{\G}^{\frac{1}{2}} \big] \ = \ h_j.
\end{eqnarray*}
Then one defines again $\Phi_j^q: L_q(\widehat{\Gamma}) \ni f \mapsto  f u_j \vert h_j \vert^{\frac{2}{q}} \in L_q(\widehat{\G})$, where $h_j = u_j \vert h_j\vert$ is the polar decomposition. The proof proceeds then exactly as in the unimodular case in Section \ref{Sect=Restriction}, which relies on two claims. We can check that Claim A and Claim B can be proved mutatis mutandis, except that we need Corollary \ref{Ty-y} for noncommutative $L_p$-spaces associated with type III algebras. Recall that we applied this result to the Fourier multiplier $T_\zeta$ for some unital, continuous, positive definite function $\zeta \in L_\infty(\G)$. Then $T_\zeta$ is a unital completely positive, $\varphi$-preserving map and by \cite[Example 5.9]{HJX} it commutes with the modular automorphism group. Thus, we may apply Corollary \ref{Ty-y arbitrary vNA} in conjunction with Remark \ref{RCHHLp}. \fin

We now illustrate Theorem A with some example, the main of which will be to show that we can apply restriction to any ADS amenable subgroup for which the modular function restricts properly. We need a preliminary technical result. 

\begin{lemma} \label{Lem=FactAboutAnL1Estimate}
Let $\G$ be a locally compact group. Let $\varepsilon>0$ and $\xi, \eta_1, \cdots, \eta_n$ be positive functions in $L_1(\G)$ satisfying $\sum_{\ell=1}^n \| \xi - \eta_\ell \|_1 < \varepsilon$ and $\Vert \xi \Vert_1 = 1$. Then there exists $t > 0$ such that $$\sum_{\ell=1}^n \big\| \ind_{\{ \xi > t \}} - \ind_{\{\eta_\ell > t\}} \big\|_{L_1(\G)} \, < \, \varepsilon \| \ind_{\{\xi > t\}} \|_{L_1(\G)}.$$
\end{lemma}
 
\dem Given $g \in \G$ and $1\leq \ell\leq n$, we have
\begin{eqnarray*}
\vert \xi(g) \vert & = & \int_0^\infty  \ind_{\{\xi(g) > t\}}dt, \\ \vert \xi(g) - \eta_\ell(g) \vert & = & \int_0^\infty \big| \ind_{\{\xi(g) > t\}} - \ind_{\{\eta_\ell(g) > t\}} \big| \, dt .
\end{eqnarray*}
Hence the hypothesis can be written as follows
\begin{eqnarray*}
\int_0^\infty \sum_{\ell=1}^n \Vert \ind_{\{ \xi > t \} } - \ind_{\{\eta_\ell > t \}} \Vert_{L_1(\G)} dt & = & \sum_{\ell=1}^n \Vert \xi - \eta_\ell \Vert_{L_1(\G)} \\ & <  & \varepsilon \|\xi\|_{L_1(\G)} \ = \ \varepsilon \int_0^\infty \Vert \ind_{\{ \xi > t \}} \Vert_{L_1(\G)} \, dt.
\end{eqnarray*}
This immediately implies the existence of some $t>0$ satisfying the assertion. \fin

\begin{theorem} \label{AmenableImpliesSAIN}
We have $$\G \in [\mathrm{SAIN}]_{\Gamma}$$ for any discrete amenable subgroup $\Gamma$ satisfying that $\Delta_{\G_{\mid_\Gamma}}=\Delta_\Gamma = 1$.
\end{theorem}

\dem Fix a finite set $\mathrm{F}\subset \Gamma$. Since $\Gamma$ is amenable and discrete, we know from the F\o{}lner condition (see Lemma \ref{amenable}) that for any $j\geq 1$ there exist a finite subset $U_{\mathrm{F}, j} \subset \Gamma$ such that $$\frac{\vert U_{\mathrm{F},j} \gamma \triangle U_{\mathrm{F},j} \vert}{\vert U_{\mathrm{F},j} \vert} \, < \, \frac{1}{j|\mathrm{F}|} \quad \mbox{for any } \quad \gamma \in \mathrm{F}.$$ Let $( Z_j )_{j \geq 1}$ be a basis of symmetric neighborhoods of $e$ such that $\mu(Z_j) < \infty$. Since the sets $(U_{\mathrm{F},j})_{j\geq 1}$ are finite, by continuity of the multiplication on $\G$ we can find a sequence $( W_j )_{j \geq 1}$ of symmetric neighborhoods of the identity such that $\bigcup_{g \in U_{\mathrm{F},j}} g^{-1} W_j g \subset Z_j$. Define for each $j \geq 1$ $$\xi_j \, = \, \frac{1}{ \vert U_{\mathrm{F}, j} \vert \mu(W_j)} \sum_{g \in U_{\mathrm{F},j} } \ind_{g^{-1} W_j g}.$$ By unimodularity of $\Gamma$ (since it is discrete), we have $\Vert \xi_j \Vert_{L_1(\G)}= 1$ and we can prove $$\big\| \xi_j( \gamma \cdot \gamma^{-1} ) - \xi_j \big\|_{L_1(\G)} \, < \, \frac{1}{j\vert \mathrm{F} \vert} \quad \mbox{ for any } \gamma \in \mathrm{F}.$$ Indeed, using that for any $j\geq 1$ we have $$\xi_j(\gamma \cdot \gamma^{-1}) - \xi_j \, = \, \frac{1}{\vert U_{\mathrm{F},j}\vert \mu(W_j)} \Big( \sum_{g \in U_{\mathrm{F},j} \gamma \setminus U_{\mathrm{F},j}}  \ind_{ g^{-1} W_j g } - \sum_{g \in U_{\mathrm{F},j}\setminus U_{\mathrm{F},j} \gamma} \ind_{g^{-1} W_j g}\Big),$$ and by construction of the F\o{}lner sets $(U_{\mathrm{F} ,j})_{j\geq 1}$ we get $$\big\| \xi_j(\gamma \cdot \gamma^{-1}) - \xi_j \big\|_{L_1(\G)} \, \le \, \frac{\vert U_{\mathrm{F},j} \gamma \triangle U_{\mathrm{F},j}\vert }{ \vert U_{\mathrm{F},j} \vert} \, < \, \frac{1}{j \vert \mathrm{F} \vert}.$$ Here we also used that $U_{\mathrm{F},j} \gamma \cup U_{\mathrm{F},j} \subset \Gamma$ for any $\gamma \in \mathrm{F}$ and the unimodularity of $\Gamma$. Hence, by applying Lemma \ref{Lem=FactAboutAnL1Estimate}, for each $j\geq 1$ we can find $t_j > 0$ such that the set $V_j = \{ \xi_j > t_j \}$ satisfies
\begin{equation}\label{Eqn=SAINestimate}
\sum_{\gamma \in \mathrm{F}} \frac{\mu \big( \gamma^{-1} V_j \gamma \, \triangle \, V_j \big)}{\mu(V_j)} \, = \, \sum_{\gamma \in \mathrm{F}} \frac{\Vert \ind_{V_j } - \ind_{\gamma^{-1} V_j \gamma} \Vert_{L_1 (\G)}}{\Vert \ind_{V_j} \Vert_{L_1 (\G)}} \, < \, \frac{1}{j}.
\end{equation}
It remains to check that $(V_j)_{j\geq 1}$ is a basis of symmetric neighborhoods of the identity. Since $W_j$ is symmetric, we have $\xi_j(g^{-1})=\xi_j(g)$ for any $g\in \G$ and $V_j$ is clearly symmetric. On the other hand, note that $\|\xi_j\|_{\infty}=\xi_j(e) = \mu(V_j)^{-1}$. Thus $\xi_j(e)>t_j$, otherwise  we would have $\ind_{V_j}=0$, which contradicts \eqref{Eqn=SAINestimate}. Finally, the inclusions $$V_j \subset \supp(\xi_j) \subset \bigcup_{g \in U_{\mathrm{F},j}}g^{-1} W_j g \subset Z_j$$ ensure that  $(V_j)_{j\geq 1}$ is a basis of neighborhoods of the identity. \fin

\begin{corollary} \label{example-restriction}
We have $$\big\| T_{m_{\mid_\H}}: L_p(\widehat{\H}) \to L_p(\widehat{\H}) \big\| \, \le \, \big\| T_m: L_p(\widehat{\G}) \to L_p(\widehat{\G}) \big\|$$ for any $\mathrm{ADS}$ amenable subgroup $\mathrm{H}$ satisfying the identity $\Delta_\mathrm{H} = \Delta_{\G_{\mid_\mathrm{H}}}$.
\end{corollary}

\dem According to Theorem C, we have $$\big\| T_{m_{\mid_\H}}: L_p(\widehat{\H}) \to L_p(\widehat{\H}) \big\| \, \le \, \sup_{j \ge 1} \big\| T_{m_{\mid_{\Gamma_j}}}: L_p(\widehat{\Gamma}_j) \to L_p(\widehat{\Gamma}_j) \big\|$$ for any family $(\Gamma_j)_{j \ge 1}$ of discrete subgroups approximating $\H$. By amenability and unimodularity of $\H$, it is easily seen that each $\Gamma_j$ satisfies the hypothesis of Theorem \ref{AmenableImpliesSAIN}, from which the assertion follows. This completes the proof. \fin

Beyond discrete amenable subgroups of unimodular groups, other pairs $(\G,\H)$ satisfying Corollary \ref{example-restriction} are given by $\G$ unimodular and $\H$ belonging to the families in Remark \ref{examples-ADS}. Corollary \ref{example-restriction} also admits pairs with $\G$ nonunimodular, consider for instance the affine group $\G = \mathbb{R}^n \rtimes {\rm GL}_n(\mathbb{R})$ which is nonunimodular \cite{F}. However $\Delta_\G$ restricts to  ${\rm SL}_n(\mathbb{R})$ (which is unimodular) trivially and hence also to every ADS subgroup. In particular ADS subgroups of ${\rm O}_n(\mathbb{R})$ will form examples of subgroups of $\G$ that satisfy the criteria of Theorem A.

\begin{remark} \label{nonunimodularD}

\emph{Outside the cb-setting (Section \ref{cb}), the compactification Theorem D i) and ii) does not have a suitable analogue in the nonunimodular setting (at least not from our techniques) since we require that $\Delta_\G = \Delta_{\G_{\mathrm{disc}}} \equiv 1$. Similarly the periodization Theorem D iii) is meaningless, since we showed that commutativity of $\G$ (hence unimodularity) is an essential assumption. However, Theorem D iv) does generalize to nonunimodular groups provided that $\Delta_{\G_{\mid_\H}} = \Delta_\H$, the proof is analogous to the one we gave for unimodular groups.}
\end{remark}

\section{{\bf Operator space results}}\label{cb}

The goal of this section is to study de Leeuw's theorems for locally compact groups in the category of operator spaces. More precisely, we are interested in restriction, compactification and periodization results under the assumption that our multipliers are not only bounded, but completely bounded when we equip our $L_p$ spaces with their natural operator space structure \cite{Pis,Pis2}. Then we aim to show that the conclusions also give cb-bounded multipliers. It is easily seen that this is the case when we keep the hypotheses of Theorems A, C and D. In other words, we have for $1 \le p \le \infty$:
\begin{itemize}
\item If $\H \in \mathrm{ADS}$ and $\G \in [\mathrm{SAIN}]_\H$, we have $$\hskip-15pt \big\| T_{m_{\mid_\H}}: L_p(\widehat{\H}) \to L_p(\widehat{\H}) \big\|_{\mathrm{cb}} \, \le \, \big\| T_m: L_p(\widehat{\G}) \to L_p(\widehat{\G}) \big\|_{\mathrm{cb}}$$ for bounded continuous symbols $m: \G \to \C$ provided $\Delta_{\G_{\mid_\H}} = \Delta_\H$.

\vskip3pt

\item If $\G \in \mathrm{ADS}$ is approximated by $(\Gamma_j)_{j \ge 1}$ $$\hskip8pt \big\| T_m: L_p(\widehat{\G}) \to L_p(\widehat{\G}) \big\|_{\mathrm{cb}} \, \le \, \sup_{j \ge 1} \big\| T_{m_{\mid_{\Gamma_j}}}: L_p(\widehat{\Gamma}_j) \to L_p(\widehat{\Gamma}_j) \big\|_{\mathrm{cb}}$$ for bounded $m: \G \to \C$ which are continuous $\mu_\G$--almost everywhere.

\vskip3pt

\item If $\G$ is $\mathrm{ADS}$, $\Gd$ is amenable and $m$ continuous $$\hskip5pt \big\| T_m: L_p(\widehat{\G}) \to L_p(\widehat{\G}) \big\|_{\mathrm{cb}} \, = \, \big\| T_m: L_p(\widehat{\Gd}) \to L_p(\widehat{\Gd}) \big\|_{\mathrm{cb}}.$$ The $\le$ holds for $\G \in \mathrm{ADS}$, the $\ge$ for $\G_{\mathrm{disc}}$ amenable and $\G$ unimodular.

\vskip3pt

\item If $\H \vartriangleleft \G$ is compact and $m_\pi(g) = m_q(g \H)$ is bounded $$\hskip10pt \big\| T_{m_\pi}: L_p(\widehat{\G}) \to L_p(\widehat{\G}) \big\|_{\mathrm{cb}} \, \ge \, \big\| T_{m_q}: L_p(\widehat{\G/\H}) \to L_p(\widehat{\G/\H}) \big\|_{\mathrm{cb}}.$$
\end{itemize}
    
Indeed, except for Theorem D iii) our results remain valid when we apply them to the cartesian product of $\G$ with any finite group, since our ADS and SAIN assumptions are stable under that operation. This operation allows to generalize our results to the cb-setting in a trivial way. 

\begin{remark}
\emph{The upper estimate $\le$ in our cb-periodization result can be extended to any pair $(\G,\H)$ as long as $\G$ is discrete, $\V$ is QWEP and $\Delta_{\G} = \Delta_\H$ on $\H$. Indeed, the discreteness of $\G$ and $\G/\H$ allows us to apply Fell's absorption principle in Lemma \ref{Fell} ii) to the strongly continuous representation $\pi: g \mapsto \lambda_{\G/\H}(g\H)$ and the existence of an invariant measure is then used to factorize the integral over $\G$ as an integral over $\G/\H \times \H$. After rearrangement and Fubini's theorem (for which we use the QWEP property following \cite{J-Fubini}) one concludes.}
\end{remark}    

Motivated by the transference results from \cite{BF,CdlS,NR} between Fourier and Schur multipliers, an alternative approach to obtain de Leeuw type theorems is to exploit that such results are much more elementary for Schur multipliers. Namely, given a bounded symbol $m: \G \to \C$, recall that the associated Herz-Schur multiplier is formally defined as the linear map $$S_m: \sum_{g_1, g_2 \in \G} a_{g_1, g_2} e_{g_1, g_2} \mapsto \sum_{g_1, g_2 \in \G} m(g_1^{-1} g_2) a_{g_1, g_2} e_{g_1, g_2}.$$ By the boundedness of $m$, it is clear that $S_m$ is (completely) bounded on the Schatten class $S_2(L_2(\G))$. When it maps $S_2(L_2(\G)) \cap S_p(L_2(\G))$ to $S_p(L_2(\G))$ and extends to a cb-map on $S_p(L_2(\G))$, we say that $S_m$ is a cb-bounded Schur multiplier on $S_p(L_2(\G))$. Let us analyze de Leeuw operations for Schur multipliers.

\begin{lemma} \label{lemmaSchur}
If $1\leq p \leq \infty$ and $m: \G \to \C$ is continuous 
\begin{eqnarray*}
\lefteqn{\hskip-52pt \big\| S_m: S_p(L_2(\G)) \to S_p(L_2(\G)) \big\|_{\mathrm{cb}}} \\ \hskip20pt & = & \big\| S_m: S_p(\ell_2(\G_{\mathrm{disc}})) \to S_p(\ell_2(\G_{\mathrm{disc}})) \big\|_{\mathrm{cb}}.
\end{eqnarray*}
Moreover, let $\H$ be a closed subgroup of $\G$. Then we additionally have 
\begin{itemize}
\item[i)] If $m: \G \to \C$ is continuous $$\hskip30pt \big\| S_{m_{\mid_\H}}: S_p(L_2(\H)) \to S_p(L_2(\H)) \big\|_{\mathrm{cb}} \le \big\| S_m: S_p(L_2(\G)) \to S_p(L_2(\G)) \big\|_{\mathrm{cb}}.$$

\item[ii)] If $\H \vartriangleleft \G$ and $m_q: \G/\H \to \C$ is continuous $$\hskip18pt \big\| S_{m_\pi}: S_p(L_2(\G)) \to S_p(L_2(\G)) \big\|_{\mathrm{cb}} = \big\| S_{m_q}: S_p(L_2(\G/\H)) \to S_p(L_2(\G/\H)) \big\|_{\mathrm{cb}}.$$

\end{itemize}
\end{lemma} 

\dem Lafforgue and de la Salle established in \cite[Theorem 1.19]{LdlS} (extending an unpublished result of Haagerup in the $L_\infty$-case) that for any locally compact group $\G$ and any continuous symbol $m:\G\to \C$, the cb-norm of the Schur multiplier is given by
\begin{eqnarray} \label{cb-norm Schur}
\lefteqn{\hskip-60pt \big\|S_m: S_p(L_2(\G)) \to S_p(L_2(\G)) \big\|_{\mathrm{cb}}} \\ \nonumber \hskip30pt & = & \sup_{\begin{subarray}{c} \mathrm{F} \subset \G \\ \mathrm{F} \ \mathrm{finite} \end{subarray}} \big\| S_{m_{\mid_{\mathrm{F}}}}:S_p(\ell_2(\mathrm{F})) \to S_p(\ell_2(\mathrm{F}))\big\|_{\mathrm{cb}}.
\end{eqnarray}
The first assertion (compactification) and property i) (restriction) follow directly from this. The cb-periodization ii) for Schur multipliers can also be deduced from \eqref{cb-norm Schur} as follows. For a fixed fundamental domain $\mathrm{X}$, we consider the natural map $\sigma:\G/\H \to \mathrm{X}$. Then we may identify the group $\G$ with the cartesian product $\G/\H \times \H$ as in the proof of Lemma \ref{LHsubalgLG} via the bijective map $$\Upsilon: \G \ni g \mapsto (g\H,h(g)) \in \G/\H \times \H$$ where $g=\sigma(g\H)h(g)$. For $1\leq p \leq \infty$, this gives a map $$\Upsilon:S_p(L_2(\G))\to S_p(L_2(\G/\H) \otimes L_2(\H))$$ which is completely isometric on finite subsets. Moreover, this map intertwines the Schur multipliers $\Upsilon \circ S_{m_\pi} \circ \Upsilon^{-1}=S_{m_q} \otimes id_{\mathcal{B}(L_2(\H))}$. Therefore, by \eqref{cb-norm Schur} we can write
\begin{eqnarray*}
\lefteqn{\big\|S_{m_q}:S_p(L_2(\G/\H))\to S_p(L_2(\G/\H))\big\|_{\mathrm{cb}}} \\ & = & \sup_{n \ge 1} \sup_{\begin{subarray}{c} (\mathrm{F}_1, \mathrm{F}_2) \in \G/\H \times \H \\ \mathrm{F}_1, \mathrm{F}_2 \ \mathrm{finite} \end{subarray}} \sup_{\begin{subarray}{c} \|\mathrm{A}\|_{S_p} \le 1 \\ \mathrm{A} \in \mathbb{M}_{|\mathrm{F}_1| |\mathrm{F}_2| n} \end{subarray}} \big\|S_{m_q}\otimes id_{\mathbb{M}_{|F_2|}}\otimes id_{\mathbb{M}_n}(\mathrm{A}) \big\|_{S_p(\ell_2(\mathrm{F}_1) \otimes \ell_2(\mathrm{F}_2) \otimes \ell_2^n)}.
\end{eqnarray*} 
On the other hand, each term of this supremum satisfies
\begin{eqnarray*}
\lefteqn{\hskip-15pt \big\| S_{m_q} \otimes id_{\mathbb{M}_{|F_2|}} \otimes id_{\mathbb{M}_n} (\mathrm{A}) \big\|_{S_p(\ell_2(\mathrm{F}_1)\otimes \ell_2(\mathrm{F}_2)\otimes \ell_2^n)}} \\ [3pt] & = & \big\|( \Upsilon \otimes id_{\mathbb{M}_n}) (S_{m_\pi} \otimes id_{\mathbb{M}_n})(\Upsilon^{-1} \otimes id_{\mathbb{M}_n})(\mathrm{A}) \big\|_{S_p(\ell_2(\mathrm{F}_1) \otimes \ell_2(\mathrm{F}_2) \otimes \ell_2^n)} \\ [3pt] & = &  \big\|(S_{m_\pi} \otimes id_{\mathbb{M}_n})(\widetilde{\mathrm{A}}) \big\|_{S_p(\ell_2(\mathrm{F}) \otimes \ell_2^n)},
\end{eqnarray*}
where $\mathrm{F} = \Upsilon^{-1}(\mathrm{F}_1 \times \mathrm{F}_2) \subset \G$ is finite and $$\widetilde{\mathrm{A}} = \Upsilon^{-1} \otimes id_{\mathbb{M}_n}(\mathrm{A}) \in S_p(\ell_2(\mathrm{F}) \otimes \ell_2^n)$$ is of norm $1$. Hence we deduce that $$\big\|S_{m_q}:S_p(L_2(\G/\H))\to S_p(L_2(\G/\H))\big\|_{\mathrm{cb}} =\big\|S_{m_\pi}:S_p(L_2(\G))\to S_p(L_2(\G))\big\|_{\mathrm{cb}}.$$ Indeed, the left hand side is clearly dominated by the right hand side. The lower estimate also holds since for any finite subset $\mathrm{F} \subset \G$, we can find finite subsets $\mathrm{F}_1\subset \G/\H$ and $\mathrm{F}_2\subset \H$ such that $$\mathrm{F} \subset \Upsilon^{-1}(\mathrm{F}_1 \times \mathrm{F}_2).$$ Thus, the result follows using that the cb-norm in \eqref{cb-norm Schur} is increasing with $\mathrm{F}$. \fin

This shows that de Leeuw theorems extend in almost full generality to the context of Schur multipliers, only continuity of the symbols is needed. In particular, we do not impose any of our former conditions like ADS, SAIN, the compatibility of modular functions or the amenability of $\G_{\mathrm{disc}}$. We now want to use certain transference results to obtain de Leeuw type theorems for Fourier multipliers from the results in Lemma \ref{lemmaSchur}. More precisely, we will use that we have 
\begin{equation}\label{Schur-Fourier}
\big\| T_m: L_p(\widehat{\G}) \to L_p(\widehat{\G}) \big\|_{\mathrm{cb}} \, = \, \big\| S_m: S_p(L_2(\G)) \to S_p(L_2(\G)) \big\|_{\mathrm{cb}}
\end{equation} 
for $1 \le p \le \infty$ under the following conditions
\begin{itemize}
\item[(i)] $\G$ is an amenable group,
\item[(ii)] $m\in L_\infty(\G)$ defines a completely bounded Fourier multiplier on $L_p(\widehat{\G})$.
\end{itemize}
When $p = 1,\infty$ this was proved by Bo\.zejko and Fendler \cite{BF}. Other values of $p$ were first considered by Neuwirth and Ricard \cite{NR}, who proved \eqref{Schur-Fourier} for amenable discrete groups. 
Caspers and de la Salle \cite{CdlS} then obtained this result for arbitrary amenable groups and $1 < p < \infty$. We shall need this identity to transfer Lemma \ref{lemmaSchur} to Fourier multipliers. Hence the price to avoid our conditions listed at the beginning of this section is to assume amenability of $\G$. 

\begin{remark}\label{Schur-Fourier-approx}
{\rm Observe that the transference theorem proved in \cite{CdlS} requires the extra assumption that the symbol $m:\G\to \C$ gives rise to a completely bounded Fourier multiplier on $\cL \G$. The set of such symbols is denoted by $M_{\mathrm{cb}}(\G)$. By approximation we may extend the identity \eqref{Schur-Fourier} to any bounded symbol $m:\G\to \C$ satisfying the above condition (ii) whenever $\G$ is amenable. Indeed, consider a symbol $m\in L_\infty(\G)$ verifying (ii). Notice that when $\G$ is amenable there is a continuous contractive approximate unit $(m_i)_{i\geq 1}$ in the Fourier algebra $A(\G)$ with compact support. Take also $(\chi_j)_{j\geq 1}$ a contractive approximate unit in $L_1(\G)$ that also belongs to $L_2(\G)$. Define $$m_{i,j}=\chi_j*(m_i m)\in L_\infty(\G).$$ Clearly $m_{i,j}\in A(\G)$ and hence lies in $M_{cb}(\G)$. On the other hand, one can check that
\begin{eqnarray*} 
\big\| T_m:L_p(\widehat{\G}) \to L_p(\widehat{\G}) \big\|_{\mathrm{cb}} & = & \lim_{i,j} \big\| T_{m_{i,j}}: L_p(\widehat{\G}) \to L_p(\widehat{\G}) \big\|_{\mathrm{cb}}, \\ \big\| S_{m}:S_p(L_2(\G)) \to S_p(L_2(\G)) \big\|_{\mathrm{cb}} & = & \lim_{i,j}\big\| S_{m_{i,j}}: S_p(L_2(\G)) \to S_p(L_2(\G)) \big\|_{\mathrm{cb}}.
\end{eqnarray*}
Indeed, the lower estimates easily follows from standard properties of Fourier and Schur multipliers, and we may deduce the upper estimates from the fact that $T_{m_{i,j}}\to T_{m}$ (resp. $S_{m_{i,j}}\to S_{m}$) pointwise in the weak-topology of  $L_p(\widehat{\G})$ (resp. $S_p(L_2(\G))$). Using Caspers and de la Salle's result for the symbols $m_{i,j}$ in $M_{\mathrm{cb}}(\G)$, this allows us to conclude that \eqref{Schur-Fourier} holds true for the symbol $m$.}
\end{remark}
 
\begin{theorem}\label{de Leeuw-cb}
Let $1 \le p \le \infty$ and $\G$ amenable\emph{:} 
\begin{itemize}
\item[i)] If $m: \G \to \C$ is bounded and continuous and $\H$ is a closed subgroup of $\G$ $$\hskip-20pt \big\| T_{m_{\mid_\H}}: L_p(\widehat{\H}) \to L_p(\widehat{\H}) \big\|_{\mathrm{cb}} \, \le \, \big\| T_m: L_p(\widehat{\G}) \to L_p(\widehat{\G}) \big\|_{\mathrm{cb}}.$$

\item[ii)] If $m:  \G \to \C$ is bounded and continuous and $\Gd$ is amenable, we have $$\big\| T_m: L_p(\widehat{\G}) \to L_p(\widehat{\G}) \big\|_{\mathrm{cb}} \, = \, \big\| T_m: L_p(\widehat{\G_{\mathrm{disc}}}) \to L_p(\widehat{\G_{\mathrm{disc}}}) \big\|_{\mathrm{cb}}.$$

\item[iii)] If $m_q: \G/\H \to \C$ is bounded and continuous and $\H$ is a normal closed subgroup of $\G$ $$\hskip2pt \big\| T_{m_\pi}: L_p(\widehat{\G}) \to L_p(\widehat{\G}) \big\|_{\mathrm{cb}} = \big\| T_{m_q}: L_p(\widehat{\G/\H}) \to L_p(\widehat{\G/\H}) \big\|_{\mathrm{cb}}.$$
\end{itemize}
\end{theorem} 

\dem It follows from Lemma \ref{lemmaSchur}, the transference theorem \eqref{Schur-Fourier} from \cite{CdlS} and Remark \ref{Schur-Fourier-approx}. \fin 


\begin{remark}
{\rm Recall that the lattice approximation Theorem C only works in the unimodular setting (since we need to assume $\G \in \mathrm{ADS}$), hence applying the transference in that case would not improve the cb-result obtained directly from Theorem C. 
In fact, applying the transference theorem from \cite{CdlS} and Remark \ref{Schur-Fourier-approx} in conjunction with Theorem \ref{de Leeuw-cb} i) to that result, we deduce the analog for Schur multipliers. Namely, for any group $\G\in\mathrm{ADS}$ approximated by $(\Gamma_j)_{j\geq 1}$, $1\leq p\leq \infty$ and any bounded a.e. continuous symbol $m:\G \to \C$, we have $$ \big\| S_m: S_p(L_2(\G)) \to S_p(L_2(\G)) \big\|_{\mathrm{cb}} \, = \, \sup_{j \ge 1} \big\| S_{m_{\mid_{\Gamma_j}}}: S_p(\ell_2(\Gamma_j)) \to S_p(\ell_2(\Gamma_j)) \big\|_{\mathrm{cb}}.$$}
\end{remark}

\section*{{A. \hskip4pt \bf Idempotent multipliers in $\R$}}

Idempotent Fourier multipliers are those whose symbols are the characteristic functions of a measurable set $\Sigma$. Intervals in $\R$ or polyhedrons in $\R^n$ are examples of idempotent symbols which yield $L_p$-bounded Fourier multipliers $(1 < p < \infty)$ as a consequence of the boundedness of the Hilbert transform. When $n > 1$, we know from the work of Fefferman \cite{Fef} a fundamental restriction for $L_p$-boundedness of idempotent Fourier multipliers over (say) convex sets $\Sigma$ with boundary $\partial \Sigma$. Namely, let $$\partial \Sigma^\perp \, = \, \Big\{ v \in \mathbb{S}^{n-1} \, \big| \, v \perp \partial \Sigma \Big\}.$$ Then, given $\Pi \subset \R^n$ any $2$-dimensional vector space, $\Omega = \partial \Sigma^\perp \cap \Pi$ can not admit Kakeya sets of directions in the sense of \cite{Fef} or \cite[Lemma 10.1.1]{Gra} when $\Sigma$ leads to an $L_p$-bounded idempotent multiplier. To be more precise we need a bit of terminology. Given a rectangle $R$ in $\R^2$, denote by $R'$ one of the two translations of $R$ which are adjacent to $R$ along its shortest side. After a careful reading of the argument in \cite{Fef, Gra}, we could say that a subset $\Omega$ of the unit circle in $\R^2$ admits Kakeya sets of directions when for every $\mathrm{N} \ge 1$ there exists a finite collection of pairwise disjoint rectangles $\RR_\Omega(\mathrm{N})$ with longest side pointing in a direction of $\Omega$ and a family $\RR'_\Omega(\mathrm{N})$ formed by rectangles $R'$ adjacent to the members of $\RR_\Omega(\mathrm{N})$ along their shortest side and such that $$\Big| \bigcup_{R \in \RR_\Omega({\mathrm{N}})} R \, \Big| \ \ge \ \mathrm{N} \ \Big| \bigcup_{R' \in \RR'_\Omega({\mathrm{N})}} R' \, \Big|.$$ The above symbol $| \ |$ refers to the Lebesgue measure. This notion is closely related to Bateman's notion of Kakeya sets of directions \cite{Bateman}. Fefferman's theorem implies that $\partial \Sigma$ must have vanishing curvature for $L_p$-boundedness, as for polyhedrons. Other regions with flat boundary ---polytopes with infinitely many faces--- may or may not admit Kakeya sets of directions. This is very connected to the boundedness of directional maximal operators \cite{Bateman,PR} but we shall not analyze these subtleties here. Apart from the geometric aspect of $\Sigma$, one may consider which topological structures of $\Sigma$ yield $L_p$-boundedness. In dimension $1$, Lebedev and Olevskii \cite{LO} showed that $\Sigma$ must be open up to a set of zero measure, see also Mockenhaupt and Ricker \cite{MR} for $L_p$-bounded idempotents which are not $L_q$-bounded.           

Our aim in this Appendix is motivated by a problem left open in \cite{JMP1}. The authors provided there a noncommutative H\"ormander-Mihlin multiplier theorem using group cocycles in discrete groups as substitutes of more standard geometric tools for Lie groups. This gave rise to some exotic Euclidean multipliers which are $L_p$-bounded in $\R^n$. Consider the cocycle $b:\R \to \R^4$ given by $$s \mapsto b(s) = \big( \cos (2\pi s) -1, \sin (2\pi s), \cos(2 \pi \beta s) - 1, \sin (2 \pi \beta s) \big)$$ associated with the action $\alpha:\R \curvearrowright \R^4 \simeq \C^2$ $$\alpha_s(x_1,x_2,x_3,x_4) \, \simeq \, \alpha_s(z_1, z_2) \, = \, \big( e^{2 \pi i s} z_1, e^{2\pi i \beta s} z_2 \big).$$ Then, any symbol of the form $m(s) = \widetilde{m}(b(s))$ satisfying that $$|\partial_s^\beta \widetilde{m}(s)| \, \lesssim \, |s|^{-|\beta|} \quad \mbox{for} \quad s \in \R^4 \setminus \{0\} \quad \mbox{and} \quad 0 \le |\beta| \le 3$$ defines an $L_p$-bounded Fourier multiplier in $\R$ for $1 < p < \infty$. Take for instance $\widetilde{m}$ a H\"ormander-Mihlin smoothing of the characteristic function of an open set $\Sigma$ in $\R^4$ intersecting the range of $b$. If $\beta \in \R \setminus \mathbb{Q}$, the cocycle $b$ has a dense orbit and $m$ oscillates from $0$ to $1$ infinitely often with no periodic pattern. A moment of thought shows that the $L_p$-boundedness of such a multiplier follows from the combination of de Leeuw's restriction and periodization theorems, but this cocycle formulation led Junge, Mei and Parcet to pose a similar problem in \cite{JMP1} when the lifted multiplier $\widetilde{m}$ is not smooth anymore. More precisely, let $\widetilde{m}$ be the characteristic function of certain set $\Sigma$ which yields an $L_p$-bounded multiplier in $\R^4$ and intersects the range of the cocycle $b$. Is $m = \widetilde{m} \circ b$ an $L_p$-bounded idempotent multiplier on $\R$?    

In order to answer the question above, let us formulate the problem in a more transparent way. The image of the cocycle $b$ is an helix in a two-dimensional torus which up to a translation we may identify with $\mathbb{T}^2 \simeq [0,1]^2$. Moreover, under this identification, the helix corresponds to the straight line $\gamma$ in $\R^2$ passing through the origin with slope $\beta$. Let us consider the set $\Omega$ which results of the intersection between $\Sigma$ and the two-dimensional torus where $b$ takes values. We shall identify this set with the corresponding set in $[0,1]^2$, still denoted by $\Omega$. According to the results in \cite{Fef}, we know that $\Sigma$ must have flat boundary. Assume for simplicity that $\Sigma$ is a simple object like a semispace or a convex polyhedron ---finite unions and certain infinite unions of this kind of sets also define $L_p$-bounded idempotent multipliers--- so that $\Omega$ is a closed simply connected set. In summary, given a simply connected set $\Omega$ in $[0,1]^2$ and certain slope $\beta$, we may consider the idempotent Fourier multiplier associated with the symbol determined by Figure I below and given by $$\mathrm{M}_{\Omega,\beta}(s) = \ind_{\Omega}\Big((s,\beta s)+\Z^2\Big) \quad \mbox{for} \quad s \in \R.$$ 

\vskip-5pt

\noindent
\begin{picture}(360,140)(-180,-60)
\linethickness{0.2pt}
    
\qbezier(10,10)(20,20)(15,0)    
\qbezier(5,-10)(0,0)(10,10)   
\qbezier(15,0)(10,-15)(5,-10)    

\qbezier(-30,10)(-20,20)(-25,0)    
\qbezier(-35,-10)(-40,0)(-30,10)   
\qbezier(-25,0)(-30,-15)(-35,-10)    

\qbezier(-70,10)(-60,20)(-65,0)    
\qbezier(-75,-10)(-80,0)(-70,10)   
\qbezier(-65,0)(-70,-15)(-75,-10)    

\qbezier(-110,10)(-100,20)(-105,0)    
\qbezier(-115,-10)(-120,0)(-110,10)   
\qbezier(-105,0)(-110,-15)(-115,-10)    

\qbezier(-110,50)(-100,60)(-105,40)    
\qbezier(-115,30)(-120,40)(-110,50)   
\qbezier(-105,40)(-110,25)(-115,30)  

\qbezier(-70,50)(-60,60)(-65,40)    
\qbezier(-75,30)(-80,40)(-70,50)   
\qbezier(-65,40)(-70,25)(-75,30)  

\qbezier(-30,50)(-20,60)(-25,40)    
\qbezier(-35,30)(-40,40)(-30,50)   
\qbezier(-25,40)(-30,25)(-35,30)  

\qbezier(10,50)(20,60)(15,40)    
\qbezier(5,30)(0,40)(10,50)   
\qbezier(15,40)(10,25)(5,30)  

\qbezier(50,50)(60,60)(55,40)    
\qbezier(45,30)(40,40)(50,50)   
\qbezier(55,40)(50,25)(45,30)  

\qbezier(90,50)(100,60)(95,40)    
\qbezier(85,30)(80,40)(90,50)   
\qbezier(95,40)(90,25)(85,30)  

\qbezier(130,50)(140,60)(135,40)    
\qbezier(125,30)(120,40)(130,50)   
\qbezier(135,40)(130,25)(125,30)  

\qbezier(130,10)(140,20)(135,0)    
\qbezier(125,-10)(120,0)(130,10)   
\qbezier(135,0)(130,-15)(125,-10)  

\qbezier(130,-30)(140,-20)(135,-40)    
\qbezier(125,-50)(120,-40)(130,-30)   
\qbezier(135,-40)(130,-55)(125,-50)  

\qbezier(90,-30)(100,-20)(95,-40)    
\qbezier(85,-50)(80,-40)(90,-30)   
\qbezier(95,-40)(90,-55)(85,-50)  

\qbezier(90,10)(100,20)(95,0)    
\qbezier(85,-10)(80,0)(90,10)   
\qbezier(95,0)(90,-15)(85,-10)  

\qbezier(50,10)(60,20)(55,0)    
\qbezier(45,-10)(40,0)(50,10)   
\qbezier(55,0)(50,-15)(45,-10)  

\qbezier(50,-30)(60,-20)(55,-40)    
\qbezier(45,-50)(40,-40)(50,-30)   
\qbezier(55,-40)(50,-55)(45,-50)  

\qbezier(10,-30)(20,-20)(15,-40)    
\qbezier(5,-50)(0,-40)(10,-30)   
\qbezier(15,-40)(10,-55)(5,-50)  

\qbezier(-30,-30)(-20,-20)(-25,-40)    
\qbezier(-35,-50)(-40,-40)(-30,-30)   
\qbezier(-25,-40)(-30,-55)(-35,-50)  

\qbezier(-70,-30)(-60,-20)(-65,-40)    
\qbezier(-75,-50)(-80,-40)(-70,-30)   
\qbezier(-65,-40)(-70,-55)(-75,-50)  

\qbezier(-110,-30)(-100,-20)(-105,-40)    
\qbezier(-115,-50)(-120,-40)(-110,-30)   
\qbezier(-105,-40)(-110,-55)(-115,-50)  

\qbezier(-140,-60)(0,-5)(150,53.9285)  

\put(-130,-40){$\Omega$}    
\put(-18,-17.5){{\tiny $\mbox{slope$(\gamma)$ } \hskip-4pt = \hskip-2pt \beta$}}    

\put(150,45){$\gamma$}

    \put(140,-60){\line(0,1){130}}
    \put(100,-60){\line(0,1){130}}
    \put(60,-60){\line(0,1){130}}
    \put(-100,-60){\line(0,1){130}}
    \put(-140,60){\line(1,0){290}}
    \put(-140,-60){\vector(0,1){130}}
    \put(-60,-60){\line(0,1){130}}
    \put(-20,-60){\line(0,1){130}}
    \put(-140,-20){\line(1,0){290}}
    \put(-140,20){\line(1,0){290}}
    \put(20,-60){\line(0,1){130}}
    \put(-140,-60){\vector(1,0){290}}

\put(-100,80){\line(-1,-2){5}}
    \put(-100,80){\line(1,-2){5}}
   \put(150,-100){\line(-2,1){10}}
    \put(150,-100){\line(-2,-1){10}}

\linethickness{.6pt}
\qbezier(-114.5,-50.1786)(-111.5,-49.0000)(-108.5,-47.8215)  

\qbezier(-72.8,-33.6000)(-67.5,-31.7142)(-63.5,-30.1428)  

\qbezier(4,-3.4285)(10,-1.0714)(15,0.8928)  

\qbezier(87,29.1785)(88.5,29.7678)(90,30.3571)  

\qbezier(126,44.3071)(131,46.2714)(136,48.2357)  
\end{picture}

\null

\vskip-25pt

\null
\begin{center}
\textsc{Figure I} \\ The idempotent symbol $\mathrm{M}_{\Omega,\beta}$ \\ $\mathrm{M}_{\Omega,\beta} = 1$ when $\gamma$ intersects $\Omega + \Z^2$ and $0$ otherwise
\end{center}

\vskip5pt

Our problem is to decide for which pairs $(\Omega,\beta)$ we get $L_p$-bounded idempotent multipliers on $\R$. There are two cases for which the answer is simple. If the slope $\beta \in \mathbb{Q}$, the helix is periodic and so is $\mathrm{M}_{\Omega,\beta}$. Therefore, the $L_p$-boundedness follows by the boundedness of the Hilbert transform for $1 < p < \infty$ (finitely many times) in conjunction with de Leeuw's periodization in $\R$. On the other hand, we also obtain $L_p$-boundedness when $\Omega$ is a polyhedron (finitely many faces) since we know its characteristic function defines an $L_p$-bounded idempotent multiplier in $\R^2$ (finitely many directional Hilbert transforms). Namely, its $\Z^2$-periodization in $L_p(\R^2)$ and its restriction to $\gamma$ in $L_p(\R)$ are still bounded by de Leeuw's periodization and restriction theorems. In particular, the interesting case arises for sets $\Omega$ admitting Kakeya sets of directions ---either having smooth boundary with non-zero curvature as in Figure I or with infinitely many flat faces admitting Kakeya sets--- and slope $\beta \in \R \setminus \mathbb{Q}$. We will answer this problem in the negative by combining de Leeuw's restriction, lattice approximation and Fefferman's construction. 

\begin{Atheorem}\label{thm:idempotent}
Assume that $$\beta \in \R \setminus \mathbb{Q} \quad \mbox{and} \quad \Omega \subset [0,1]^2 \ \mbox{admits Kakeya sets of directions}.$$ Then $\mathrm{M}_{\Omega,\beta}$ does not give rise to a bounded multiplier in $L_p(\R)$ for $1<p\neq 2 <\infty$.  
\end{Atheorem}

\renewcommand{\theequation}{A.1}
\addtocounter{equation}{-1}

\dem Assume there exists $1 < p_0 \neq 2 < \infty$ such that 
\begin{equation} \label{eqA1}
\big\| T_{\mathrm{M}_{\Omega,\beta}}: L_{p_0}(\R) \to L_{p_0}(\R) \big\| \, \le \, \mathrm{C}_0 \, < \, \infty.
\end{equation}
Let $\mathrm{M}_{\Omega,\beta}(\mathrm{L}) = \mathrm{M}_{\Omega,\beta} \, \ind_{[0,\mathrm{L}]}$ be the $\mathrm{L}$-truncation of our multiplier. If $\mathrm{HT}_{p_0}$ denotes the norm of the Hilbert transform on $L_{p_0}(\R)$, it is clear that we have the following bound $$\sup_{\mathrm{L} > 0} \big\| T_{\mathrm{M}_{\Omega,\beta}(\mathrm{L})}: L_{p_0}(\R) \to L_{p_0}(\R) \big\| \, \le \, 2 \, \mathrm{C}_0 \, \mathrm{HT}_{p_0}.$$ On the other hand, we may consider the polygon $\Pi_\Omega(\mathrm{L}, \beta)$ determined by the crossing points of $\gamma + \Z^2$ with $\Omega$ in $[0, \mathrm{L}]$. More precisely, let us set $$\Pi_\Omega(\mathrm{L}, \beta) \, = \, \mathrm{Conv \Big(\Omega \cap \big\{ (s, \beta s) + \Z^2 \, : \, s \in [0,\mathrm{L}] \big\} \Big)}.$$ It is illustrated in Figure II and $\mathrm{M}_{\Omega,\beta}(\mathrm{L})=\mathrm{M}_{\Pi_\Omega(\mathrm{L}, \beta),\beta}(\mathrm{L})$. By the irrationality of $\beta$, the set $\gamma+\Z^2$ is dense in $[0,1]^2$ and $\Pi_\Omega(\mathrm{L}, \beta)$ converges uniformly to $\Omega$ as $\mathrm{L} \to \infty$. In particular, by constructing finer and finer Kakeya sets of directions, we may pick $\mathrm{L}_0$ large enough so that the following inequality holds 
\renewcommand{\theequation}{A.2}
\addtocounter{equation}{-1}
\begin{equation} \label{eqA2}
\inf_{\mathrm{L} \ge \mathrm{L}_0} \big\| T_{\ind_{\Pi_\Omega(\mathrm{L}, \beta)}}: L_{p_0}(\mathbb{R}^2) \to L_{p_0}(\mathbb{R}^2) \big\| \, > \, 4 \, \mathrm{C}_0 \, \mathrm{HT}_{p_0}.
\end{equation} 
We will complete the proof by showing that \eqref{eqA1} and \eqref{eqA2} produce a contradiction.

\vskip-5pt

\noindent
\begin{picture}(360,140)(-180,-60)
\linethickness{0.2pt}
    
\qbezier(30,30)(60,60)(45,0)    
\qbezier(15,-30)(0,0)(30,30)   
\qbezier(45,0)(30,-45)(15,-30)    

\put(16.5, 5){{\tiny $\Pi_\Omega(\mathrm{L},\beta)$}}

\put(-45,-17.5){{\tiny $\mbox{slope } \beta$}}    

\qbezier(-60,-60)(0,-36.4285)(60,-12.8571)  
\qbezier(-60,-12.8571)(0,10.7143)(60,34.2857)  
\qbezier(-60,34.2857)(-27.2728,47.1428)(5.4545,60)  
\qbezier(5.4545,-60)(32.7272,-49.2857)(60,-38.5714)  
\qbezier(-60,-38.5714)(0,-15)(60,8.5714)  
\qbezier(-60,8.5714)(0,32.1428)(60,55.7142)  
\qbezier(-60,55.7142)(-54.5453,57.8571)(-49.0906,60) 
\qbezier(-49.0906,-60)(5.4547,-38.5714)(60,-17.1429) 
\qbezier(-60,-17.1429)(0,6.4285)(60,29.9999) 
   
    \put(60,-60){\line(0,1){120}}
    \put(-60,-60){\vector(0,1){130}}
    \put(-60,60){\line(1,0){120}}
    \put(-60,-60){\vector(1,0){130}}
    
    \put(-65,-66){{\tiny 1}}
    \put(62,-12){{\tiny 2}}
    \put(-65,-12){{\tiny 2}}
    \put(62,35){{\tiny 3}}
    \put(-65,35){{\tiny 3}}
    \put(3,63){{\tiny 4}}
    \put(3,-66){{\tiny 4}}
    \put(62,-40){{\tiny 5}}
    \put(-65,-40){{\tiny 5}}
    \put(62,10){{\tiny 6}}
   \put(-65,10){{\tiny 6}}
   \put(62,55){{\tiny 7}}
   \put(-65,55){{\tiny 7}}
   \put(-52,63){{\tiny 8}}
   \put(-52,-66){{\tiny 8}}
   \put(62,-20){{\tiny 9}}
   \put(-65,-20){{\tiny 9}}
   \put(62,28){{\tiny 10}}

\linethickness{.4pt}

\qbezier(20.5,18.7677)(35,24.4642)(50,30.3571) 

\qbezier(16.5,12.7677)(18.5,15.7677)(20.5,18.7677) 

\qbezier(16.5,12.7677)(13.25,0.9)(10,-11.1) 

\qbezier(15.5,-30)(12.75,-20.55)(10,-11.1) 

\qbezier(15.5,-30)(17.5,-31.5)(19.5,-33) 

\qbezier(31.3,-28.2)(25.4,-30.6)(19.5,-33) 

\qbezier(31.3,-28.2)(33.55,-25.2)(35.6,-22.2) 

\qbezier(45.5,3)(40.5,-9.6)(35.6,-22.2) 

\qbezier(45.5,3)(47.5,14.25)(49.5,25.5) 

\qbezier(49.8,30)(49.65,27.75)(49.5,25.5) 

\end{picture}

\null

\vskip-15pt

\null
\begin{center}
\textsc{Figure II} \\ The polygon $\Pi_\Omega (\mathrm{L},\beta)$ determined by a $\mathrm{L}$-truncation \\ If $\beta \in \R \setminus \mathbb{Q}$, the polygon $\Pi_\Omega (\mathrm{L},\beta)$ converges to $\Omega$ as $\mathrm{L} \to \infty$ \\ Pick coprimes $p \neq q$ such that $\frac{p}{q} \sim \beta$ \hskip-2pt / \hskip-2pt Dilate $\mathrm{M}_{\Omega,\beta}$ and approximate it by $\Z_{pq}$
\end{center}
\renewcommand{\theequation}{A.3}
\addtocounter{equation}{-1}

\vskip5pt

\noindent Indeed, according to Dirichlet's diophantine approximation, since $\beta$ is irrational we may find infinitely many coprime integers $p,q$ so that $| \beta - p/q | < 1/q^2$. Denote by $\mathcal{I}$ the set of such pairs of coprime integers and pick $(p,q) \in \mathcal{I}$. On the other hand by dilation-invariance of the $L_{p_0}$-operator norm of $T_{\mathrm{M}_{\Omega,\beta}}$, \eqref{eqA1} implies
\begin{equation} \label{eqA3}
\begin{array}{c} \big\| T_{\mathrm{M}_{\Omega,\beta}^{p,q}}: L_{p_0}(\R) \to L_{p_0}(\R) \big\| \, \le \, 2 \, \mathrm{C}_0 \, \mathrm{HT}_{p_0} \\ [5pt] \mbox{for} \ \mathrm{M}_{\Omega,\beta}^{p,q} (s) = \mathrm{M}_{\Omega,\beta} \Big( \frac{\sqrt{p^2+q^2}}{\mathrm{L}_0 \sqrt{1 + \beta^2}} s \Big) \, \ind_{[0, \mathrm{L}_0]}(s).  \end{array}
\end{equation}
Divide the segment in $\gamma$ running from the origin to the point $(\mathrm{L}_0, \beta \mathrm{L}_0)$ into $pq$ equidistributed points. Formally, we identify this segment with the torus $\mathbb{T}$ and the set of points $$\Big\{ \Big( \frac{\mathrm{L}_0 k}{pq}, \beta \frac{\mathrm{L}_0 k}{pq} \Big) \; : \; 0 \le k \le pq - 1 \Big\} \simeq \big\{ k/pq \; : \; 0 \le k \le pq-1 \big\}$$ with the cyclic group $\Z_{pq}$. According to \eqref{eqA3} and de Leeuw's restriction 
\renewcommand{\theequation}{A.4}
\addtocounter{equation}{-1}
$$\begin{array}{c} \big\| T_{(\mathrm{M}_{\Omega,\beta}^{p,q})_{\mid \Z_{pq}}}: L_{p_0}(\widehat{\Z_{pq}}) \to L_{p_0}(\widehat{\Z_{pq}}) \big\| \, \le \, 2 \, \mathrm{C}_0 \, \mathrm{HT}_{p_0}. \end{array}$$ Since $p$ and $q$ are coprime, we may consider the group isomorphism $$\Lambda: \Z_{pq} \ni \frac{k}{pq} \, \mapsto \, \Big(\frac{k}{p},\frac{k}{q}\Big) \in \Z_p \times \Z_q,$$ where $\Z_p \times \Z_q$ is viewed as a lattice of $[0,1]^2$. It is clear that $\Lambda$ extends to an isometry on $L_{p_0}$ of the corresponding dual groups (still denoted by $\Lambda$) and we obtain  
\begin{equation} \label{eqA4}
\big\| T_{m_{p,q}}: L_{p_0}(\widehat{\Z_{p} \times \Z_q}) \to L_{p_0}(\widehat{\Z_{p} \times \Z_q}) \big\| \, \le \, 2 \, \mathrm{C}_0 \, \mathrm{HT}_{p_0}, 
\end{equation}
where $m_{p,q}(s_1,s_2) = \mathrm{M}_{\Omega,\beta}^{p,q} (\mathrm{L}_0 \Lambda^{-1}(s_1,s_2))$. Given $0\leq k_1 \leq p-1$ and $0\leq k_2 \leq q-1$ let $k = k(k_1,k_2)$ be the only integer $0 \le k \leq pq-1$ satisfying that $k \, \mathrm{mod} \, p \, = \, k_1$ and $k \, \mathrm{mod} \, q \, = \, k_2$. Then, we can write 
\begin{eqnarray*}
m_{p,q} \Big( \frac{k_1}{p} , \frac{k_2}{q} \Big) & = & \mathrm{M}_{\Omega,\beta}^{p,q} \Big(\frac{\mathrm{L}_0 k}{pq} \Big) \ = \ \mathrm{M}_{\Omega,\beta}  \Big( \frac{\sqrt{p^2+q^2}}{pq \sqrt{1 + \beta^2}} k \Big) \, \ind_{[0, \mathrm{L}_0]} \Big( \frac{\mathrm{L}_0 k}{pq} \Big) \\ & = & \mathrm{M}_{\Omega,\beta}(\mathrm{L}_{p,q})  \Big( \frac{\sqrt{p^2+q^2}}{pq \sqrt{1 + \beta^2}} k \Big) \quad \mbox{with} \quad \mathrm{L}_{p,q} =  \frac{\sqrt{p^2 + q^2}}{\sqrt{1 + \beta^2}}.
\end{eqnarray*}
Letting $e_\beta$ and $e_{\frac{p}{q}}$ be the unit vectors in the directions of $\gamma$ and $(\frac1p,\frac1q)$ respectively 
\begin{eqnarray*}
\Big| \frac{\sqrt{p^2+q^2}}{pq} k \frac{(1,\beta)}{\sqrt{1+\beta^2}} - \Big( \frac{k}{p}, \frac{k}{q} \Big) \Big| & = & \Big| \frac{\sqrt{p^2+q^2}}{pq} k e_\beta - \frac{\sqrt{p^2+q^2}}{pq} k e_{\frac{p}{q}} \Big| \\ & \le & \sqrt{p^2 + q^2} \, \big| e_\beta - e_{\frac{p}{q}} \big| \, \lesssim \, \frac{\sqrt{p^2 + q^2}}{q^2} \, \lesssim \, \frac{1}{q}
\end{eqnarray*}
since we may assume with no loss of generality that $\beta < 1$ and $p < q$. We obtain 
\begin{eqnarray*}
m_{p,q} \Big( \frac{k}{p},\frac{k}{q} \Big) & = & \ind_{\Omega}\Big(\Big( \frac{k}{p}, \frac{k}{q} \Big)+\alpha(k)+\Z^2\Big) \\ & = & \ind_{\Pi_\Omega(\mathrm{L}_{p,q},\beta)}\Big(\Big( \frac{k}{p}, \frac{k}{q} \Big)+\alpha(k)+\Z^2\Big) \quad \mbox{with} \quad |\alpha(k)| \lesssim \frac{1}{q}.
\end{eqnarray*}
We deduce that there must exist a small perturbation $\Omega(p,q)$ of $\Omega$ so that $$m_{p,q}\, = \, \ind_{\Omega(p,q)_{\mid_{\Z_p\times \Z_q}}} \qquad \mbox{and} \qquad \Omega(p,q) \to \Omega \ \mbox{ uniformly as } \ p,q\to \infty.$$ By considering the symbol $\widetilde{m}_{p,q} = \ind_{\Omega(p,q)}: \mathbb{T}^2 \to \C$ for $(p,q)\in \mathcal{I}$, we get a sequence of symbols which converges uniformly to $\ind_\Omega$ and satisfy the uniform estimate below $$\sup_{(p,q)\in \mathcal{I}}\big\|T_{(\widetilde{m}_{p,q})_{\mid_{\Z_p \times \Z_q}}}:L_{p_0}(\widehat{\Z_p\times \Z_q})\to L_{p_0}(\widehat{\Z_p\times \Z_q})\big\|\leq \, 2 \, \mathrm{C}_0 \, \mathrm{HT}_{p_0}.$$ By the lattice approximation result obtained in Remark  \ref{Rk-Igari}, this would imply that $\ind_\Omega$ yields an $L_{p_0}$-bounded Fourier multiplier in $\mathbb{T}^2$, and also in $\R^2$ by standard periodization and Hilbert transform truncation. This is a contradiction since $\Omega$ admits Kakeya sets of directions. The proof is complete. \fin

\begin{Aremark}
{\rm This result also holds in higher-dimensions by using cocycles into higher-dimensional spaces, essentially the same argument applies. On the other hand, when $\beta \in \R \setminus \mathbb{Q}$ and $\Omega$ is a polytope with infinitely many faces not admitting Kakeya sets of directions, the conjecture is that such $\Omega$ should define a $L_p$-bounded Fourier multiplier $(1 < p < \infty)$ so that we may argue as we did for polyhedrons with finitely many faces. In dimension 2 this is supported by the results in \cite{Bateman,CF} as for higher-dimensions by \cite{PR}.}
\end{Aremark}

\section*{{B. \hskip4pt \bf Noncommutative Jodeit theorems}}

Jodeit's theorem \cite{Jo} provides another approach to de Leeuw's compactification by looking at extensions of Fourier multipliers. He proved that any $L_p$-bounded Fourier multiplier on $\Z^n$ is the restriction of a $L_p$-bounded Fourier multiplier on $\R^n$. To be more precise, define $$M_p^q(\G) \, = \, \Big\{ m: \G \to \C \ \big| \, T_m: L_p(\widehat{\G}) \to L_q(\widehat{\G}) \Big\}$$ for $1\leq p\leq q\leq \infty$. One of the results in \cite{Jo} is that there is a bounded linear map $\phi: M_p^q(\Z^n)\to M_p^q(\R^n)$ so that the restriction of $\phi(m)$ to $\Z^n$ is $m$. When $n=1$, the symbol $\phi(m) = \widetilde{m}$ is just the multiplier given by the piecewise linear extension $\widetilde m= \ind_{[-\frac 12, \frac 12]} * m*\ind_{[-\frac12, \frac 12]}$ of $m$. Then the ADS property readily gives compactification but one looses on the norm by some  constant depending on $n$.
 
This question of extending multipliers from a subgroup makes sense for general LCA groups and suits in our framework. A commutative solution was provided by Fig\`a-Talamanca and Gaudry in \cite{FTG} by extending Jodeit's result to arbitrary discrete subgroups $\Gamma$ of LCA groups $\G$. Given any such pair, they construct a contractive map $\phi: M_p^q(\Gamma) \to M_p^q(\G)$ so that $\phi(m)=\widetilde{m}$ with $\widetilde{m} = \Delta * m * \Delta$ where $\Delta$ is a positive definite function with small support relative to $\Gamma$. This is not the exact analogue of Jodeit's result (as $\Delta=\ind_{[-\frac 12, \frac 12]}*\ind_{[-\frac 12, \frac 12]}$) but one gains on the constants. Shortly after, Cowling \cite{Cowling} generalized it to all pairs $\H \subset \G$ where $\H$ is closed but not open, $\G$ LCA and $m \in \mathcal{C}_c(\H)$. In the same paper, he also looked at periodization. The underlying idea is to use suitably the disintegration theory and with that respect are of commutative nature.

If we restrict ourselves only to discrete subgroups, such a result would perfectly fit in our framework. In full generality, we yet do not have the right tools to extend Fourier multipliers. However, for the completely bounded ones, we can use transference from Schur multipliers as in Section \ref{cb}. Indeed, the latter are much more flexible and it is proved in \cite[Lemma 2.6]{LdlS} that a Jodeit's theorem  for them is elementary. More precisely, if $\Gamma \subset \G$ is a lattice with a symmetric fundamental domain $\mathrm{X}$ and $m: \Gamma \to \C$ is a cb-bounded Schur multipliers on $S_p(\ell_2(\Gamma))$, then $\widetilde{m} = \ind_\mathrm{X} * m * \ind_\mathrm{X}$ is a cb-bounded Schur multiplier on $S_p(L_2(\G))$. In particular we obtain the following extension result.

\begin{Btheorem}
Let $\Gamma \subset \G$ be a lattice in an amenable locally compact group $\G$ with a symmetric fundamental domain $\mathrm{X}$. For any $m: \Gamma \to \C$ with $\widetilde{m} = \ind_\mathrm{X} * m * \ind_\mathrm{X}$ $$ \big\| T_{\widetilde{m}} : L_p(\widehat{\G}) \to L_p(\widehat{\G}) \big\|_{\rm cb} \, \le \, \big\| T_{m}: L_p(\widehat{\Gamma}) \to L_p(\widehat{\Gamma}) \big\|_{\rm cb}.$$
\end{Btheorem}

\noindent In particular, the cb-bounded version of Jodeit's theorem holds with constant 1.

\vskip3pt

\noindent \textbf{Acknowledgement.} M. Caspers is partially supported by the grant SFB 878 \lq\lq Groups, geometry and actions"; J. Parcet and M. Perrin are partially supported by the ERC StG-256997-CZOSQP (UE) and ICMAT Severo Ochoa SEV-2011-0087 (Spain); and \'E. Ricard by ANR-2011-BS01-008-11 (France).

\bibliography{biblio}
\bibliographystyle{amsplain}

\null

\hfill \noindent \textbf{Martijn Caspers} \\
\null \hfill Fachbereich Mathematik und Informatik
\\ \null \hfill Westf\"alische Wilhelmsuniversit\"at  M\"unster \\
\null \hfill Einsteinstrasse 62, 48149 M\"unster, Germany \\
\null \hfill\texttt{martijn.caspers@uni-muenster.de}

\vskip2pt

\hfill \noindent \textbf{Javier Parcet} \\
\null \hfill Instituto de Ciencias Matem{\'a}ticas \\ \null \hfill
CSIC-UAM-UC3M-UCM \\ \null \hfill Consejo Superior de
Investigaciones Cient{\'\i}ficas \\ \null \hfill C/ Nicol\'as Cabrera 13-15.
28049, Madrid. Spain \\ \null \hfill\texttt{javier.parcet@icmat.es}

\vskip2pt

\hfill \noindent \textbf{Mathilde Perrin} \\
\null \hfill Instituto de Ciencias Matem{\'a}ticas \\ \null \hfill
CSIC-UAM-UC3M-UCM \\ \null \hfill Consejo Superior de
Investigaciones Cient{\'\i}ficas \\ \null \hfill C/ Nicol\'as Cabrera 13-15.
28049, Madrid. Spain \\ \null \hfill\texttt{mathilde.perrin@icmat.es}

\vskip2pt

\hfill \noindent \textbf{\'Eric Ricard} \\
\null \hfill Laboratoire de Math\'ematiques Nicolas Oresme \\
\null \hfill Universit\'e de Caen Basse-Normandie \\ 
\null \hfill 14032 Caen Cedex, France \\
\null \hfill \texttt{eric.ricard@unicaen.fr}
\end{document}